\documentclass{article} 
\usepackage{amssymb}
\usepackage{amsmath}
\usepackage{amsthm}
\usepackage{geometry}
\usepackage[numbers,sort&compress]{natbib}
\usepackage[colorlinks=true, linkcolor=blue, citecolor=blue]{hyperref}

\usepackage{amssymb}
\usepackage{amsmath,amscd}
\usepackage{amsthm} 
\usepackage{cancel}     
\usepackage{centernot}  
\usepackage{graphicx}   
\usepackage{soul}       
\usepackage{color}      

\newtheorem{thm}{Theorem}[section]
\newtheorem{lemma}[thm]{Lemma}
\newtheorem{cor}[thm]{Corollary}
\newtheorem{prop}[thm]{Proposition}

\newtheorem*{thm*}{Theorem}
\newtheorem*{lemma*}{Lemma}
\newtheorem*{cor*}{Corollary}
\newtheorem*{prop*}{Proposition}
\newtheorem*{conjecture*}{Conjecture}

\theoremstyle{definition}
\newtheorem{defn}{Definition}[section]
\newtheorem*{defn*}{Definition}
\theoremstyle{definition}

\theoremstyle{definition}
\newtheorem{ex}{Example}
\theoremstyle{remark}
\newtheorem*{ex*}{Example}
\theoremstyle{definition}

\theoremstyle{definition}
\newtheorem*{assm*}{Assumption}
\theoremstyle{remark}
\newtheorem{remark}{Remark}[section]
\theoremstyle{remark}
\newtheorem*{remark*}{Remark}

\DeclareFontFamily{U}{mathx}{\hyphenchar\font45}
\DeclareFontShape{U}{mathx}{m}{n}{
      <5> <6> <7> <8> <9> <10> gen * mathx
      <10.95> mathx10 <12> <14.4> <17.28> <20.74> <24.88> mathx12
      }{}
\DeclareSymbolFont{mathx}{U}{mathx}{m}{n}
\DeclareMathSymbol{\intop}  {1}{mathx}{"B3}

\newcommand{\wt}{\widetilde}
\newcommand{\wh}{\widehat}
\newcommand{\cl}{\overline}

\let\temp\phi
\let\phi\varphi
\let\varphi\temp

\newcommand{\pr}{\mathbb{P}}

\newcommand{\R}{\mathbb{R}}

\newcommand{\E}{\mathbb{E}}

\newcommand{\normalN}{\mathcal{N}}
\newcommand{\given}{\,|\,}  

\newcommand{\norm}[1]{\Vert#1\Vert}

\newcommand{\ip}[1]{\langle#1\rangle}

\newcommand{\eps}{\varepsilon}

\newcommand{\iid}{\overset{\text{iid}}{\sim}}
\newcommand{\symdiff}{\triangle}

\DeclareMathOperator{\tr}{tr}
\DeclareMathOperator{\supp}{supp}

\DeclareMathOperator{\conv}{conv}

\DeclareMathOperator*{\argmax}{arg\,max}

\DeclareMathOperator{\GammaDist}{Gamma}

\DeclareMathOperator{\GumbelDist}{Gumbel}

\usepackage{enumitem}
\usepackage{pgfplots}
\usetikzlibrary{patterns}
\usepackage{subcaption}
\usepackage{cleveref}
\usepackage{booktabs}   
\usepackage{commath}

\DeclareMathOperator{\PD}{PD}

\newcommand{\probs}{\mathcal{P}}
\newcommand{\coupling}{\sigma}
\newcommand{\mixing}{\probs^{2}}
\newcommand{\gen}{\mathcal{K}}
\newcommand{\mixgen}{\mathfrak{K}}
\newcommand{\mix}{\mathcal{M}}
\newcommand{\mixfcn}{m}
\newcommand{\mixset}{\mathfrak{L}}

\newcommand{\gaussians}{\mathcal{G}}
\newcommand{\mixgaussians}{\mathfrak{G}}
\newcommand{\lebesgue}{\zeta}

\newcommand{\Hdiam}{\Delta}
\newcommand{\separ}{\eta}
\newcommand{\gensep}{\overline{\separ}}
\newcommand{\sepconst}{\xi}
\newcommand{\normalmean}{a}
\newcommand{\normalvar}{v}
\newcommand{\normalmeanset}{A}
\newcommand{\normalvarset}{V}

\newcommand{\true}{\Gamma}
\newcommand{\truecmp}{\gamma}
\newcommand{\truemix}{\Lambda}
\newcommand{\truewgt}{\lambda}
\newcommand{\truedens}{F}
\newcommand{\truedenscmp}{f}

\newcommand{\genproj}{Q}
\newcommand{\gencmp}{q}
\newcommand{\genmix}{\Omega}
\newcommand{\genwgt}{\omega}
\newcommand{\gendist}{D}

\newcommand{\projset}{\mathfrak{\genproj}}
\newcommand{\proj}{\genproj^{*}}
\newcommand{\projcmp}{\gencmp^{*}}
\newcommand{\projmix}{\genmix^{*}}
\newcommand{\projwgt}{\genwgt^{*}}
\newcommand{\projdist}{\gendist^{*}}

\newcommand{\genclustwgt}{\varpi} %
\newcommand{\projclustwgt}{\genclustwgt^{*}}
\newcommand{\estclustwgt}{\wh{\genclustwgt}}

\newcommand{\gendens}{G}
\newcommand{\gendenscmp}{g}
\newcommand{\projdens}{\gendens^{*}}
\newcommand{\projdenscmp}{\gendenscmp^{*}}

\newcommand{\estproj}{\wh{\genproj}}
\newcommand{\estcmp}{\wh{\gencmp}^{\vphantom{*}}}
\newcommand{\estmix}{\wh{\genmix}^{\vphantom{*}}}
\newcommand{\estwgt}{\wh{\genwgt}^{\vphantom{*}}}
\newcommand{\estdens}{\wh{\gendens}}
\newcommand{\estdenscmp}{\wh{\gendenscmp}}
\newcommand{\estdist}{\wh{\gendist}}

\newcommand{\projmap}{T}
\newcommand{\mixingmap}{M}
\newcommand{\clustmap}{\chi}

\newcommand{\assg}{\alpha}
\newcommand{\assgseq}{\{\assg_{L}\}}
\newcommand{\assignsym}{\mathbb{A}}
\newcommand{\assignmaps}{\assignsym_{L\to K}}
\newcommand{\assignseq}{\assignsym_{K}^{\infty}}
\newcommand{\assignproj}{\projset_{L\to K}}
\newcommand{\estassg}{\wh{\assg}}

\newcommand{\rv}{Z}
\newcommand{\base}{X}
\newcommand{\basem}{d}
\newcommand{\probm}{\rho}
\newcommand{\mixm}{W_{r}} %

\newcommand{\eqset}{E_{0}}
\newcommand{\partn}{\Pi}
\newcommand{\pt}{\pi}

\newcommand{\tick}[1]{#1'}

\newcommand{\regassg}{\assignseq}
\newcommand{\tseq}{t_{L,n}}

\newcommand{\gaussmix}{P}
\newcommand{\dH}{\probm_{\textup{H}}}
\newcommand{\dTV}{\probm_{\textup{TV}}}

\newcommand{\mixfamily}{\mathfrak{F}}

\newcommand{\disjmix}{\mixfamily_{1}}
\newcommand{\metamix}{\mixfamily_{2}}

\newcommand{\projmean}{\normalmean^{*}}
\newcommand{\projvar}{\normalvar^{*}}

\newcommand{\wconv}{\Rightarrow}
\newcommand{\minsep}{\probm_{*}}
\newcommand{\bb}{4}

\newcommand{\intlo}{b}
\newcommand{\intup}{c}

\newcommand{\dx}{\,d\lebesgue}

\title{Identifiability of Nonparametric Mixture Models and Bayes Optimal Clustering}
\author{Bryon Aragam \and Chen Dan \and Eric P. Xing \and Pradeep Ravikumar}

\begin{document}
\maketitle

{\let\thefootnote\relax\footnote{Contact: \texttt{bryon@chicagobooth.edu, cdan@andrew.cmu.edu, epxing@cs.cmu.edu, pradeepr@cs.cmu.edu}}}

\begin{abstract}
Motivated by problems in data clustering, we establish general conditions under which families of nonparametric mixture models are identifiable, by introducing a novel framework involving clustering overfitted \emph{parametric} (i.e. misspecified) mixture models. These identifiability conditions generalize existing conditions in the literature, and are flexible enough to include for example mixtures of Gaussian mixtures. In contrast to the recent literature on estimating nonparametric mixtures, we allow for general nonparametric mixture components, and instead impose regularity assumptions on the underlying mixing measure. As our primary application, we apply these results to partition-based clustering, generalizing the notion of a Bayes optimal partition from classical parametric model-based clustering to nonparametric settings. Furthermore, this framework is constructive so that it yields a practical algorithm for learning identified mixtures, which is illustrated through several examples on real data. The key conceptual device in the analysis is the convex, metric geometry of probability measures on metric spaces and its connection to the Wasserstein convergence of mixing measures. The result is a flexible framework for nonparametric clustering with formal consistency guarantees.
\end{abstract}

\section{Introduction}
\label{sec:intro}

In data clustering, we provide a grouping of a set of data points, or more generally, a partition of the input space from which the data points are drawn~\citep{hartigan1975}. The many approaches to formalize the learning of such a partition from data include mode clustering \citep{chen2016}, density clustering \citep{rinaldo2010,sriperumbudur2012,steinwart2011,steinwart2015}, spectral clustering \citep{ng2001,yan2009,schiebinger2015}, $K$-means \citep{macqueen1967,steinhaus1956,lloyd1982}, stochastic blockmodels \citep{holland1983,rohe2011,airoldi2013,eldridge2016}, and hierarchical clustering \citep{hartigan1981,chaudhuri2010,thomann2015}, among others. In this paper, we are interested in so-called model-based clustering where the data points are drawn i.i.d. from some distribution, the most canonical instance of which is arguably Gaussian model-based clustering, in which points are drawn from a Gaussian mixture model \citep{dasgupta1999,arora2005}. This mixture model can then be used to specify a natural partition over the input space, specifically into regions where each of the Gaussian mixture components is most likely. When the Gaussian mixture model is appropriate, this provides a simple, well-defined partition, and has been extended to various parametric and semi-parametric models \citep{wolfe1970,bock1996,fraley2002}. However, the extension of this methodology to general nonparametric settings has remained elusive. This is largely due to the extreme non-identifiability of nonparametric mixture models, a problem which is well-studied but for which existing results require strong assumptions~\citep{teicher1967,bordes2006,hunter2007,jochmans2017}. It has been a significant open problem to generalize these assumptions to a more flexible class of nonparametric mixture models. 

Unfortunately, without the identifiability of the mixture components, we cannot extend the notion of the input space partition used in Gaussian mixture model clustering.
Nonetheless, there are many practical clustering algorithms used in practice, such as $K$-means and spectral techniques, that do estimate a partition even when the data arises from ostensibly unidentifiable nonparametric mixture models, such as mixtures of sub-Gaussian or log-concave distributions \cite{achlioptas2005,kannan2008,shi2009,mixon2017}.
A crucial motivation for this paper is in addressing this gap between theory and practice: This entails demonstrating that nonparametric mixture models might actually be identifiable given additional \emph{side information}, such as the number of clusters $K$ and the separation between the mixture components, used for instance by algorithms such as $K$-means.

Let us set the stage for this problem in some generality. Suppose $\true$ is a probability measure over some metric space $\base$, and that $\true$ can be written as a finite mixture model 
\begin{align}
\label{eq:true:mixmodel}
\true
= \sum_{k=1}^{K}\truewgt_{k}\truecmp_{k},
\quad
\truewgt_{k}>0
\;\;\text{and}\;\;
\sum_{k=1}^{K}\truewgt_{k}=1,
\end{align}

\noindent
where $\truecmp_{k}$ are also probability measures over $\base$. The $\truecmp_{k}$ represent distinct subpopulations belonging to the overall heterogeneous population $\true$. Given observations from $\true$, we are interested in classifying each observation into one of these $K$ subpopulations \emph{without} labels. When the mixture components $\truecmp_{k}$ and their weights $\truewgt_{k}$ are identifiable, we can expect to learn the model \eqref{eq:true:mixmodel} from this unlabeled data, and then obtain a partition of $\base$ into regions where one of the mixture components is most likely. This can also be cast as using Bayes' rule to classify each observation, thus defining a target partition that we call the \emph{Bayes optimal partition} (see Section~\ref{sec:bayesopt} for formal details). 
Thus, in studying these partitions, a key question is \emph{when is the mixture model \eqref{eq:true:mixmodel} identifiable?} Motivated by the aforementioned applications to clustering, this question is the focus of this paper. Under parametric assumptions such as Gaussianity of the $\truecmp_{k}$, it is well-known that the representation \eqref{eq:true:mixmodel} is unique and hence identifiable \citep{teicher1963,barndorff1965,holzmann2006}. These results mostly follow from an early line of work on the general identification problem \citep{teicher1961,teicher1963,yakowitz1968,ahmad1982}. 

Such parametric assumptions rarely hold in practice, however, and thus it is of interest to study \emph{nonparametric} mixture models of the form \eqref{eq:true:mixmodel}, i.e. for which each $\truecmp_{k}$ comes from a flexible, nonparametric family of probability measures.
In the literature on nonparametric mixture models, a common assumption is that the component measures $\truecmp_{k}$ are multivariate with independent marginals \citep{teicher1967,hall2003,hall2005mixture,levine2011,dhaultfuille2015}, which is particularly useful for statistical problems involving repeated measurements \citep{hettmansperger2000,bonhomme2016repeated}. This model also has deep connections to the algebraic properties of latent structure models \citep{allman2009,bonhomme2016latent}. Various other structural assumptions have been considered including symmetry \citep{bordes2006,hunter2007}, tail conditions \citep{jochmans2017}, 
and translation invariance \citep{gassiat2013}. 
The identification problem in discrete mixture models is also a central problem in topic models which are popular in machine learning \citep{anandkumar2013,tang2014,anandkumar2015}. 
Most notably, this existing literature imposes structural assumptions on the components $\truecmp_{k}$ (e.g. independence, symmetry), which are difficult to satisfy in clustering problems. Are there reasonable constraints that ensure the uniqueness of \eqref{eq:true:mixmodel}, while avoiding restrictive assumptions on the $\truecmp_{k}$?

In this paper, we establish a series of positive results in this direction, and as a bonus that arises naturally from our theoretical results, we develop a practical algorithm for nonparametric clustering. In contrast to the existing literature, we allow each $\truecmp_{k}$ to be an arbitrary probability measure over $\base$. We propose a novel framework for reconstructing nonparametric mixing measures by using simple, overfitted mixtures (e.g. Gaussian mixtures) as mixture density estimators, and then using clustering algorithms to partition the resulting estimators. This construction implies a set of regularity conditions on the mixing measure that suffice to ensure that a mixture model is identifiable. As our main application of interest, we apply this to problems in nonparametric clustering.

In the remainder of this section, we outline our major contributions and present a high-level geometric overview of our method.
Section~\ref{sec:prelim} covers the necessary background required for our abstract framework. 
In Section~\ref{sec:ident}, we introduce two important concepts, \emph{regularity} and \emph{clusterability}, that are crucial to our identifiability results, along with several examples.
In Section~\ref{sec:estimation} we show how these concepts are sufficient to ensure identifiability of a nonparametric mixture model and consistency of a minimum distance estimator.
In Section~\ref{sec:bayesopt} we apply these results to the problem of clustering and prove a consistency theorem for this problem. Section~\ref{sec:exp} introduces a simple algorithm for nonparametric clustering along with some experiments, and Section~\ref{sec:conc} concludes the paper with some discussion and extensions. All proofs are deferred to the appendices.

\paragraph{Contributions} 
\label{sec:intro:contrib}

Our main results can be divided into three main theorems:
\begin{enumerate}[label=\arabic*.,leftmargin=0.5cm,itemindent=.5cm, itemsep=.5em]
\item \emph{Nonparametric identifiability} (Section~\ref{sec:ident:main}). We formulate a general set of assumptions that guarantee a family of nonparametric mixtures will be identifiable (Theorem~\ref{thm:main:ident}), based on two properties introduced in Section~\ref{sec:ident}: \emph{regularity} (Definition~\ref{defn:regular}) and \emph{clusterability} (Definition~\ref{defn:clusterable}).
\item \emph{Estimation and consistency} (Section~\ref{sec:estimation:est}). We show that a simple clustering procedure will correctly identify the mixing measure that generates $\true$ as long as the $\truecmp_{k}$ are sufficiently well-separated, and this procedure defines an estimator that consistently recovers the nonparametric clusters given i.i.d. observations from $\true$ (Theorem~\ref{thm:main:learn}). 
\item \emph{Application to nonparametric clustering} (Section~\ref{sec:bayesopt}). We extend the notion of a \emph{Bayes optimal partition} (Definition~\ref{defn:bayesopt}) to general nonparametric settings and prove a consistency theorem for recovering such partitions when they are identified (Theorem~\ref{thm:main:part}). 
\end{enumerate}

\noindent
Each of these contributions builds on the previous one, and provides an overall narrative that strengthens the well-known connections between identifiability in mixture models, cluster analysis, and nonparametric density estimation. 
We conclude our study by applying these results to construct an intuitive algorithm for nonparametric clustering, which is investigated in Section~\ref{sec:exp}.

\paragraph{Overview}
\label{sec:overview}

\begin{figure}[t]
\centering
\begin{subfigure}[t]{0.5\textwidth}
\centering
\includegraphics[width=\textwidth]{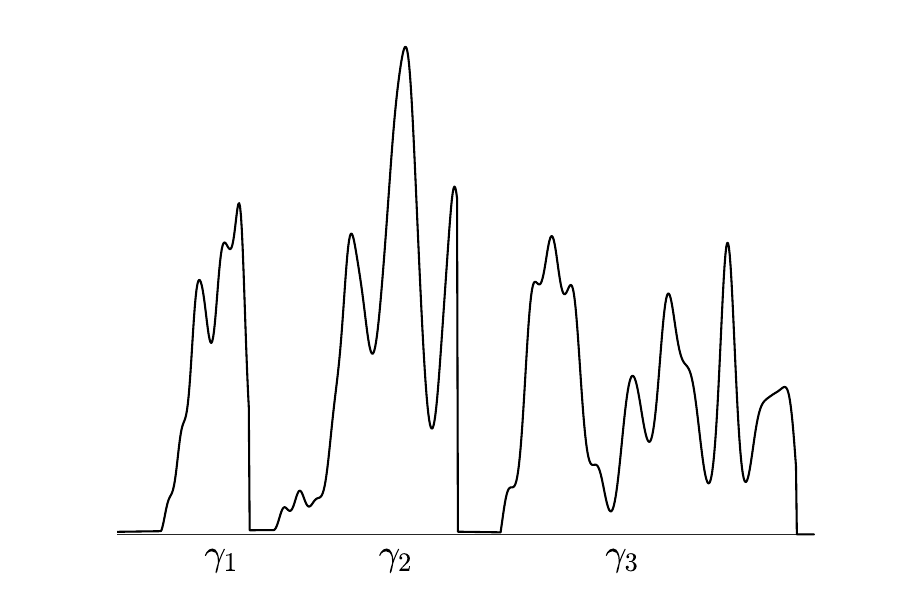}
\caption{Original mixture $\true=\mixfcn(\truemix)=\sum_{k}\truewgt_{k}\truecmp_{k}$ and $K=3$.}
\label{fig:overview:true}
\end{subfigure}%
~
\begin{subfigure}[t]{0.5\textwidth}
\centering
\includegraphics[width=\textwidth]{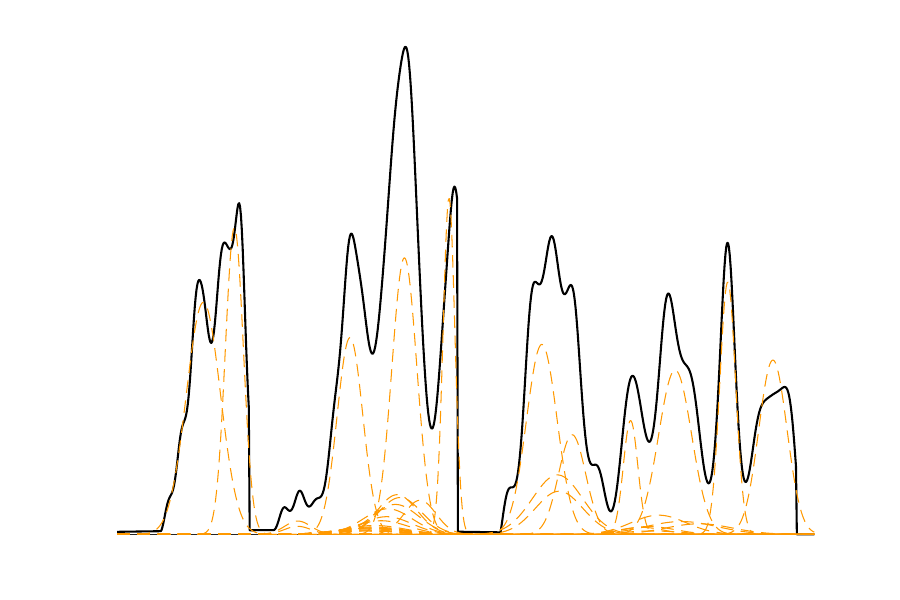}
\caption{Approximate mixture of Gaussians $\proj=\mixfcn(\projmix)=\sum_{\ell}\projwgt_{\ell}\projcmp_{\ell}$.}
\label{fig:overview:proj}
\end{subfigure}%
\\
\begin{subfigure}[t]{0.5\textwidth}
\centering
\includegraphics[width=\textwidth]{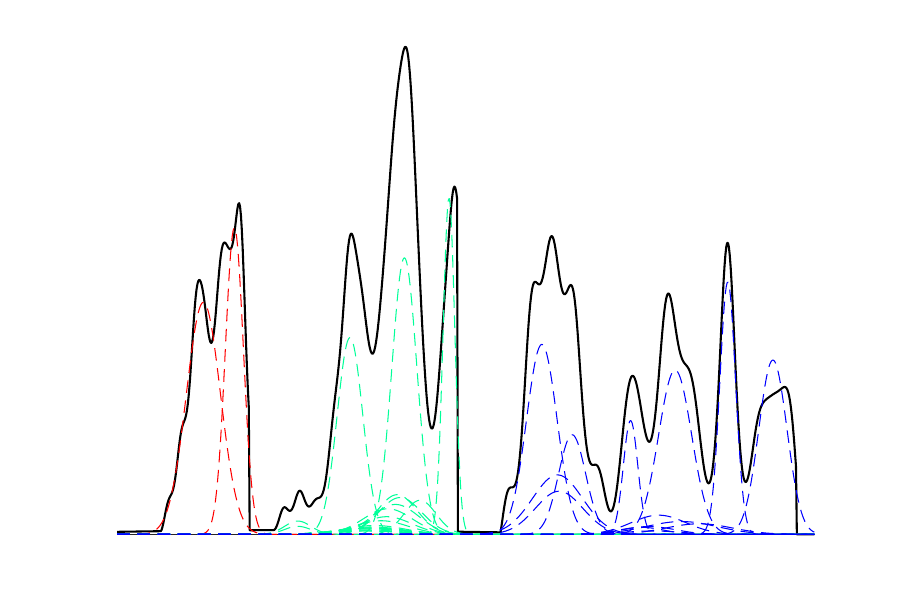}
\caption{Components $\projcmp_{\ell}$ grouped by clustering.}
\label{fig:overview:clust}
\end{subfigure}%
~
\begin{subfigure}[t]{0.5\textwidth}
\centering
\includegraphics[width=\textwidth]{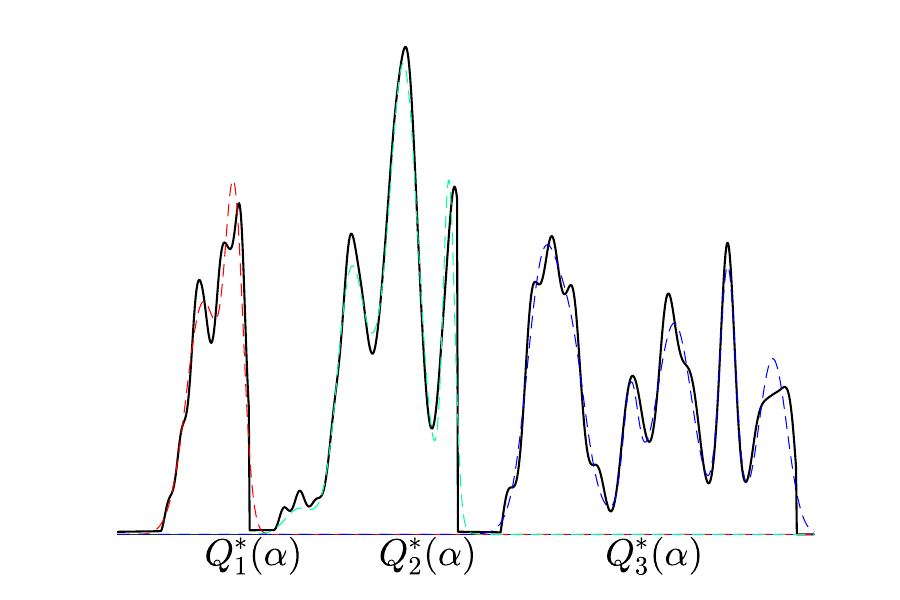}
\caption{Final approximate nonparametric mixing measure $\genmix(\assg)$.}
\label{fig:overview:mix}
\end{subfigure}%
\caption{Overview of the method.}
\label{fig:overview}
\end{figure}

Before outlining the formal details, we present an intuitive geometric picture of our approach in Figure~\ref{fig:overview}. This same example is developed in more detail in the experiments (see Section~\ref{sec:exp}, Figure~\ref{fig:1d:sobolev}).
At a high-level, our strategy for identifying the mixture distribution \eqref{eq:true:mixmodel} is the following:
\begin{enumerate}[label=(\arabic*), itemsep=0mm]
\item\label{steps:overfit} Approximate $\true$ with an overfitted mixture of $L\gg K$ Gaussians (Figure~\ref{fig:overview:proj});
\item\label{steps:cluster} Cluster these $L$ Gaussian components into $K$ groups such that each group roughly approximates some $\truecmp_{k}$ (Figure~\ref{fig:overview:clust});
\item\label{steps:mixing} Use this clustering to define a new mixing measure (Figure~\ref{fig:overview:mix});
\item\label{steps:conv} Show that this new mixing measure converges to the true mixing measure $\truemix$ as $L\to\infty$.
\end{enumerate}

\noindent
Of course, this construction is not guaranteed to succeed for arbitrary mixing measures $\truemix$, which will be illustrated by the examples in Section~\ref{sec:prelim:ident}. This is a surprisingly subtle problem and requires careful consideration of the various spaces involved. Thus, a key aspect of our analysis will be to provide assumptions that ensure the success of this construction. Intuitively, it should be clear that as long as the $\truecmp_{k}$ are well-separated, the corresponding mixture approximation will consist of Gaussian components that are also well-separated. Unfortunately, this is not quite enough to imply identifiability, as illustrated by Example~\ref{ex:danger:mix3}. This highlights some of the subtleties inherent in this construction. 
Furthermore, although we have used mixtures of Gaussians to approximate $\true$ in this example, our main results will apply to any properly chosen family of base measures.

\section{Preliminaries}
\label{sec:prelim}

Our approach is general, built on the theory of abstract measures on metric spaces \citep{parthasarathy1967}. In this section we introduce this abstract setting, outline our notation, and discuss the general problem of identifiability in mixture models. 
For a more thorough introduction to the general topic of mixture models in statistics, see \citet{lindsay1995,ritter2014,titterington1985}.

\subsection{Nonparametric mixture models}
\label{sec:prelim:npmix}

Let $(\base, \basem)$ be a metric space and $(\probs(\base),\probm)$ denote the space of regular Borel probability measures on $\base$ with finite $r$th moments ($r\ge 1$) metrized by a metric $\probm$. Common choices for $\probm$ include the Hellinger and variational metrics, however, our results will apply to any metric on $\probs(\base)$.
Define $\mixing(\base)=\probs(\probs(\base))$, the space of (infinite) mixing measures over $\probs(\base)$. In this paper, we study \emph{finite} mixture models, i.e. mixtures with a finite number of atoms. To this end, define for $s\in\{1,2,\ldots\}$
\begin{align*}
\mixing_{s}(\base)
&:= \{\truemix\in\mixing(\base) : |\supp(\truemix)|\le s\},
\qquad
\mixing_{0}(\base)
:=\bigcup_{s=1}^{\infty}\mixing_{s}(\base).
\end{align*}

\noindent
We treat $\mixing(\base)$ as a metric space by endowing it with the $L_{r}$-Wasserstein metric $\mixm$ ($r\ge 1$). 
When $\truemix\in\mixing_{K}(\base)$ and $\truemix'\in\mixing_{K'}(\base)$, 
this is given by the optimal value of the transport problem
\begin{align}
\label{eq:defn:wasserstein}
\begin{aligned}
\mixm(\truemix, \truemix')
= \inf\Bigg\{
\Big[\sum_{i,j}\coupling_{ij}\probm^{r}&(\truecmp_{i},\truecmp_{j}')\Big]^{1/r}
:
0\le \coupling_{ij}\le 1,\, \\
&\sum_{i,j}\coupling_{ij}=1,\,
\sum_{i}\coupling_{ij}=\truewgt_{j}',\,
\sum_{j}\coupling_{ij}=\truewgt_{i}
\Bigg\}.
\end{aligned}
\end{align}

\noindent
where the infimum is taken over all couplings $\coupling$, i.e. probability measures on $\probs(\base)\times\probs(\base)$ with marginals $\truemix$ and $\truemix'$. For more on Wasserstein distances and their importance in mixture models, see \citet{nguyen2013}.

Given $\truemix\in\mixing_{0}(\base)$, define a new probability measure $\mixfcn(\,\cdot\,;\truemix)\in\probs(\base)$ by %
\begin{align}
\label{eq:def:mixmeasure}
\mixfcn(A;\truemix)
= \int\truecmp(A)\,d\truemix(\truecmp)
= \sum_{k=1}^{K}\truewgt_{k}\truecmp_{k}(A),
\quad
K := |\supp(\truemix)|,
\end{align}

\noindent
where $\truecmp_{1},\ldots,\truecmp_{K}$ are the mixture components (i.e. a particular enumeration of $\supp(\truemix)$) and $\truewgt_{1},\ldots,\truewgt_{K}$ are the corresponding weights. 
Formally, for any Borel set $A\subset \base$ we have a function $h_{A}:\probs(\base)\to\R$ defined by $h_{A}(\truecmp)=\truecmp(A)$, and $\mixfcn(A;\truemix)=\int\truecmp(A)\,d\truemix(\truecmp)=\int h_{A}\,d\truemix$. This uniquely defines a measure called a \emph{mixture distribution} over $\base$. In a slight abuse of notation, we will write $\mixfcn(\truemix)$ as shorthand for $\mixfcn(\,\cdot\,;\truemix)$ when there is no confusion between the arguments. An element $\truecmp_{k}$ of $\supp(\truemix)$ is called a \emph{mixture component}. Given a Borel set $\mixset\subset\mixing(\base)$, define in analogy with $\mixing_{s}(\base)$ the subsets of finite mixtures by
\begin{align}
\label{eq:def:finitemixtures}
\mixset_{s}=\mixset\cap\mixing_{s}(\base)
\end{align}

\noindent
and
\begin{align}
\label{eq:def:mixturemodel}
\mix(\mixset)
:= \{\mixfcn(\truemix) : \truemix\in\mixset\},
\end{align}

\noindent
i.e. the family of mixture distributions over $\base$ induced by $\mixset$, which can be regarded as a formal representation of a statistical mixture model.

\begin{remark}
This abstract presentation of mixture models is needed for two reasons: (i) To emphasize that $\truemix$ is the statistical parameter of interest, in contrast to the usual parametrization in terms of atoms and weights; and (ii) To emphasize that our approach works for general measures on metric spaces. This will have benefits in the sequel, albeit at the cost of some extra abstraction here at the onset. 
For the most part, we will work with finite mixtures, i.e. $\mixing_{0}(\base)$, a space which should be contrasted with the more complex space of infinite measures $\mixing(\base)$, although some of the examples and proofs will invoke infinite mixtures.
\end{remark}

\begin{remark}
\label{rem:notation}
As a convention, we will use upper case letters for mixture distributions (e.g. $\true$, $\genproj$) and mixing measures (e.g. $\truemix$, $\genmix$), and lower case letters for mixture components (e.g. $\truecmp_{k}$, $\gencmp_{k}$) and weights (e.g. $\truewgt_{k}$, $\genwgt_{k}$).
\end{remark}

We conclude this subsection with some examples.

\begin{ex}[Parametric mixtures]
\label{ex:param:mix}
Let $\mathcal{Q}=\{q_{\theta}:\theta\in\Theta\}$ be a family of measures parametrized by $\theta$. Then any mixing measure whose support is contained in $\mathcal{Q}$ defines a parametric mixture distribution. For example, let $\mixgaussians\subset\mixing(\R^{p})$ denote the subset of mixing measures whose support is contained in the family of $p$-dimensional Gaussian measures. Then $\mix(\mixgaussians)$ is the family of Gaussian mixtures, and $\mix(\mixgaussians_{0})$ is the family of finite Gaussian mixtures. It is well-known that $\mix(\mixgaussians_{0})$ is identifiable \citep{teicher1961,teicher1963}. Other examples include certain exponential family mixtures \citep{barndorff1965} and translation families \citep{teicher1961} (i.e. $q_{\theta}(A)=\mu(A-\theta)$ for some known measure $\mu\in\probs(\R^{d})$).
\end{ex}

\begin{ex}[Sub-Gaussian mixtures]
\label{ex:np:mix}
\sloppypar{Let $\gen$ be the collection sub-Gaussian measures on $\R$, i.e.}
\begin{align*}
\gen = \{
    \truecmp\in\probs(\R) : \truecmp(\{x:|x|>t\}) \le e^{1-t^{2}/c^{2}}\,\text{for some $c>0$ and all $t>0$}
\},
\end{align*}
and $\mixgen\subset\mixing(\R)$ be the subset of mixing measures whose support is a subset of $\gen$. Then $\mix(\mixgen)$ is the family of sub-Gaussian mixture models and $\mix(\mixgen_{0})$ is the family of finite sub-Gaussian mixtures. This is a nonparametric mixture model, since the base measures $\gen$ do not belong to a parametric family. Extensions to sub-Gaussian measures on $\R^{p}$ are natural.
\end{ex}

\noindent
Our definition of mixtures over subsets of mixing measures---as opposed to over families of component distributions---makes it easy to encode additional constraints, as in the following example.

\begin{ex}[Constrained mixtures]
\label{ex:constrain:mix}
Continuing the previous example, suppose we wish to impose additional constraints on the family of mixture distributions. For example, we might be interested in Gaussian mixtures 
whose means are contained within some set $\normalmeanset\subset\R^{p}$ and whose covariance matrices are contained within another set $\normalvarset\subset\PD(p)$, where $\PD(p)$ is the set of $p\times p$ positive-definite matrices. Define $\gaussians({\normalmeanset,\normalvarset}):= \{\normalN(\normalmean,\normalvar): \normalmean\in\normalmeanset, \normalvar\in\normalvarset\}$ and 
\begin{align}
\label{eq:defn:mixgaussians}
\mixgaussians({\normalmeanset,\normalvarset})
&:= \{\truemix\in\mixing(\base) : \supp(\truemix)\subset \gaussians({\normalmeanset,\normalvarset})\}.
\end{align}

\noindent
Then $\mix(\mixgaussians({\normalmeanset,\normalvarset}))$ is the desired family of mixture models. A special case of interest is $V=\{\normalvar\}$ for some fixed $\normalvar\in\PD(p)$, which we denote by $\mixgaussians({\normalmeanset,\normalvar})$, also known as a convolutional (Gaussian) mixture model.
Finite mixtures from these families are denoted by $\mix(\mixgaussians_{0}({\normalmeanset,\normalvarset}))$ and $\mix(\mixgaussians_{0}({\normalmeanset,\normalvar}))$.
\end{ex}

\begin{ex}[Mixture of regressions]
Suppose $\pr(Y\given Z)=\int\truecmp(Z)\,d\truemix(\truecmp)$ is a mixture model depending on some covariates $Z$. We assume here that $(Z,Y)\in W\times X$ where $(W,\basem_{W})$ and $(\base,\basem_{X})$ are metric spaces. This is a nonparametric extension of the usual mixed linear regression model. To recover the mixed regression model, assume $\truemix$ has at most $K$ atoms and $\truecmp_{k}(Z)\sim\normalN(\ip{\theta_{k},Z},\omega_{k}^{2})$, so that
\begin{align*}
\pr(Y\given Z)
=\int\truecmp(Z)\,d\truemix(\truecmp)
=\sum_{k=1}^{K}\truewgt_{k}\normalN(\ip{\theta_{k},Z},\omega_{k}^{2}).
\end{align*}

\noindent
By further allowing the mixing measure $\truemix=\truemix(Z)$ to depend on the covariates, we obtain the nonparametric generalization of a mixture of experts model \citep{jacobs1991,hunter2012,bordes2013}.
\end{ex}

\subsection{Identifiability in mixture models}
\label{sec:prelim:ident}

A mixture model $\mix(\mixset)$ is identifiable if the map $\mixfcn:\mixset\to\mix(\mixset)$ that sends $\truemix\mapsto\mixfcn(\truemix)$ is injective. For an overview of this problem, see \citet{hunter2007} and \citet{allman2009}. The main purpose of this section is to highlight some of the known subtleties in identifying nonparametric mixture models.

Unsurprisingly, whether or not a specific mixture $\mixfcn(\truemix)$ is identified depends on the choice of $\mixset$. If we allow $\mixset$ to be all of $\mixing(\base)$, then it is easy to see that $\mix(\mixset)$ is not identifiable, and this continues to be true even if the number of components $K$ is known in advance (i.e. $\mixset=\mixing_{K}(\base)$). Indeed, for any partition $\{A_{k}\}_{k=1}^{K}$ of $\base$ and any Borel set $B\subset\base$, we can write
\begin{align}
\label{eq:np:mix:nonident}
\true(B)
= \sum_{k=1}^{K} \underbrace{\true(A_{k})\vphantom{\frac{\true( B \cap A_{k})}{\true(A_{k})}}}_{\wt{\truewgt}_{k}}\cdot\underbrace{\frac{\true( B \cap A_{k})}{\true(A_{k})}}_{\wt{\truecmp}_{k}} %
= \sum_{k=1}^{K}\wt{\truewgt}_{k}\wt{\truecmp}_{k}(B),
\end{align}

\noindent
and thus there cannot be a unique decomposition of the measure $\true$ into the sum \eqref{eq:true:mixmodel}. Although this example allows for arbitrary, pathological decompositions of $\true$ into conditional measures, the following concrete example shows that solving the nonidentifiability issue is more complicated than simply avoiding certain pathological partitions of the input space.

\begin{figure}
\centering
\begin{tikzpicture}
\begin{axis}[
  no markers, domain=0:10, samples=200,
  axis lines*=left, %
  every axis y label/.style={at=(current axis.above origin),anchor=south},
  every axis x label/.style={at=(current axis.right of origin),anchor=west},
  height=3cm, width=5cm,
  xtick={1,3,9}, xticklabels={$\normalmean_{1}$, $\normalmean_{2}$, $\normalmean_{3}$},
  ytick=\empty,
  enlargelimits=false, clip=false, axis on top,
  grid = major
  ]
  \addplot [draw=none, domain=0:2] {mix1} \closedcycle;
  \addplot [draw=none, domain=2:4] {mix1} \closedcycle;
  \addplot [draw=none, domain=4:10] {mix1} \closedcycle;
  \addplot [black!80] {mix1};
\end{axis}
\end{tikzpicture}
\qquad

\begin{tikzpicture}
\begin{axis}[
  no markers, domain=0:10, samples=200,
  axis lines*=left, %
  every axis y label/.style={at=(current axis.above origin),anchor=south},
  every axis x label/.style={at=(current axis.right of origin),anchor=west},
  height=3cm, width=5cm,
  xtick={1,3,9}, xticklabels={$\normalmean_{1}$, $\normalmean_{2}$, $\normalmean_{3}$},
  ytick=\empty,
  enlargelimits=false, clip=false, axis on top,
  grid = major
  ]
  \addplot [pattern=north west lines, pattern color=red, draw=none, domain=0:5] {mix1} \closedcycle;
  \addplot [pattern=crosshatch dots, pattern color=blue, domain=5:10] {mix1} \closedcycle;
  \addplot [black!80] {mix1};
\end{axis}
\end{tikzpicture}
\qquad
\begin{tikzpicture}
\begin{axis}[
  no markers, domain=0:10, samples=200,
  axis lines*=left, %
  every axis y label/.style={at=(current axis.above origin),anchor=south},
  every axis x label/.style={at=(current axis.right of origin),anchor=west},
  height=3cm, width=5cm,
  xtick={1,3,9}, xticklabels={$\normalmean_{1}$, $\normalmean_{2}$, $\normalmean_{3}$},
  ytick=\empty,
  enlargelimits=false, clip=false, axis on top,
  grid = major
  ]
  \addplot [pattern=north west lines, pattern color=red, draw=none, domain=0:2] {mix1} \closedcycle;
  \addplot [pattern=crosshatch dots, pattern color=blue, draw=none, domain=2:4] {mix1} \closedcycle;
  \addplot [pattern=north west lines, pattern color=red, draw=none, domain=4:10] {mix1} \closedcycle;
  \addplot [black!80] {mix1};
\end{axis}
\end{tikzpicture}
\qquad
\begin{tikzpicture}
\begin{axis}[
  no markers, domain=0:10, samples=200,
  axis lines*=left, %
  every axis y label/.style={at=(current axis.above origin),anchor=south},
  every axis x label/.style={at=(current axis.right of origin),anchor=west},
  height=3cm, width=5cm,
  xtick={1,3,9}, xticklabels={$\normalmean_{1}$, $\normalmean_{2}$, $\normalmean_{3}$},
  ytick=\empty,
  enlargelimits=false, clip=false, axis on top,
  grid = major
  ]
  \addplot [pattern=crosshatch dots, pattern color=blue, draw=none, domain=0:2] {mix1} \closedcycle;
  \addplot [pattern=north west lines, pattern color=red, draw=none, domain=2:10] {mix1} \closedcycle;
  \addplot [black!80] {mix1};
\end{axis}
\end{tikzpicture}

\caption{(top) Mixture of three Gaussians. (bottom) Different representations of a mixture of Gaussians as a mixture of two sub-Gaussians. Different fill patterns and colours represent different assignments of mixture components.}
\label{fig:threegaussians}
\end{figure}
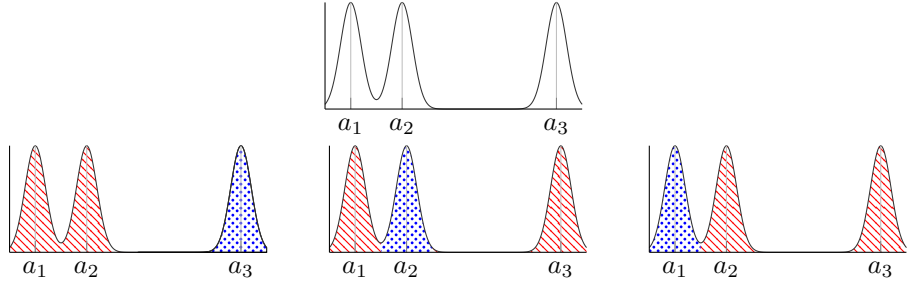

\begin{ex}[Sub-Gaussian mixtures are not identifiable]
\label{ex:danger:mix3}
Consider the mixture of three Gaussians $\mixfcn(\truemix)$
in Figure~\ref{fig:threegaussians}. We can write $\mixfcn(\truemix)$ as a mixture in four ways: In the top panel, $\mixfcn(\truemix)$ is represented uniquely as a mixture of three Gaussians. If we allow sub-Gaussian components, however, then the bottom panel shows three equally valid representations of $\mixfcn(\truemix)$ as a mixture of \emph{two} sub-Gaussians. 
Indeed, even if we assume the number of components $K$ is known and the component means are well-separated, $\mixfcn(\truemix)$ is non-identifiable 
as a mixture of sub-Gaussians: Just take $K=2$, $|\normalmean_{1}-\normalmean_{2}|>0$ and move $\normalmean_{3}$ arbitrarily far to the right. 
\end{ex}

\pgfmathdeclarefunction{gauss}{2}{%
  \pgfmathparse{1/(#2*sqrt(2*pi))*exp(-((x-#1)^2)/(2*#2^2))}%
}
\pgfmathdeclarefunction{mix1}{0}{%
  \pgfmathparse{0.33*gauss(3,0.4) + 0.33*gauss(1,0.4) + 0.33*gauss(9,0.4)}%
}

Much of the existing literature makes assumptions on the structure of the allowed $\truecmp_{k}$, which is evidently equivalent to restricting the supports of the mixing measures in $\mixset$ (e.g. Example~\ref{ex:param:mix}). Our focus, by contrast, will be to allow the components to take on essentially any shape while imposing regularity assumptions on the mixing measures $\truemix\in\mixset$. In this sense, we shift the focus from the properties of the ``local'' mixture components to the ``global'' properties of the mixture itself.

\section{Regularity and clusterability}
\label{sec:ident}

Fix an integer $K$ and let $\mixset\subset\mixing_{K}(\base)$ be a family of mixing measures. In particular, we assume that  $K$---the number of nonparametric mixtures---is known; in Section~\ref{sec:conc} we discuss the case where $K$ is unknown. In this section we study conditions that guarantee the injectivity of the embedding $\mixfcn:\mixset\to\mix(\mixset)$ using the procedure described in Section~\ref{sec:overview}. Throughout this section, it will be helpful to keep Figure~\ref{fig:overview} in mind for intuition.

\subsection{Projections}
\label{sec:ident:notion:proj}

Let $\{\projset_{L}\}_{L=1}^{\infty}$ be an indexed collection of families of mixing measures that satisfies the following:
\begin{enumerate}[label=(A\arabic*)]
\item\label{condn:A1} $\projset_{L}\subset\mixing_{L}(\base)$ for each $L$;
\item\label{condn:A2} $\{\projset_{L}\}$ is monotonic, i.e. $\projset_{L}\subset\projset_{L+1}$;
\item\label{condn:A3} $\mix(\projset_{L})$ is identifiable for each $L$.
\end{enumerate}

\noindent
The purpose of $\{\projset_{L}\}$ is to approximate $\true$ with a sequence of mixture distributions of increasing complexity, as quantified by the maximum number of atoms $L$, which will be taken to be much larger than $K$. Although our results apply to generic collections satisfying Conditions \ref{condn:A1}-\ref{condn:A3}, in the sequel we will consider the collection induced by a single subset $\projset\subset\mixing(\base)$ and defined by $\projset_{L}=\projset\cap\mixing_{L}(\base)$ (cf. \eqref{eq:def:finitemixtures}). We make the following assumption on $\projset$:
\begin{enumerate}[label=(A)]
\item\label{condn:A} The collection $\{\projset_{L}\}_{L=1}^{\infty}$ defined by $\projset_{L}=\projset\cap\mixing_{L}(\base)$ satisfies Condition \ref{condn:A3} for the family $\projset\subset\mixing(\base)$.
\end{enumerate}

\noindent
If $\projset$ satisfies Condition~\ref{condn:A}, then $\{\projset_{L}\}$ automatically satisfies Conditions \ref{condn:A1}-\ref{condn:A3}. 
Examples of families that satisfy Condition~\ref{condn:A} include exponential family mixture models under certain conditions \citep{barndorff1965}, for example Gaussian or Gamma mixtures \citep{teicher1963}.

Unless otherwise mentioned, we will assume $\projset$ satisfies Condition~\ref{condn:A}, with $\projset_{L}$ as defined therein.
Define the usual $\probm$-projection by
\begin{align}
\label{eq:defn:projmap}
\projmap_{L}\true
= \big\{Q\in \mix(\projset_{L}) : \probm(Q,\true) \le \probm(P,\true)\quad\forall P\in\mix(\projset_{L})\big\}.
\end{align}

\noindent
As long as $\projset$ is compact, the projection $\projmap_{L}\true$ is nonempty.
Furthermore, Condition \ref{condn:A3} implies that there exists a well-defined map $\mixingmap_{L}:\mix(\projset_{L})\to\projset_{L}$ that sends a mixture distribution to its mixing measure. 
With some abuse of notation, we will write $\mixingmap_{L}\true$ for $\mixingmap_{L}(\projmap_{L}\true)$, i.e.
\begin{align}
\label{eq:defn:mixmap}
\mixingmap_{L}\true
= \big\{\genmix\in \projset_{L} : \mixfcn(\genmix)\in\projmap_{L}\true\big\}.
\end{align}

Thus for any $\proj \in \projmap_{L}\true$, we can unambiguously define 
\begin{align}
\label{eq:defn:proj}
\proj 
=\sum_{\ell=1}^{L}\projwgt_{\ell}\projcmp_{\ell}
=\mixfcn(\projmix)
\quad\text{ and }\quad
\projmix
=\mixingmap_{L}(\proj).
\end{align}

\noindent
An example of the measure $\proj$ and its mixing measure $\projmix$
are depicted in Figure~\ref{fig:overview:proj}. 

\begin{remark}
\label{rem:proj:uniqueness}
We do not assume that $\projmap_{L}\true$ is unique, i.e. there may be more than one projection. This is because $\mix(\projset_{L})$ is a nonconvex set. We present our results in this setting, however, it may be simpler on a first reading to consider the special case where the projection is unique, i.e. $\projmap_{L}\true=\proj$ for each $L$. In this case, many of the definitions simplify: Consider for example Definition~\ref{defn:regular} and \eqref{eq:defn:sep:limsup} in the sequel.
\end{remark}

\begin{remark}
The number of overfitted mixture components $L$ will play an important but largely unheralded role in the sequel. For the most part, we will suppress the dependence of various quantities (e.g. $\proj$, $\projmix$) on $L$ for notational simplicity. In Section~\ref{sec:estimation:est}, we discuss how to choose $L$ given the sample size $n$; see Corollary~\ref{cor:main:conv}.
\end{remark}

\subsection{Assignment functions}
\label{sec:ident:notion:assign}

Any projection $\proj=\mixfcn(\projmix)=\sum_{\ell=1}^{L}\projwgt_{\ell}\projcmp_{\ell}$ as defined in \eqref{eq:defn:proj} is the best approximation to $\true$ from $\mix(\projset_{L})$, however, it contains many more components $L$ than the true number of \emph{nonparametric} components $K$. The next step is to find a way to ``cluster'' the components of $\proj$ into $K$ subgroups in such a way that each subgroup approximates some $\truecmp_{k}$. This is the second step \ref{steps:cluster} in our construction from Section~\ref{sec:overview}. To formalize this, we introduce the notion of \emph{assignment functions}.

Denote the set of all maps $\assg:[L]\to[K]$ by $\assignmaps$---a function $\assg\in\assignmaps$ represents a particular assignment of $L$ mixture components into $K$ subgroups. Thus, we will call $\assg$ an \emph{assignment function} in the sequel and a sequence $\assgseq$ of assignment functions such that $\assg_{L}\in\assignmaps$ will be called an \emph{assignment sequence}. The set of all assignment sequences is denoted by $\assignseq$.
For any $\genmix\in\projset_{L}$, write $\genproj=\mixfcn(\genmix)=\sum_{\ell=1}^{L}\genwgt_{\ell}q_{\ell}$. %
Given some $\assg\in\assignmaps$, define normalizing constants by
\begin{align}
\label{eq:defn:normalization}
\genclustwgt_{k}(\assg)
:= \sum_{\ell\in\assg^{-1}(k)}\genwgt_{\ell},
\quad k=1,\ldots,K.
\end{align}

\noindent
Denote the point mass concentrated at $\gencmp_{\ell}$ by $\delta_{\gencmp_{\ell}}$ and define 
\begin{align}
\label{eq:defn:mix:partition}
\genmix_{k}(\assg)
:=\frac1{\genclustwgt_{k}(\assg)}\sum_{\ell\in\assg^{-1}(k)}\genwgt_{\ell}\delta_{\gencmp_{\ell}},
\quad
\genproj_{k}(\assg)
:= \mixfcn(\genmix_{k}(\assg)).
\end{align}

\noindent
These quantities define a single, aggregate $K$-mixture by
\begin{align}
\label{eq:defn:mix:partition:2}
\genmix(\assg)
&:= \sum_{k=1}^{K} \genclustwgt_{k}(\assg)\delta_{\genproj_{k}(\assg)}, 
\quad
\genproj(\assg)
:= \mixfcn(\genmix(\assg))
=\sum_{k=1}^{K} \genclustwgt_{k}(\assg)\genproj_{k}(\assg).
\end{align}

\noindent
Since $\genproj_{k}(\assg)\in\mix(\projset_{0})$, $\genmix(\assg)$ is an atomic mixing measure whose atoms come from $\mix(\projset_{0})$. Informally, we hope that $\genproj_{k}(\assg)$ is able to approximate $\truecmp_{k}$, in a sense that will be made precise in the next section.

\subsection{Regular mixtures}
\label{sec:ident:regularity}

Given a nonparametric mixture $\mixfcn(\truemix)$, its $\probm$-proj-ection $\proj=\sum_{\ell=1}^{L}\projwgt_{\ell}\projcmp_{\ell}$, and an assignment function $\assg$, define $\projclustwgt_{k}(\assg)$ as in \eqref{eq:defn:normalization} and $\proj_{k}(\assg)$ and $\projmix_{k}(\assg)$ as in \eqref{eq:defn:mix:partition}. 
We'd like $\proj_{k}(\assg)$ to approximate $\truecmp_{k}$, but this is certainly not guaranteed for any $\assg$. %
A key step in our construction is to find such an assignment. Before finding such an assignment, however, we must first ask whether or not such an assignment \emph{exists}. The following notion of regularity encodes this assumption:

\begin{defn}[Regularity]
\label{defn:regular}
Suppose $\truemix\in\mixing_{K}(\base)$ and $\true=\mixfcn(\truemix)$. The mixing measure $\truemix$ is called \emph{$\projset$-regular} if: %
\begin{enumerate}[label=(\alph*)]
\item\label{defn:regular:a} 
There exists $L_{0}\ge 0$ such that $\projmap_{L}\true\ne\emptyset$ for each $L\ge L_{0}$ and\\$\lim_{L\to\infty}\proj=\true$ for every $\proj\in\projmap_{L}\true$;
\item\label{defn:regular:b} There exists an assignment sequence $\assgseq\in\assignseq$ such that %
\begin{align*}
\lim_{L\to\infty}\proj_{k}(\assg_{L})
&=\truecmp_{k} 
\quad\text{and}\quad
\lim_{L\to\infty}\projclustwgt_{k}(\assg_{L})
=\truewgt_{k}
\quad \forall\,k,
\, \forall\proj\in\projmap_{L}\true.
\end{align*}
\end{enumerate}

\noindent
When $\truemix$ is $\projset$-regular, we will also call $\mixfcn(\truemix)$ $\projset$-regular. 
\end{defn}

\begin{defn}[Regular assignment sequences]
\label{defn:regularseq}
Given a regular mixing measure $\truemix$, denote set of all assignment sequences $\assgseq$ such that Definition~\ref{defn:regular}\ref{defn:regular:b} holds by $\regassg(\truemix)$. An arbitrary assignment sequence $\assgseq\in\regassg(\truemix)$ will be called a \emph{regular assignment sequence}, or \emph{$\truemix$-regular} when we wish to emphasize the underlying mixing measure.
\end{defn}

\noindent
Whether or not a mixing measure is regular depends on both $\projset$ and $\probm$, although the dependence on $\probm$ will typically be suppressed. When we wish to emphasize this dependence, we will say $\truemix$ is \emph{$\projset$-regular under $\probm$}. Clearly, $\truemix$ is $\projset$-regular under the Hellinger metric if and only if it is $\projset$-regular under the variational metric.

The following examples construct several families of regular mixing measures, as well as an example where regularity fails. Proofs of these claims can be found in Appendix~\ref{app:ex}. 

\begin{ex}[Disjoint components]
\label{ex:disjmix}
Let $\base=\R$. Assume each $\truecmp_{k}$ has a density $\truedenscmp_{k}$ with respect to some dominating measure $\lebesgue$, and there exist disjoint intervals $E_{k}:=[\intlo_{k},\intup_{k}]\subset\R$ such that $\supp(\truedenscmp_{k}) \subseteq E_{k}$. Then the resulting mixing measure $\truemix$ is $\mixgaussians$-regular under both the Hellinger and variational metrics (Lemma~\ref{lem:disjmix:reg}). Furthermore, this example can be generalized to measures on $\R^{d}$ whose supports are contained in disjoint convex sets.
\end{ex}

\begin{ex}[Mixtures of finite mixtures]
\label{ex:metamix}
Fix $\projset\subset\mixing(\base)$ satisfying Condition~\ref{condn:A} and assume $\truecmp_{k}=\mixfcn(\gaussmix_{k})$ for each $k$, where $\gaussmix_{k}\in\projset_{0}$, i.e. $\gaussmix_{k}$ is a finite mixture model, but note that no upper bound is imposed on the number of components in each $\gaussmix_{k}$. Define $B_{k}:=\supp(\gaussmix_{k})$ and assume that $B_{1},\ldots,B_{K}$ are disjoint. Then $\truemix$ is $\projset$-regular under any metric $\probm$ (Lemma~\ref{lem:metamix:reg}). 
\end{ex}

\begin{ex}[Mixtures of infinite mixtures]
\label{ex:identmix}
In fact, the previous example can be generalized quite substantially.
Let $\projset\subset\mixing(\base)$ be compact and identifiable.
Assume that $\truecmp_{k}=\mixfcn(\gaussmix_{k})$ for each $k$, where $\gaussmix_{k}\in\projset$. For example, $\truecmp_{k}$ could be a potentially infinite convolutional mixture (see \citep{nguyen2013} for details), such as an infinite mixture of Gaussians with $\gaussmix_{k}\in\mixgaussians({\normalmeanset,\normalvar})$ (Example~\ref{ex:constrain:mix}). Define $B_{k}:=\supp(\gaussmix_{k})$ and assume that (a) $B_{1},\ldots,B_{K}$ are disjoint, compact sets and (b) Each $B_{k}$ is a $\gaussmix$-continuity set where $\gaussmix:=\sum_{k}\truewgt_{k}\gaussmix_{k}$. Then $\truemix$ is $\projset$-regular in both the Hellinger and variational metrics (Lemma~\ref{lem:identmix:reg}). 
\end{ex}

\begin{ex}[Failure of regularity]
\label{ex:regularityfail}
Let $g_{\pm}\sim\normalN(\pm \normalmean,1)$ and $G\sim\normalN(0,\sigma^{2})$ where $\sigma^{2}>0$, and define for some $0<\beta_{1}<\beta_{2}<1$,
$\true = \tfrac12\truecmp_{1} + \tfrac12\truecmp_{2}$,
$\truecmp_{1} \propto (1-\beta_{1}-\beta_{2})g_{+} + \tfrac{\beta_{1}}{2}G$, and
$\truecmp_{2} \propto \beta_{2} g_{-} + \tfrac{\beta_{1}}{2}G$.
In this example, $K=2$. If $\projset_{L}=\mixgaussians_{L}$, then for any $L>3$, $\proj=\true$, and there is no way to cluster the 3 components into 2 mixtures of Gaussians that approximate the $\truecmp_{k}$. The problem here is that $\truecmp_{1}$ and $\truecmp_{2}$ ``share'' the same Gaussian component $G$, which evidently cannot be assigned to both $\truecmp_{1}$ and $\truecmp_{2}$.
\end{ex}

We conclude by pointing out that, in addition to the concrete examples discussed above, in general the set of regular mixing measures is quite large:

\begin{lemma}
\label{lem:reg:dense}
Let $\dTV$ be the variational distance on $\probs(\R^{p})$ and let $\mixm$ be the induced Wasserstein metric on $\mixing_{K}(\R^{p})$. 
Then for any $\truemix\in\mixing_{K}(\R^{p})$ and $\eps>0$, there exists a $\mixgaussians$-regular mixing measure $\tick{\truemix}\in\mixing_{K}(\R^{p})$ such that $\mixm(\tick{\truemix},\truemix)<\eps$. In particular, the set of  $\mixgaussians$-regular mixing measures is dense in $\mixing_{K}(\R^{p})$.
\end{lemma}

\noindent
In fact, the proof is constructive: The family defined in Example~\ref{ex:metamix} with $\projset=\mixgaussians$ is dense in $\mixing_{K}(\R^{p})$.

\subsection{Clusterable families}
\label{sec:ident:clust}

If a mixing measure $\truemix$ is $\projset$-regular, then the $\probm$-projections of $\mixfcn(\truemix)$ can always be grouped in such a way that each group approximates the nonparametric component $\truecmp_{k}$ and its mixing weight $\truewgt_{k}$. We have not said anything yet about \emph{how} one might find such an assignment; only that it exists. The following condition asserts that regular assignments can be determined from the projections $\proj_{L}$:

\begin{defn}[Clusterable family]
\label{defn:clusterable}
A family of mixing measures $\mixset\subset\mixing(\base)$ is called a \emph{$\projset$-clusterable family}, or just a \emph{clusterable family}, if 
\begin{enumerate}[label=(\alph*)]
\item\label{defn:clusterable:a} $\truemix$ is $\projset$-regular for all $\truemix\in\mixset$;
\item\label{defn:clusterable:b} There exists a function $\clustmap_{L}:\mixingmap_{L}(\mixset)\to\assignmaps$ such that $\{\clustmap_{L}(\projmix)\}\in\assignseq(\truemix)$ for every $\truemix\in\mixset$.
\end{enumerate}

\noindent
The resulting mixture model $\mix(\mixset)$ is called a \emph{clusterable mixture model}. If $\truemix$ belongs to a clusterable family, we shall call both $\truemix$ and $\true=\mixfcn(\truemix)$ \emph{clusterable measures}.
\end{defn}

\noindent
As with regularity, clusterability depends on both $\projset$ and $\probm$. When we wish to emphasize this dependence, we will say $\truemix$ is \emph{$\projset$-clusterable under $\probm$}.
The terminology ``clusterable'' is intended to provoke the reader into imagining $\clustmap_{L}$ as a cluster function that ``clusters'' the $L$ components and $L$ weights of $\proj$ together in such a way that $\projmix(\assg)$ approximates $\truemix$.
More precisely, Definition~\ref{defn:clusterable}\ref{defn:clusterable:b} means that for every $\truemix\in\mixset$, if we let $\projmix=\mixingmap_{L}(\projmap_{L}(\mixfcn(\truemix)))$, then $\assg_{L}=\clustmap_{L}(\projmix)$ defines a regular assignment sequence (Definition~\ref{defn:regularseq}).

\subsection{Separation and clusterability}
\label{sec:estimation:sep}

In this section, we construct an explicit cluster function $\clustmap_{L}$ via single-linkage clustering.

Given $\genmix\in\projset_{L}$ with atoms $\gencmp_{\ell}$, define the $\probm$-diameter of $\genmix$ by
\begin{align*}
\Hdiam(\genmix)
:= \sup\{\probm(\gencmp,\gencmp') : \gencmp,\gencmp'\in\conv(\supp(\genmix))\}
\end{align*}

\noindent
where $\conv(\,\cdot\,)$ is the convex hull in $\probs(\base)$.
Recalling \eqref{eq:defn:mix:partition:2}, define for any $\assg\in\assignmaps$
\begin{align}
\label{eq:defn:sep}
\separ(\genmix(\assg))
:= \sup_{k}\Hdiam(\genmix_{k}(\assg)) + \sup_{k}\probm(\truecmp_{k},\genproj_{k}(\assg)).
\end{align}

\noindent
We will be interested in the special case $\genmix=\projmix$: $\Hdiam(\projmix_{k}(\assg))$ quantifies how ``compact'' the mixture component $\proj_{k}(\assg)$ is and $\separ(\projmix(\assg))$ is a measure of separation between the mixture components $\truecmp_{k}$.
Finally, define the $\probm$-distance matrix by
\begin{align}
\label{eq:defn:distmat}
\gendist(\genmix)
=(\probm(\gencmp_{i},\gencmp_{j}))_{i,j=1}^{L}.
\end{align}

Our goal is to show that if the atoms of $\truemix$ are sufficiently well-separated, then the cluster assignment $\assg$ can be reconstructed by clustering the distance matrix $\projdist=\gendist(\projmix)=(\probm(\projcmp_{i},\projcmp_{j}))_{i,j=1}^{L}$ (hence the choice of terminology \emph{clusterable}). 
More precisely, we make the following definition:
\begin{defn}[Separation]
\label{defn:separated}
A mixing measure $\truemix\in\mixing_{0}(\base)$ is called \emph{$\delta$-separated} if $\inf_{i\ne j}\probm(\truecmp_{i},\truecmp_{j}) > \delta$ for some $\delta>0$.

\end{defn}

\noindent
It turns out that separation of the order $\separ(\projmix(\assg))$ (cf. \eqref{eq:defn:sep}) is sufficient to define a cluster function:

\begin{prop}
\label{prop:pop:cluster}
Let $\truemix\in\mixing_{K}(\base)$. Let $\proj\in\projmap_{L}\true$ be a $\probm$-projection of $\true$ for some $L\ge K$. 
Then for any $\assg\in\assignmaps$ such that $\truemix$ is $4\separ(\projmix(\assg))$-separated, 
\begin{align}
\label{eq:prop:pop:cluster:1}
\assg(i)=\assg(j)
&\iff
\probm(\projcmp_{i}, \projcmp_{j}) \le \separ(\projmix(\assg)), \\
\label{eq:prop:pop:cluster:2}
\assg(i)\ne\assg(j)
&\iff
\probm(\projcmp_{i}, \projcmp_{j}) \ge 2\separ(\projmix(\assg)).
\end{align}

\noindent
Moreover, $\assg$ can be recovered by single-linkage clustering on $\projdist$.
\end{prop}

\noindent
Thus, the assignment $\assg$ can be recovered by single-linkage clustering of $\projdist$ \emph{without knowing the optimal threshold $\separ(\projmix(\assg))$}. 

Now suppose $\truemix$ is a regular mixing measure and let $\assgseq\in\assignseq(\truemix)$. Define
\begin{align}
\label{eq:defn:sep:limsup}
\separ(\truemix)
:= \limsup_{L\to\infty}\sup_{\projmix\in\mixingmap_{L}\true}\sup_{\assgseq\in\assignseq(\truemix)}\separ(\projmix(\assg_{L})).
\end{align}

\noindent
As a consequence of regularity, the second term in \eqref{eq:defn:sep} tends to zero as $L\to\infty$, so that $\separ(\truemix)$ can be interpreted as a measure of the asymptotic diameter of the approximating mixtures $\proj_{k}(\assg_{L})$. For example, when the $\probm$-projection $\proj$ is unique (Remark~\ref{rem:proj:uniqueness}) the definition in \eqref{eq:defn:sep:limsup} simplifies to $\separ(\truemix)= \limsup_{L}\,\sup_{k}\Hdiam(\projmix_{k}(\assg_{L}))$. The following corollary, which is an immediate consequence of Proposition~\ref{prop:pop:cluster}, shows that control over $\separ(\truemix)$ is sufficient for $\mixset$ to be clusterable:

\begin{cor}
\label{cor:sep:ident}
Suppose $\mixset\subset\mixing_{K}(\base)$ is a family of regular mixing measures such that for every $\truemix\in\mixset$ there exists $\sepconst>0$ such that $\truemix$ is $(4+\sepconst)\separ(\truemix)$-separated. Then $\mixset$ is clusterable.
\end{cor}

\noindent
Thus, we have a practical separation condition under which a regular mixture model becomes identifiable:
\begin{align}
\label{eq:condn:sep}
\inf_{i\ne j}\probm(\truecmp_{i},\truecmp_{j})
> (4+\sepconst)\separ(\truemix).
\end{align}

\noindent
In the limit $L\to\infty$, the nonparametric components $\truecmp_{k}$ must be separated by a gap proportional to the $\probm$-diameters of the approximating mixtures $\proj_{k}(\assg_{L})$. This highlights the issue in Example~\ref{ex:danger:mix3}---although the means can be arbitrarily separated, as we increase the separation, the diameter of the components continues to increase as well. Thus, the $\truecmp_{k}$ cannot be chosen in a haphazard way (see also Example~\ref{ex:regularityfail}). Crucially, however, we make no assumptions on the shape of the mixture components.

\begin{ex}[Example of separation]
\label{ex:clust:family}
Take $\base=\R$ and let $\projset=\mixgaussians({\normalmeanset,\normalvar})$ be a family of convolutional mixtures of Gaussians (Example~\ref{ex:constrain:mix}). In Example~\ref{ex:identmix}, we claimed that as long as the $\gaussmix_{k}$ have disjoint supports, this family is $\mixgaussians({\normalmeanset,\normalvar})$-regular. To determine when a mixing measure $\truemix$ is $\mixgaussians({\normalmeanset,\normalvar})$-clusterable, it suffices to check \eqref{eq:condn:sep}. For this, we bound the $\probm$-diameters $\sup_{k}\Hdiam(\projmix_{k}(\assg_{L}))$ for large $L$. If $\gencmp_{\ell}\sim\normalN(\normalmean_{\ell},\normalvar^{2})$ and $\gencmp_{\tick{\ell}}\sim\normalN(\normalmean_{\tick{\ell}},\normalvar^{2})$ are both in $\supp(\gaussmix_{k})$, then it is easy to check that as long as
\begin{align*}
|\normalmean_{\ell}-\normalmean_{\tick{\ell}}|
\le \sqrt{8\normalvar^{2}
    \log\big(1 + \tfrac{\minsep}{\bb-\minsep}\big)
},
\quad
\minsep:=\inf_{i\ne j}\probm(\truecmp_{i},\truecmp_{j}),
\end{align*}

\noindent
the separation condition \eqref{eq:condn:sep} holds. 
\end{ex}

The separation condition \eqref{eq:condn:sep} is quite weak, but no attempt has been made here to optimize this lower bound. For example, a minor tweak to the proof can reduce the constant of 4 to any constant $b>2$. Although we expect that a more careful analysis can weaken this condition, our main focus here is to present the main idea behind identifiability and its connection to clusterability and separation, so we save such optimizations for future work. Further, although Proposition~\ref{prop:pop:cluster} justifies the use of single-linkage clustering in order to group the components $\{\projcmp_{\ell}\}$, one can easily imagine using other clustering schemes. Indeed, since the distance matrix $\projdist$ is always well-defined, we could have applied other clustering algorithms such as complete-linkage hierarchical clustering, $K$-means, or spectral clustering to $\projdist$ to define an assignment sequence $\assgseq$. Any condition on $\projdist$ that ensures a clustering algorithm will correctly reconstruct a regular assignment sequence then yields an identification condition in the spirit of Proposition~\ref{prop:pop:cluster}. For example, if the means of the overfitted components $\projcmp_{\ell}$ are always well-separated, then simple algorithms such as $K$-means could suffice to identify a regular assignment sequence. 
This highlights the advantage of our abstract viewpoint, in which the specific forms of both the assignment sequence $\assgseq$ and the cluster functions $\clustmap_{L}$ are left unspecified.

\section{Identifiability and estimation}
\label{sec:estimation}

We now turn our attention to the problem of identifying and learning a mixing measure $\truemix$ from data.

\subsection{Identifiability of nonparametric mixtures}
\label{sec:ident:main}

According to the next theorem, clusterability is sufficient to identify a nonparametric mixture model.

\begin{thm}
\label{thm:main:ident}
If $\mixset$ is a $\projset$-clusterable family then the mixture model $\mix(\mixset)$ is identifiable. 
\end{thm}

\noindent
As illustrated by the cautionary tales from Examples~\ref{ex:danger:mix3} and~\ref{ex:regularityfail}, identification in nonparametric mixtures is a subtle problem, and this theorem thus provides a powerful general condition for identifiability in nonparametric problems.

\begin{figure}[t]
\centering
\begin{subfigure}{1\textwidth}
\centering
\includegraphics[width=\textwidth]{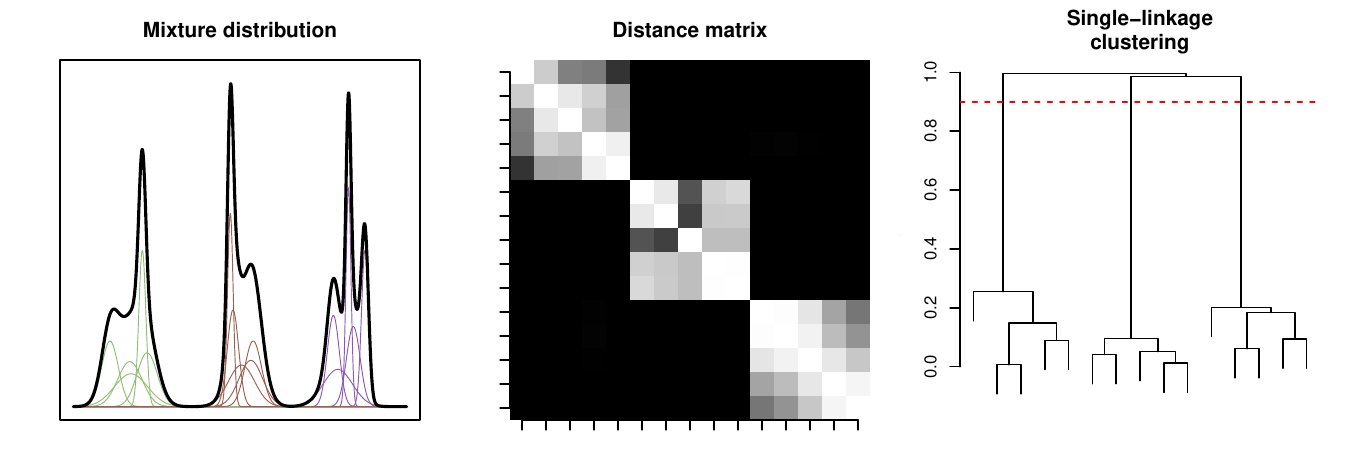}
\caption{$K=3$ and well-separated means.}
\label{fig:mog:sep}
\end{subfigure}%
\\
\begin{subfigure}{1\textwidth}
\centering
\includegraphics[width=\textwidth]{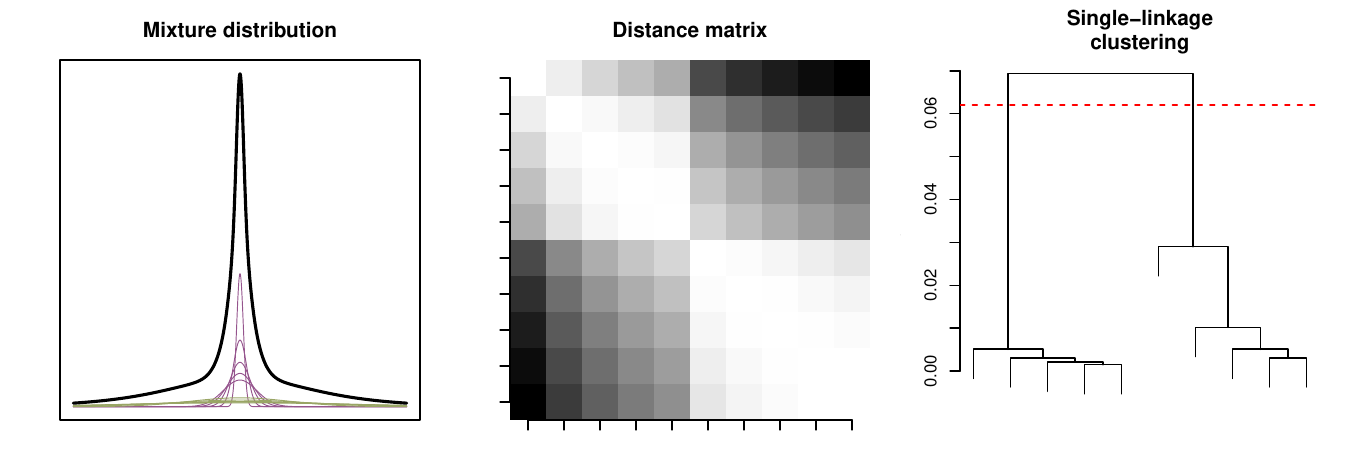}
\caption{$K=2$ and identical means.}
\label{fig:mog:overlap}
\end{subfigure}%
\caption{Illustrating Theorem~\ref{thm:main:ident} with Example~\ref{ex:metamix}. (left) The original mixture distribution (thick black line) is a mixture of finite Gaussian mixtures. Each Gaussian component is coloured according to its membership in different $\gaussmix_{k}$. (middle) The true distance matrix $\projdist$. (right) Results of single-linkage clustering on $\projdist$, cut to find the correct number of clusters.}
\label{fig:mog}
\end{figure}

Two examples of Theorem~\ref{thm:main:ident} are illustrated in Figure~\ref{fig:mog}. When the means are well-separated as in Figure~\ref{fig:mog:sep}, it is easy to see how single-linkage clustering is able to discover a correct assignment. 
Since $\probm$-separation is a weaker criterion than mean separation, however, Theorem~\ref{thm:main:ident} does not require that the mixture distributions in $\mix(\mixset)$ have components with well-separated means. In fact, each $\truecmp_{k}$ could have identical means (but different variances) and still be well-separated. This is illustrated in Figure~\ref{fig:mog:overlap}. This suggests that identifiability in mixture models is more general than what is needed in typical clustering applications, where a model such as Figure~\ref{fig:mog:overlap} would not be considered to have two distinct clusters. The subtlety here lies in interpreting clustering in $\probs(\base)$ (i.e. of the $\projcmp_{\ell}$) vs. clustering in $\base$ (i.e. of samples $\rv^{(i)}\sim\true$), the latter of which is the interpretation used in data clustering.

\subsection{Estimation of clusterable mixtures}
\label{sec:estimation:est}

We now discuss how to estimate $\truemix$ from data $\rv^{(1)},\ldots,\rv^{(n)}\iid\true$. Throughout this section, we assume that $\projmix=\projmix_{L}\in\mixingmap_{L}\true$ is arbitrary.

For each $L\ge K$, let $\estmix_{L,n}=\estmix_{L}(\rv^{(1)},\ldots,\rv^{(n)})$ be a $\mixm$-consistent estimator of $\projmix_{L}$, where we have written $\estmix_{L,n}$ and $\projmix_{L}$ to emphasize the dependence on $L$ and $n$. That is, $\{\estmix_{L,n}\}$ is a sequence of estimators and for each $L$, $\lim_{n\to\infty}\mixm(\estmix_{L,n},\projmix_{L})=0$. For example, $\estmix_{L,n}$ could be the minimum Hellinger distance estimator (MHDE) from \citet{beran1977} (see Appendix~\ref{app:mdhe} for details). Since $L$ is a known quantity, the corresponding estimation problems are always well-specified, i.e. both $\estmix_{L,n}$ and $\projmix_{L}$ have the same, known number of components. In the sequel, we will omit the dependence of $\projmix=\projmix_{L}$ on $L$ and $\estmix=\estmix_{L,n}$ on $L$ and $n$ for brevity.
Write
\begin{align}
\label{eq:def:estproj}
\estproj
:=\mixfcn(\estmix)
=\sum_{\ell=1}^{L}\estwgt_{\ell}\estcmp_{\ell}.
\end{align}

\noindent
Without loss of generality, assume that the atoms are re-arranged so that $\sup_{\ell}\probm(\estcmp_{\ell},\projcmp_{\ell})\to0$ (see Lemma~\ref{lem:wasserstein:atom}).

\begin{prop}
\label{prop:est:cluster}
Let $\truemix\in\mixing_{K}(\base)$. 
Let $\proj\in\projmap_{L}\true$ be a $\probm$-projection of $\true$ for some $L\ge K$.
Suppose further that $L$, $\assg\in\assignmaps$, and $n$ satisfy
\begin{align}
\label{eq:prop:est:cluster:1}
3\,\sup_{\ell}\probm(\estcmp_{\ell},\projcmp_{\ell})
- 2\,\sup_{k}\probm(\proj_{k}(\assg),\truecmp_{k})
< \sup_{k}\Delta(\projmix_{k}(\assg)).
\end{align}

\noindent
Define
\begin{align*}
\wh{\separ}
:= 2\,\sup_{\ell}\probm(\estcmp_{\ell},\projcmp_{\ell}) 
+ \sup_{k}\Hdiam(\projmix_{k}(\assg)).
\end{align*}

\noindent
If $\truemix$ is $4\separ(\projmix(\assg))$-separated,
then $\probm(\estcmp_{i},\estcmp_{j})\le\wh{\separ}$ if and only if $\assg(i)=\assg(j)$, and the assignment function $\assg$ can be recovered by single-linkage clustering on $\estdist=\gendist(\estmix)$.
\end{prop}

\noindent
Proposition~\ref{prop:est:cluster} is a finite sample result that holds as long as $L$ and $n$ satisfy \eqref{eq:prop:est:cluster:1}, which is guaranteed 
as long as $\truemix$ is $\projset$-regular (i.e. since in this case the left side tends to zero).

For each $L$ and $n$, let $\estassg=\estassg_{L,n}\in\assignmaps$ denote the assignment map defined in Proposition~\ref{prop:est:cluster}. With this notation, another way to phrase this result is that under \eqref{eq:prop:est:cluster:1}, we have $\estassg=\assg$. In other words, single-linkage clustering of $\estdist$ yields the same clusters as the assignment $\assg$. This suggests we use $\estmix(\estassg)$ as an estimator of $\truemix$. More precisely:
\begin{enumerate}
\item Choose $L\ge K$ sufficiently large;
\item Estimate $\estmix=\estmix(\rv^{(1)},\ldots,\rv^{(n)})$;
\item Define $\estassg=\estassg(\rv^{(1)},\ldots,\rv^{(n)})$ by single-linkage clustering on $\estdist$;
\item Return $\estmix(\estassg)$.
\end{enumerate}

\noindent
In order for this to be an estimator, we must have a precise rule for selecting $L=L_{n}$; see Corollary~\ref{cor:main:conv} below and its discussion for details.

The following theorem provides conditions under which $\estmix(\estassg)$ consistently estimates $\truemix$:
\begin{thm}
\label{thm:main:learn}
Suppose $\truemix$ is a regular mixing measure such that $\truemix$ is $(4+\sepconst)\separ(\truemix)$-separated for some $\sepconst>0$. Then 
\begin{align}
\label{eq:thm:main:learn:1a}
\lim_{L\to\infty}\lim_{n\to\infty}\mixm(\estmix(\estassg), \truemix)
&=0.
\end{align}

\end{thm}

\noindent
In particular, \eqref{eq:thm:main:learn:1a} implies  that
\begin{align*}
\lim_{L\to\infty}\lim_{n\to\infty}\probm(\estproj_{k}(\estassg), \truecmp_{k})
&= 0
\quad\text{and}\quad
\lim_{L\to\infty}\lim_{n\to\infty}|\projclustwgt_{k}(\estassg) - \truewgt_{k}| = 0.
\end{align*}

\noindent
Thus, we have a Wasserstein consistent estimate of $\truemix$ and $\probm$-consistent estimates of the component measures $\truecmp_{k}$. As stated, Theorem~\ref{thm:main:learn} has an important drawback: Without a rule for choosing $L=L_{n}$ as a function of the sample size $n$, $\estmix(\estassg)$ is not a proper estimator. This is the cost of abstraction that allows us to state such a theorem for general metric spaces and probability measures. Fortunately, in special cases we can make the dependence on $n$ explicit: Recall the convolutional mixture model described in Example~\ref{ex:constrain:mix}. We have already shown that this family is both regular (Example~\ref{ex:identmix}) and clusterable (Example~\ref{ex:clust:family}). Combining these results with a rule for choosing $L=L_{n}$, the following corollary provides a practical setting in which all of the assumptions laid out in Theorem~\ref{thm:main:learn} are satisfied.

\begin{cor}
\label{cor:main:conv}
Let $\truemix=\sum_{k=1}^{K}\truewgt_{k}\gaussmix_{k}$ be as in Example~\ref{ex:identmix} with $\gaussmix_{k}\in\mixgaussians({\normalmeanset,\normalvar})$ for each $k$. Define $B_{k}:=\supp(\gaussmix_{k})$ and assume that (a) $B_{1},\ldots,B_{K}$ are disjoint, compact sets and (b) Each $B_{k}$ is a $\gaussmix$-continuity set where $\gaussmix:=\sum_{k}\truewgt_{k}\gaussmix_{k}$. Define $\minsep:=\inf_{i\ne j}\probm(\mixfcn(\gaussmix_{i}),\mixfcn(\gaussmix_{j}))$ and assume further that 
\begin{align*}
\sup_{\gencmp_{\ell},\gencmp_{\tick{\ell}}\in\supp(\gaussmix_{k})}|\E\gencmp_{\ell}-\E\gencmp_{\tick{\ell}}|
\le \sqrt{8\normalvar^{2}
    \log\big(1 + \tfrac{\minsep}{\bb-\minsep}\big)
},
\quad
\text{for all $k$.}
\end{align*}
Then taking $L_{n}\asymp n^{2/3}/\log^{1/3}n$, we have
\begin{align}
\label{eq:cor:main:conv:1a}
\lim_{n\to\infty}\mixm(\estmix(\estassg_{L_{n},n}), \truemix)
&= 0. %
\end{align}
\end{cor}

\noindent
The proof follows immediately from Theorem~2 in \citet{nguyen2013} and the results (e.g. Theorem~5) in \citet{genovese2000}.

Finally, in applications, it will often be useful to strengthen $\probm$-convergence to \emph{uniform} convergence of the densities (assuming they exist). When the families $\projset_{L}$ are equicontinuous, this is guaranteed by Theorem~1 of \citet{sweeting1986}. We store this corollary away here for future use:

\begin{cor}
\label{cor:conv:uniform}
Let $\estdens_{k}(\estassg)$ be the density of $\estproj_{k}(\estassg)$ and $\truedenscmp_{k}$ be the density of $\truecmp_{k}$. If the families $\projset_{L}$ are equicontinuous for all $L$ and $\estproj_{k}(\estassg)$ converges weakly to $\truecmp_{k}$, then $\lim_{L\to\infty}\lim_{n\to\infty}\estdens_{k}(\estassg)=\truedenscmp_{k}$, where the limits are understood both pointwise and uniformly over compact subsets of $\base$.
\end{cor}

\noindent
The assumption that $\estproj_{k}(\estassg)$ converges weakly to $\truecmp_{k}$ restricts the choice of $\probm$, although it allows most reasonable metrics including Hellinger, variational, and Wasserstein, for example. Moreover, even weaker assumptions than equicontinuity are possible \citep{cuevas1991}.

\section{Bayes optimal clustering}
\label{sec:bayesopt}

As an application of the theory developed in Sections~\ref{sec:ident} and~\ref{sec:estimation}, we extend model-based clustering \citep{bock1996,fraley2002} to the nonparametric setting. Given samples from $\truemix$, we seek to partition these samples into $K$ clusters. More generally, $\truemix$ defines a partition of the input space $\base$, which can be formalized as a function $c:\base\to[K]$, where $K$ is the number of partitions or ``clusters''. 
First, let us recall the classical Gaussian mixture model (GMM): If $f_{1}(\cdot;\normalmean_{1},\normalvar_{1}),\ldots,f_{K}(\cdot;\normalmean_{K},\normalvar_{K})$ is a collection of Gaussian density functions, then for any choice of $\truewgt_{k}\ge0$ such that $\sum_{k}\truewgt_{k}=1$ the combination
\begin{align}
\label{eq:gmm}
F(z)
= \sum_{k=1}^{K}\truewgt_{k}f_{k}(z;\normalmean_{k},\normalvar_{k}); 
\quad z\in\R^{d}
\end{align}
is a GMM. The model \eqref{eq:gmm} is of course equivalent to the integral \eqref{eq:def:mixmeasure} (see also Example~\ref{ex:param:mix}), and the Gaussian densities $f_{k}(z;\normalmean_{k},\normalvar_{k})$ can obviously be replaced with any family of parametric densities. %

Intuitively, the density $F$ has $K$ distinct clusters given by the $K$ Gaussian densities $f_{k}$, defining what we call the \emph{Bayes optimal partition} over $\base$ into regions where each of the Gaussian components is most likely. It should be obvious that as long as a mixture model $\mix(\mixset)$ is identifiable, the Bayes optimal partition will be well-defined and has a unique interpretation in terms of distinct clusters of the input space $\base$. Thus, the theory developed in the previous sections can be used to extend these ideas to the nonparametric setting. Since the clustering literature is full of examples of datasets that are not well-approximated by parametric mixtures \citep[e.g.][]{ng2001,ultsch2005}, there is significant interest in such an extension.
In the remainder of this section, we will apply our framework to this problem. First, we discuss identifiability issues with the concept of a Bayes optimal partition (Section~\ref{sec:bayesopt:part}). Then, we provide conditions under which a Bayes optimal partition can be learned from data (Section~\ref{sec:bayesopt:learn}).

\subsection{Bayes optimal partitions}
\label{sec:bayesopt:part}

Throughout the rest of this section, we assume that $\base$ is compact and all probability measures are absolutely continuous with respect to some base measure $\lebesgue$, and hence have density functions. 
Assume $\true$ is fixed and write $\truedens=\truedens_{\true}$ for the density of $\true$ and $\truedenscmp_{k}$ for the density of $\truecmp_{k}$. Thus whenever $\true$ is a finite mixture we can write 
\begin{align}
\label{eq:basic:mix}
\truedens
= \int \truedenscmp_{\truecmp}\,d\truemix(\truecmp)
= \sum_{k=1}^{K}\truewgt_{k}\truedenscmp_{k}.
\end{align}

\noindent
For any $\truemix\in\mixing_{K}(\base)$, define the usual Bayes classifier \citep[e.g.][]{devroye2013}:
\begin{align}
\label{eq:defn:partition}
c_{\truemix}(x)
:= \argmax_{k\in[K]}\truewgt_{k}\truedenscmp_{k}(x).
\end{align}

\noindent
The classifier $c_{\truemix}$ is only well-defined up to a permutation of the labels (i.e. any labeling of $\supp(\truemix)$ defines an equivalent classifier). Furthermore, $c_{\truemix}(x)$ not properly defined when $\truewgt_{i}\truedenscmp_{i}(x)=\truewgt_{j}\truedenscmp_{j}(x)$ for $i\ne j$. To account for this, define an exceptional set
\begin{align}
\label{eq:defn:eqset}
\eqset
:= \bigcup_{i\ne j}\{x\in\base: \truewgt_{i}\truedenscmp_{i}(x) = \truewgt_{j}\truedenscmp_{j}(x)\},
\end{align}

\noindent
In principle, $\eqset$ should be small---in fact it will typically have measure zero---hence we will be content to partition $\base_{0}=\base-\eqset$. Recall that a \emph{partition} of a space $\base$ is a family of subsets $A_{k}\subset \base$ such that $A_{k}\cap A_{k'}=\emptyset$ for all $k\ne k'$ and  $\cup_{k}A_{k}=\base$. We denote the space of all partitions of $\base$ by $\partn(\base)$. 

The following definition is standard \citep[e.g. ][]{fraley2002,chacon2015}:

\begin{defn}[Bayes optimal partition]
\label{defn:bayesopt}
Define an equivalence relation on $\base_{0}$ by declaring 
\begin{align}
\label{eq:defn:bayesopt}
x\sim y
\iff
c_{\truemix}(x)=c_{\truemix}(y).
\end{align}

\noindent
This relation induces a partition on $\base_{0}$ which we denote by $\pt_{\truemix}$ or $\pt(\truemix)$. This partition is known as the \emph{Bayes optimal partition}.
\end{defn}

\begin{remark}
Although the function $c_{\truemix}$ is only unique up to a permutation, the partition defined by \eqref{eq:defn:bayesopt} is always well-defined and independent of the permutation used to label the $\truecmp_{k}$.
\end{remark}

Given samples from the mixture distribution $\true=\mixfcn(\truemix)$, we wish to learn the Bayes optimal partition $\pt_{\truemix}$. Unfortunately, there is---yet again---an identifiability issue: If there is more than one mixture measure $\truemix$ that represents $\true$, the Bayes optimal partition is not well-defined.

\begin{ex}[Non-identifiability of Bayes optimal partition]
\label{ex:danger:part}
In Example~\ref{ex:danger:mix3} and Figure~\ref{fig:threegaussians}, we have four valid representations of $\true$ as a mixture of sub-Gaussians. In all four cases, each representation leads to a different Bayes optimal partition, even though they each represent the same mixture distribution.
\end{ex}

\noindent
Clearly, if $\truemix$ is identifiable, then the Bayes optimal partition is automatically well-defined. Thus Theorem~\ref{thm:main:ident} immediately implies the following:

\begin{cor}
\label{cor:bayesopt:ident}
If $\mix(\mixset)$ is a clusterable mixture model, then there is a well-defined Bayes optimal partition $\pt_{\true}$ for any $\true\in\mix(\mixset)$.
\end{cor}

\noindent
In particular, whenever $\mix(\mixset)$ is clusterable it makes sense to write $c_{\true}$ and $\pt_{\true}$ instead of $c_{\truemix}$ and $\pt_{\truemix}$, respectively.
This provides a useful framework for discussing and analyzing partition-based clustering in nonparametric settings. As discussed previously, a $K$-clustering of $\base$ is equivalent to a function that assigns each $x\in\base$ an integer from $1$ to $K$, where $K$ is the number of clusters. Clearly, up to the exceptional set $\eqset$, \eqref{eq:defn:partition} is one such function. Thus, the Bayes optimal partition $\pt_{\true}$ can be interpreted as a valid $K$-clustering.

\subsection{Learning partitions from data}
\label{sec:bayesopt:learn}

Write $\true=\mixfcn(\truemix)$ and assume that $\truemix$ is identifiable from $\true$. Suppose we are given i.i.d. samples $\rv^{(1)},\ldots,\rv^{(n)}\iid\true$ and that we seek the Bayes optimal partition $\pt_{\true}=\pt_{\truemix}$. Our strategy will be the following:
\begin{enumerate}
\item Use a consistent estimator $\estmix$ to learn $\projmix$ for some $L\gg K$;
\item Theorem~\ref{thm:main:learn} guarantees that we can learn a cluster assignment $\estassg$ such that $\estmix(\estassg)$ consistently estimates $\truemix$;
\item Use $\pt(\estmix(\estassg))$ to approximate $\pt_{\truemix}=\pt_{\true}$.
\end{enumerate}

\noindent
The hope, of course, is that $\pt(\estmix(\estassg))\to\pt_{\true}$. 
There are, however, complications: What do we mean by convergence of partitions? Does $\pt(\estmix(\estassg))$ even converge, let alone converge to $\pt_{\true}$? 

Instead of working directly with the partitions $\pt(\estmix(\estassg))$, we will work with the Bayes classifier \eqref{eq:defn:partition}. 
Write $\estdenscmp_{\ell}$ and $\estdens$ for the densities of $\estcmp_{\ell}$ and $\estproj$, respectively, and 
\begin{align}
\label{eq:est:dens}
\estdens_{k}(\estassg)
:= \frac1{\estclustwgt_{k}}\sum_{\ell\in\estassg^{-1}(k)}\estwgt_{\ell}\estdenscmp_{\ell},
\quad
\estclustwgt_{k}(\estassg)
:= \sum_{\ell\in\estassg^{-1}(k)}\estwgt_{\ell}.
\end{align}

\noindent
Then $\estdens_{k}(\estassg)$ is the density of $\estproj_{k}(\estassg)$, where here and above we have suppressed the dependence on $\estassg$. Now define the estimated classifier (cf. \eqref{eq:defn:partition})
\begin{align}
\label{eq:est:clustmap}
\wh{c}(x)
:= c_{\estmix(\estassg)}(x)
= \argmax_{k\in[K]}\estclustwgt_{k}[\estdens_{k}(\estassg)](x).
\end{align}

\noindent
By considering classification functions as opposed to the partitions themselves, we may consider ordinary convergence of the function $\wh{c}$ to $c_{\true}$, which gives us a convenient notion of consistency for this problem. 
Furthermore, we can compare partitions by comparing the Bayes optimal equivalence classes
$A_{k}:=c^{-1}(k)=\{x\in\base:c(x)=k\}$ to the estimated equivalence classes $\wh{A}_{L,n,k}:=\wh{c}^{-1}(k)$ by controlling $A_{k}\symdiff \wh{A}_{L,n,k}$, where $A\symdiff B = (A-B)\cup(A-B)$ is the usual symmetric difference of two sets. Specifically, we'd like to show that the difference $A_{k}\symdiff \wh{A}_{L,n,k}$ is small. To this end, define a fattening of $\eqset$ by
\begin{align}
\label{eq:defn:fateqset}
\eqset(t)
:= \bigcup_{i\ne j}\{x\in\base: |\truewgt_{i}\truedenscmp_{i}(x) - \truewgt_{j}\truedenscmp_{j}(x)| \le t\},
\quad t>0.
\end{align}

\noindent
Then of course $\eqset=\eqset(0)$. When the boundaries between classes are sharp, this set will be small, however, if two classes have substantial overlap, then $\eqset(t)$ can be large even if $t$ is small. In the latter case, the equivalence classes $A_{k}$ (and hence the clusters) are less meaningful. The purpose of $\eqset(t)$ is to account for sampling error in the estimated partition.

\begin{thm}
\label{thm:main:part}
Assume that $\lim_{L\to\infty}\lim_{n\to\infty}\estdens_{k}(\estassg)=\truedenscmp_{k}$ uniformly on $X$ and $\upsilon$ is any measure on $\base$. Then there exists a sequence $\tseq\to 0$ such that $\wh{c}(x) = c_{\truemix}(x)$ for all $x\in\base-\eqset(\tseq)$ and
\begin{align}
\label{eq:thm:main:part}
\upsilon\Bigg(\bigcup_{k=1}^{K}A_{k}\symdiff \wh{A}_{L,n,k}\Bigg) 
\le \upsilon(\eqset(\tseq))
\to \upsilon(\eqset).
\end{align}
\end{thm}

\noindent
As in Corollary~\ref{cor:main:conv}, under the same assumptions we may take $L=L_{n}\asymp n^{2/3}/\log^{1/3}n$ in Theorem~\ref{thm:main:part} when $\projset = \mixgaussians_0(A,v)$.

The uniform convergence assumption in Theorem~\ref{thm:main:part} may seem strong, however, recall Corollary~\ref{cor:conv:uniform}, which guarantees uniform convergence whenever $\projset_{L}$ is equicontinuous.
For example, recalling Examples~\ref{ex:param:mix} and~\ref{ex:constrain:mix}, it is straightforward to show the following:

\begin{cor}
\label{cor:gauss:part}
Suppose $\base\subset\R^{d}$, $\projset$ is a compact subset of $\mixgaussians$ and $\upsilon$ is any measure on $\base$. If $\truemix$ is $\projset$-clusterable measure under the Hellinger or variational metric, then there exists a sequence $\tseq\to 0$ such that $\wh{c}(x) = c_{\truemix}(x)$ for all $x\in\base-\eqset(\tseq)$ and
\begin{align}
\label{eq:thm:main:part}
\upsilon\Bigg(\bigcup_{k=1}^{K}A_{k}\symdiff \wh{A}_{L,n,k}\Bigg) 
\le \upsilon(\eqset(\tseq))
\to \upsilon(\eqset).
\end{align}
\end{cor}

We can interpret Theorem~\ref{thm:main:part} as follows: As long as we take $L$ and $n$ large enough and the boundaries between each pair of classes is sharp (in the sense that $\upsilon(\eqset(\tseq))$ is small), the difference between the true Bayes optimal partition and the estimated partition becomes negligible. In fact, it follows trivially from Theorem~\ref{thm:main:part} that $\wh{c}\to c_{\truemix}$ uniformly on $\base-\eqset(t)$ for any fixed $t>0$. Thus, Theorem~\ref{thm:main:part} gives rigourous justification to the approximation heuristic outlined above, and establishes precise conditions under which \emph{nonparametric} clusterings can be learned from data.

\begin{remark}
The sequence $\tseq$ is essentially the rate of convergence of $\estdens_{k}\to \truecmp_{k}$. It is an interesting question to quantify this convergence rate more precisely, which we have left to future work.
\end{remark}

\section{Experiments}
\label{sec:exp}

The theory developed so far suggests an intuitive meta-algorithm for nonparametric clustering. This algorithm can be implemented in just a few lines of code, making it a convenient alternative to more complicated algorithms in the literature. The purpose of this section is merely to illustrate how our theory can be translated into a simple and effective meta-algorithm for nonparametric clustering, which should be understood as a complement to and not a replacement for existing methods that work well in practice.

As in Section~\ref{sec:bayesopt}, we assume we have i.i.d. samples $\rv^{(1)},\ldots,\rv^{(n)}\iid\true=\mixfcn(\truemix)$. Given these samples, we propose the following meta-algorithm:
\begin{enumerate}
\item Estimate an overfitted GMM $\wh{Q}$ with $L\gg K$ components;
\item Define an estimated assignment function $\estassg$ by using single-linkage clustering to group the components of $\wh{Q}$ together;
\item Use this clustering to define $K$ mixture components $\wh{Q}_{k}(\estassg)$;
\item Define a partition on $\base$ by using Bayes' rule, e.g. (\ref{eq:est:dens}-\ref{eq:est:clustmap}).
\end{enumerate}

\noindent
Figure~\ref{fig:mog} has already illustrated two examples where this procedure succeeds in the limit as $n\to\infty$. To further assess the effectiveness of this meta-algorithm in practice, we evaluated its performance on simulated data. In our implementation we used the EM algorithm with regularization and weight clipping to learn the GMM $\wh{Q}$ in step 1, although clearly any algorithm for learning a GMM can be used in this step. The details of these experiments can be found in Appendix~\ref{app:expts}. 

We call the resulting algorithm NPMIX (for \emph{N}on\emph{P}arametric \emph{MIX}ture modeling). To illustrate the basic idea, we first implemented four simple one-dimensional models:
\begin{enumerate}[label=(\roman*)]
\item \textsc{GaussGamma} ($K=4$): A mixture of two Gaussian distributions, one gamma distribution, and a Gaussian mixture.
\item \textsc{Gumbel} ($K=3$): A GMM with three components that has been contaminated with non-Gaussian, Gumbel noise.
\item \textsc{Poly} ($K=2$): A mixture of two polynomials with non-overlapping supports.
\item \textsc{Sobolev} ($K=3$): A mixture of three random nonparametric densities, generated from random expansions of an orthogonal basis for the Sobolev space $H^{1}(\R)$. This is the same example used in Figure~\ref{fig:overview}.
\end{enumerate}

\begin{figure}[!h]
\renewcommand\thesubfigure{\roman{subfigure}}
\centering
\begin{subfigure}[t]{0.5\textwidth}
\centering
\includegraphics[width=\textwidth]{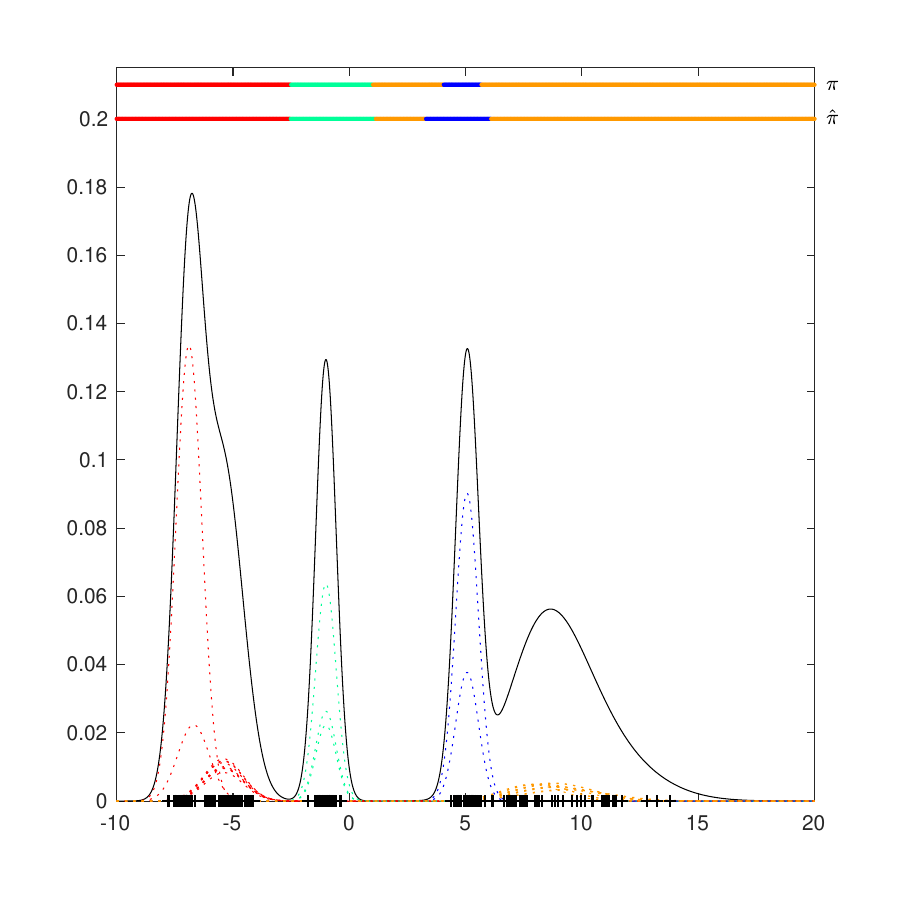}
\caption{\textsc{GaussGamma}.}
\label{fig:1d:gaussgamma}
\end{subfigure}%
\begin{subfigure}[t]{0.5\textwidth}
\centering
\includegraphics[width=\textwidth]{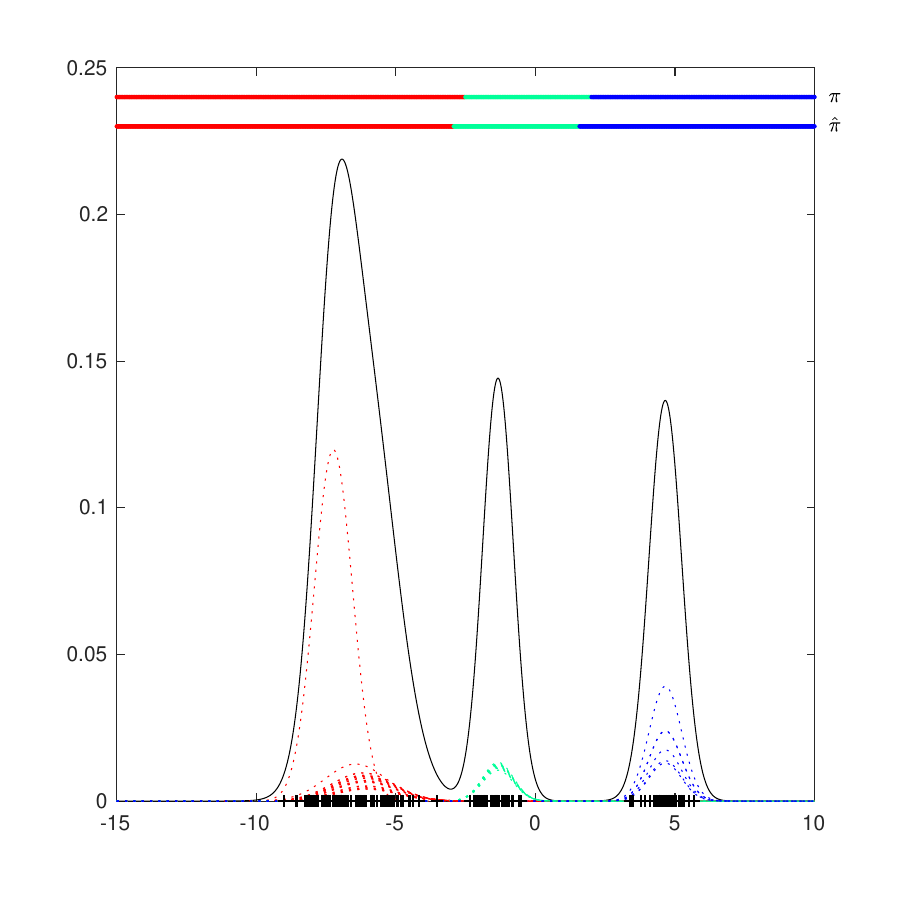}
\caption{\textsc{Gumbel}.}
\label{fig:1d:gumbel}
\end{subfigure}%
\\
\begin{subfigure}[t]{0.5\textwidth}
\centering
\includegraphics[width=\textwidth]{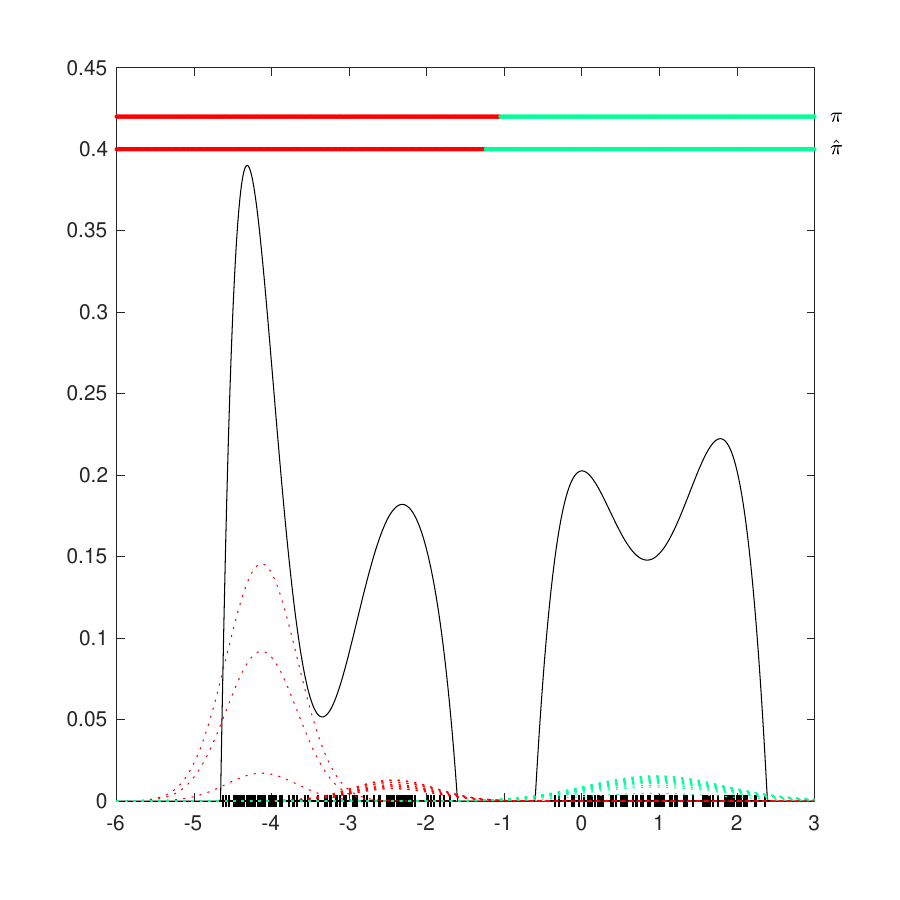}
\caption{\textsc{Poly}.}
\label{fig:1d:poly}
\end{subfigure}%
\begin{subfigure}[t]{0.5\textwidth}
\centering
\includegraphics[width=\textwidth]{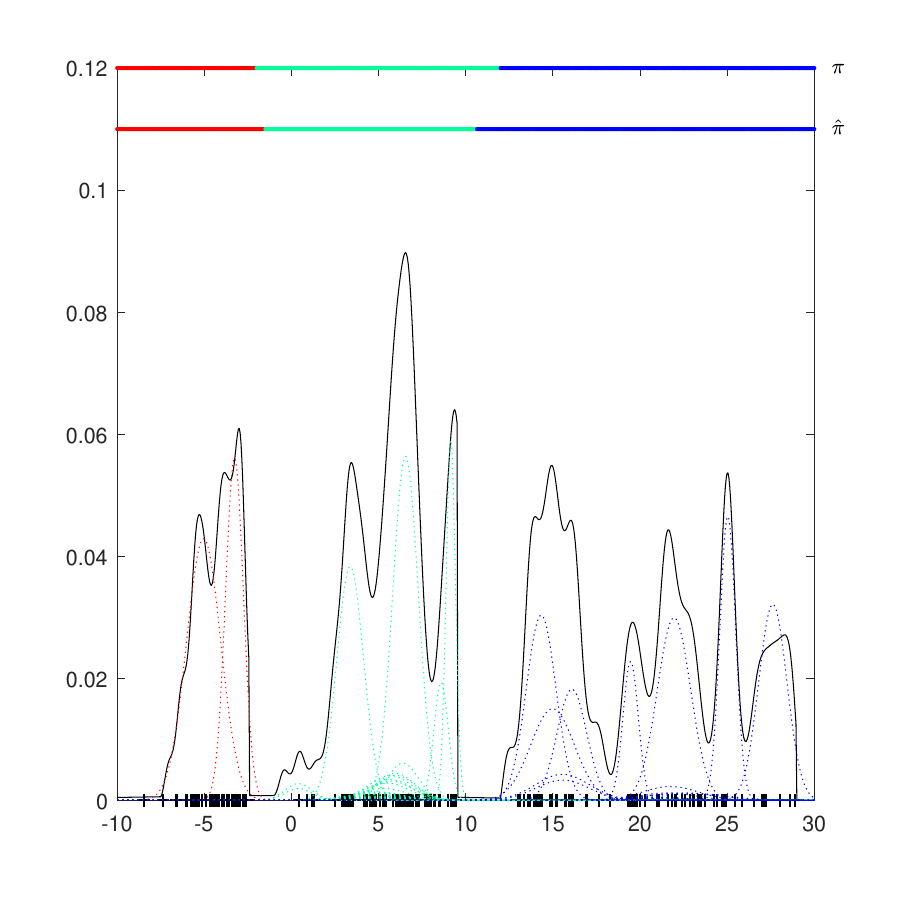}
\caption{\textsc{Sobolev}.}
\label{fig:1d:sobolev}
\end{subfigure}%
\caption{Examples (i)-(iv) of one-dimensional mixture models. The original mixture density is depicted as a solid black line, with the overfitted Gaussian mixture components as dotted lines, coloured according to the cluster they are assigned to. The true Bayes optimal partition $\pt$ and the estimated partition $\wh{\pt}$ are depicted by the horizontal lines at the top, and the raw data are plotted on the $x$-axis for reference.}
\label{fig:1d}
\end{figure}

\noindent
The results are shown in Figure~\ref{fig:1d}. These examples illustrate the basic idea behind the algorithm: Given samples, overfitted mixture components (depicted by dotted lines in Figure~\ref{fig:1d}) are used to approximate the global nonparametric mixture distribution (solid black line). Each of these components is then clustered, with the resulting partition of $\base=\R$ depicted alongside the true Bayes optimal partition. In each case, cutting the cluster tree to produce $K$ components provides sensible and meaningful approximations to the true partitions.

\begin{figure}[t]
\centering
\includegraphics[width=0.64\textwidth]{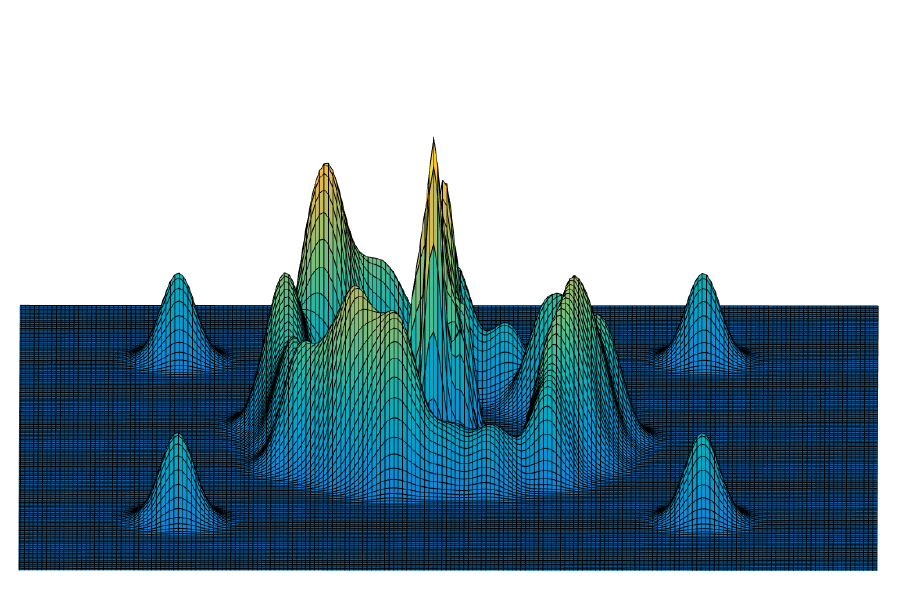}\\
\includegraphics[width=0.3\textwidth]{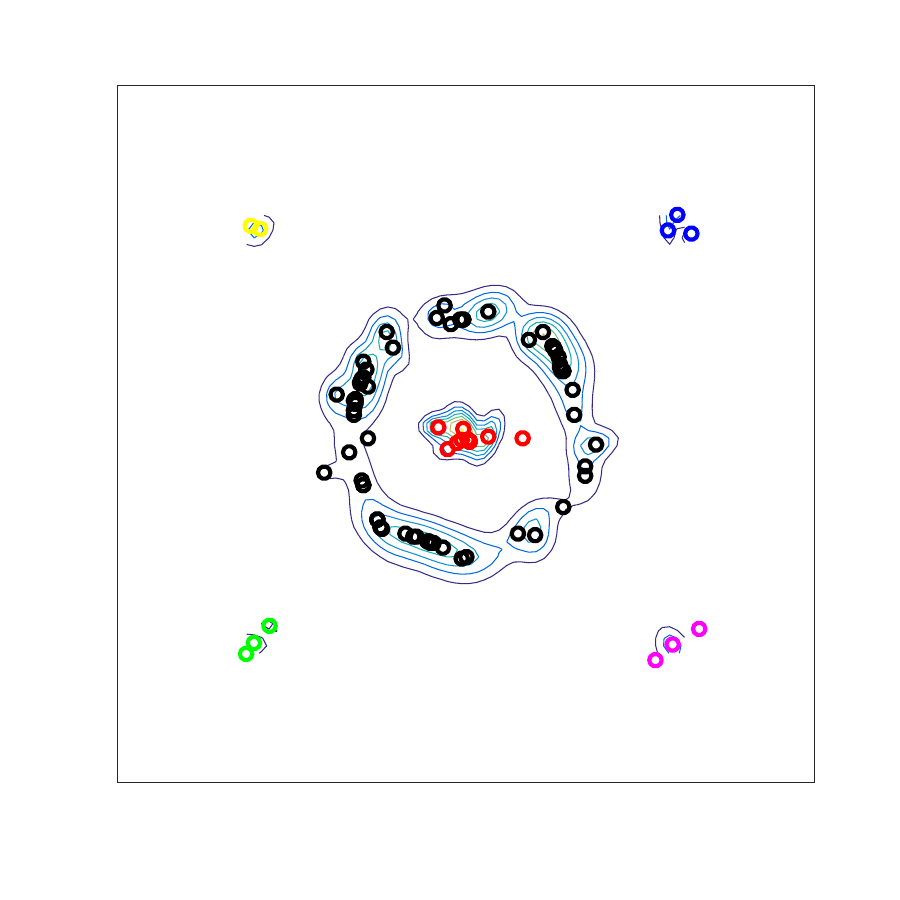}
\includegraphics[width=0.3\textwidth]{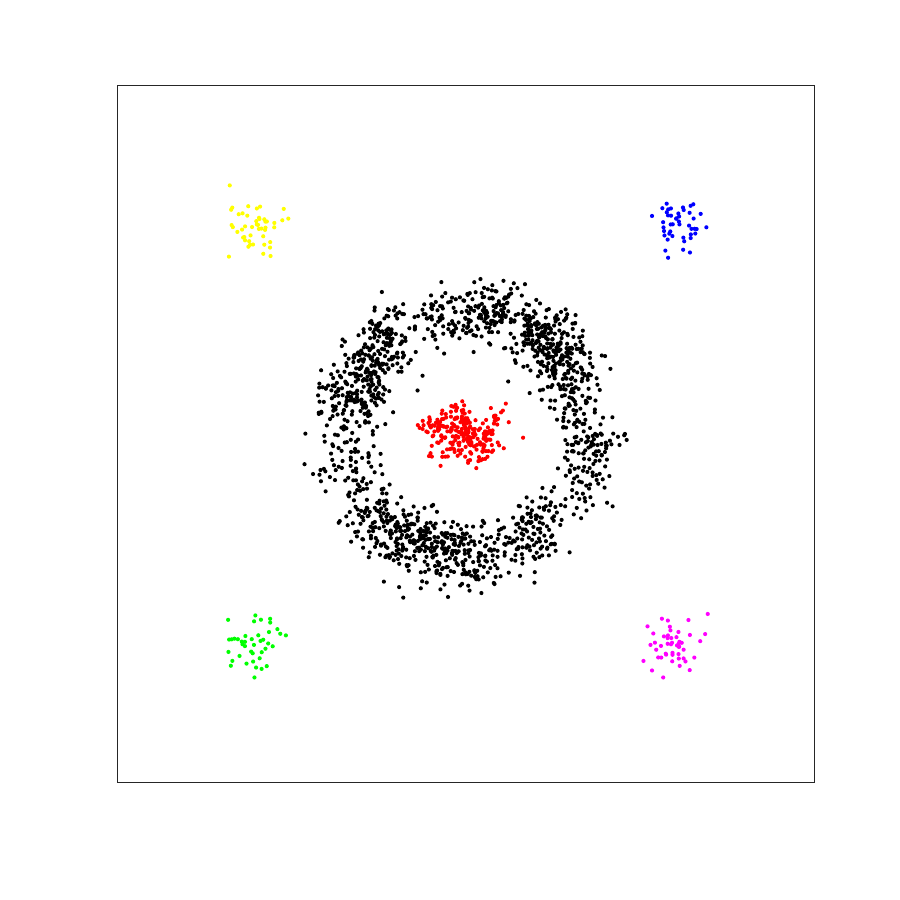}
\includegraphics[width=0.3\textwidth]{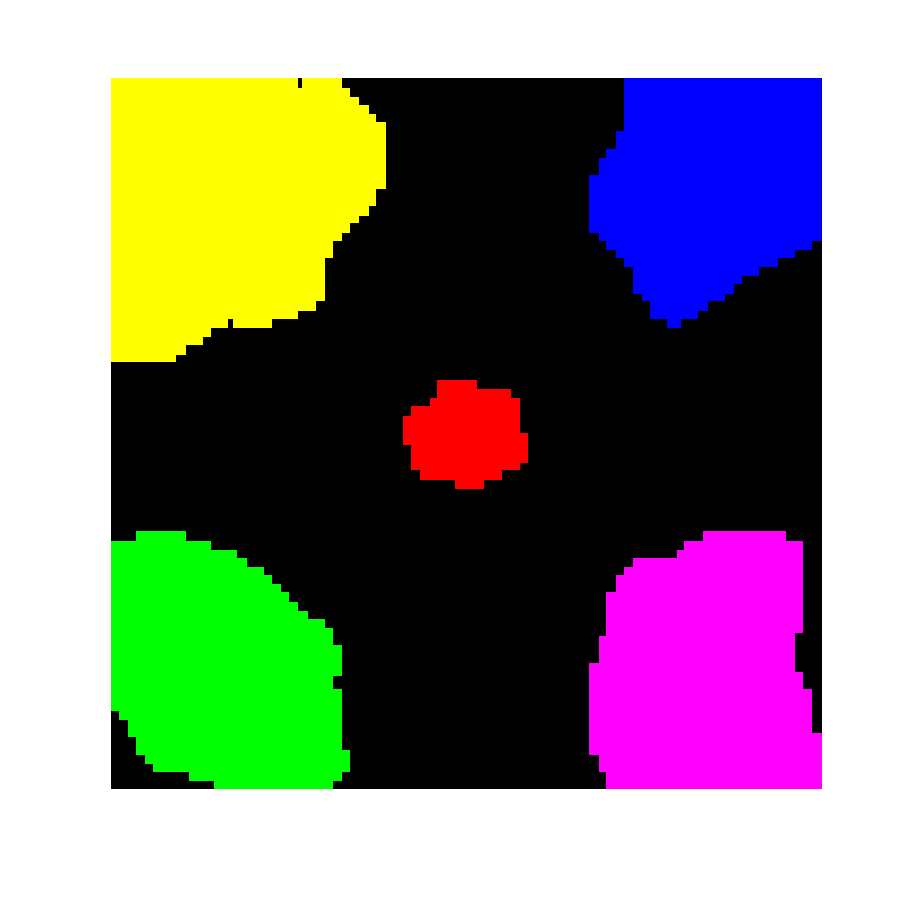}
\caption{Example of a successful clustering on the \textsc{Target} mixture model using NPMIX. (Top) Density plot of the original mixture density. (Left) Contour plot of overfitted Gaussian mixture approximation, centers marked with $\circ$'s. (Middle) Original data colour coded by the approximate Bayes optimal partition. (Right) Estimated Bayes optimal partition, visualized as the input space $\base$ colour-coded by estimated cluster membership.}
\label{fig:target}
\end{figure}

To further validate the proposed algorithm, we implemented the following two-dimensional mixture models and compared the cluster accuracy to existing clustering algorithms on simulated data:
\begin{enumerate}[label=(\roman*)]
\setcounter{enumi}{4}
\item \textsc{Moons} ($K=2$): A version of the classical \textsc{moons} dataset in two-dimensions. This model exhibits a classical failure case of spectral clustering, which is known to have difficulties when clusters are unbalanced (i.e. $\lambda_{1}\ne\lambda_{2}$). For this reason, we ran experiments with both balanced and unbalanced clusters.
\item \textsc{Target} ($K=6$): A GMM derived from the \textsc{target} dataset (Figure~\ref{fig:target}). The GMM has 143 components that are clustered into 6 groups based on the original \textsc{Target} dataset from \citep{ultsch2005}.
\end{enumerate}

\begin{figure}[t]
\centering
\includegraphics[width=0.3\textwidth]{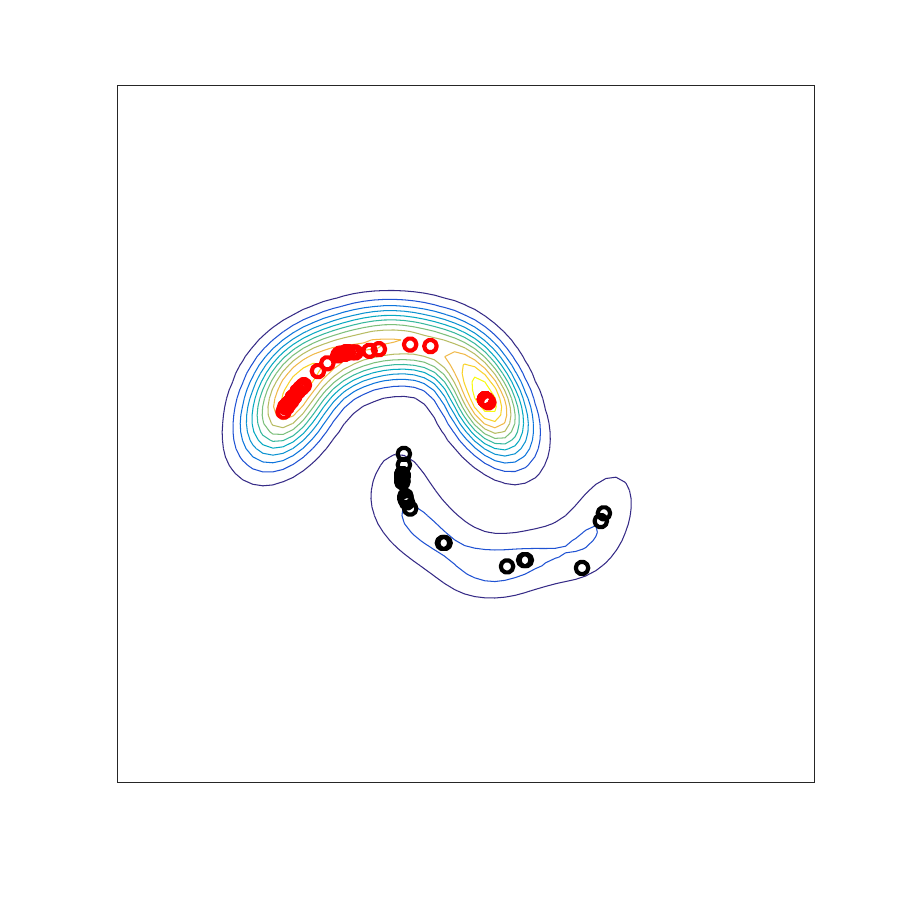}
\includegraphics[width=0.3\textwidth]{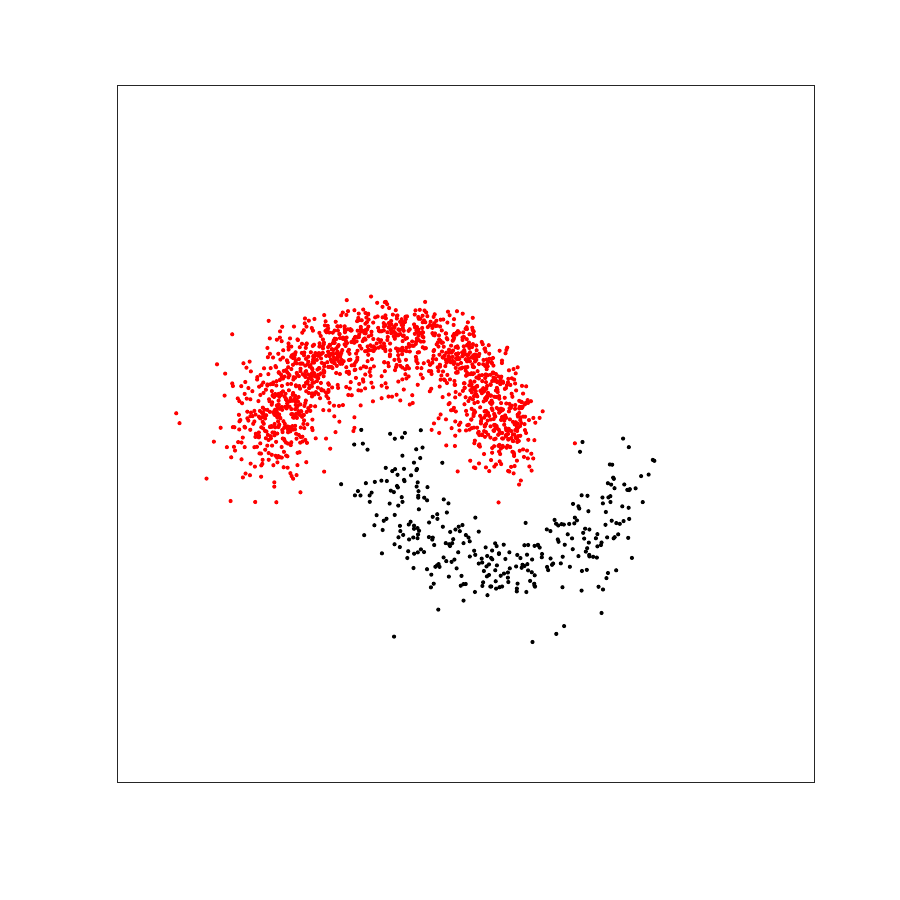}
\includegraphics[width=0.3\textwidth]{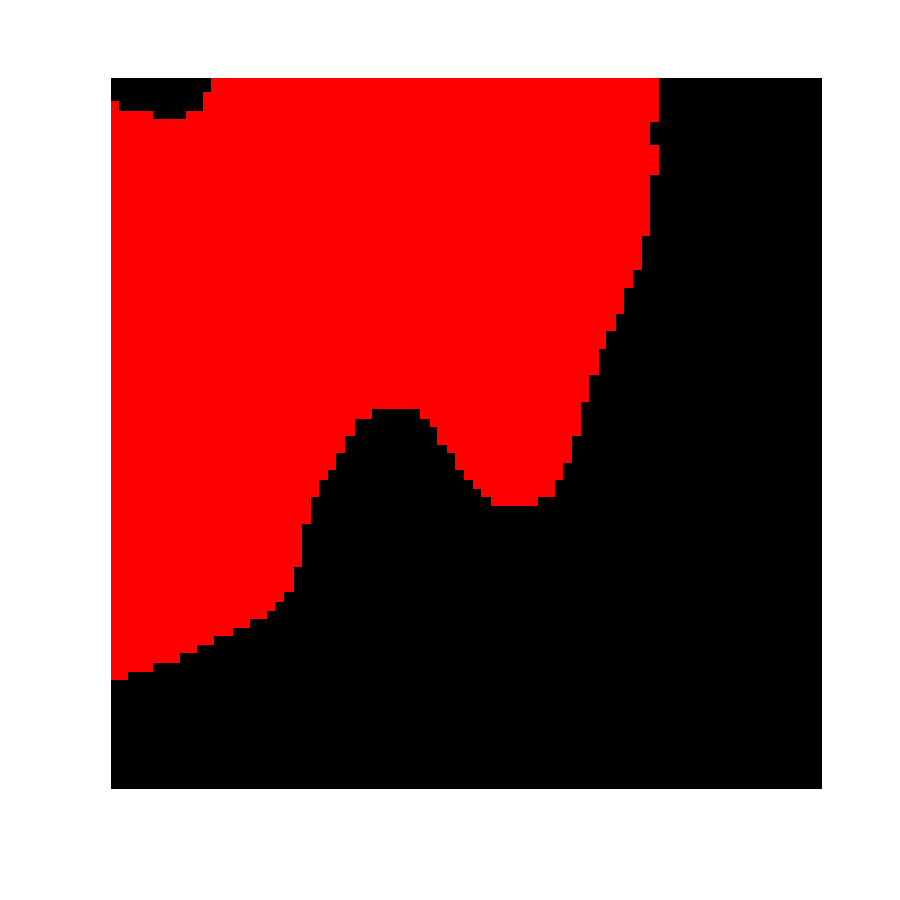}
\caption{Example of a successful clustering on the unbalanced \textsc{Moons} mixture model using NPMIX. (Left) Contour plot of overfitted Gaussian mixture approximation, centers marked with $\circ$'s. (Middle) Original data colour coded by the approximate Bayes optimal partition. (Right) Estimated Bayes optimal partition, visualized as the input space $\base$ colour-coded by estimated cluster membership.}
\label{fig:moons:unbalanced}
\end{figure}

\begin{figure}[t]
\centering
\includegraphics[width=0.3\textwidth]{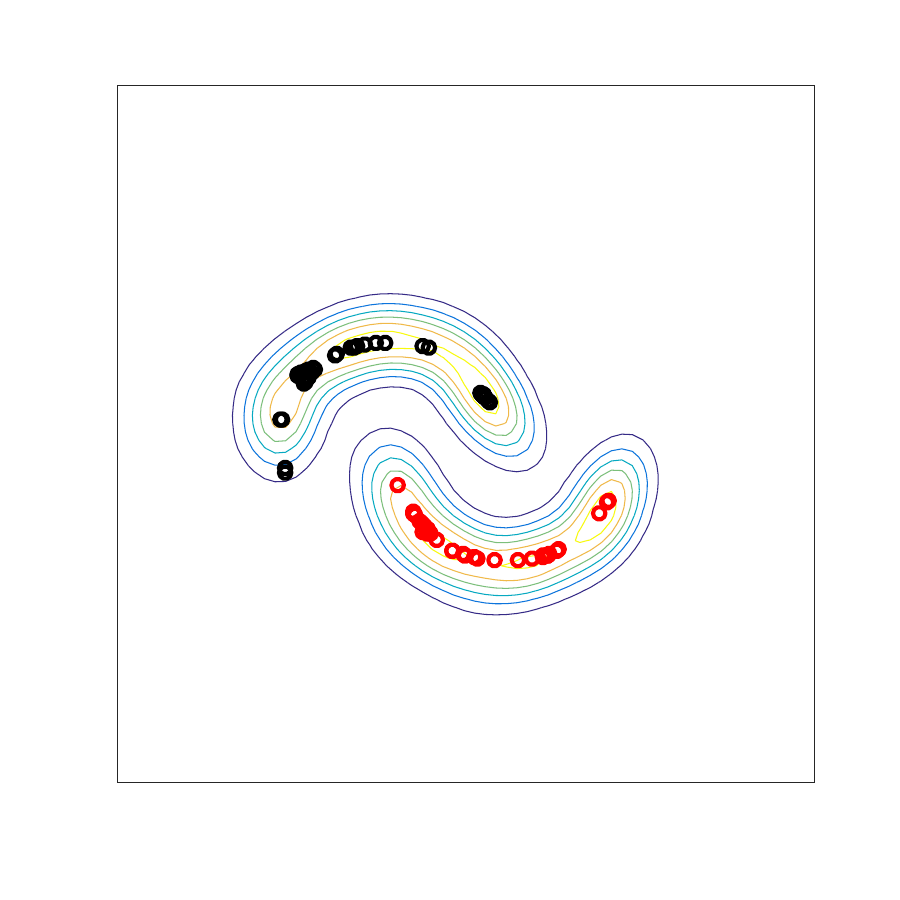}
\includegraphics[width=0.3\textwidth]{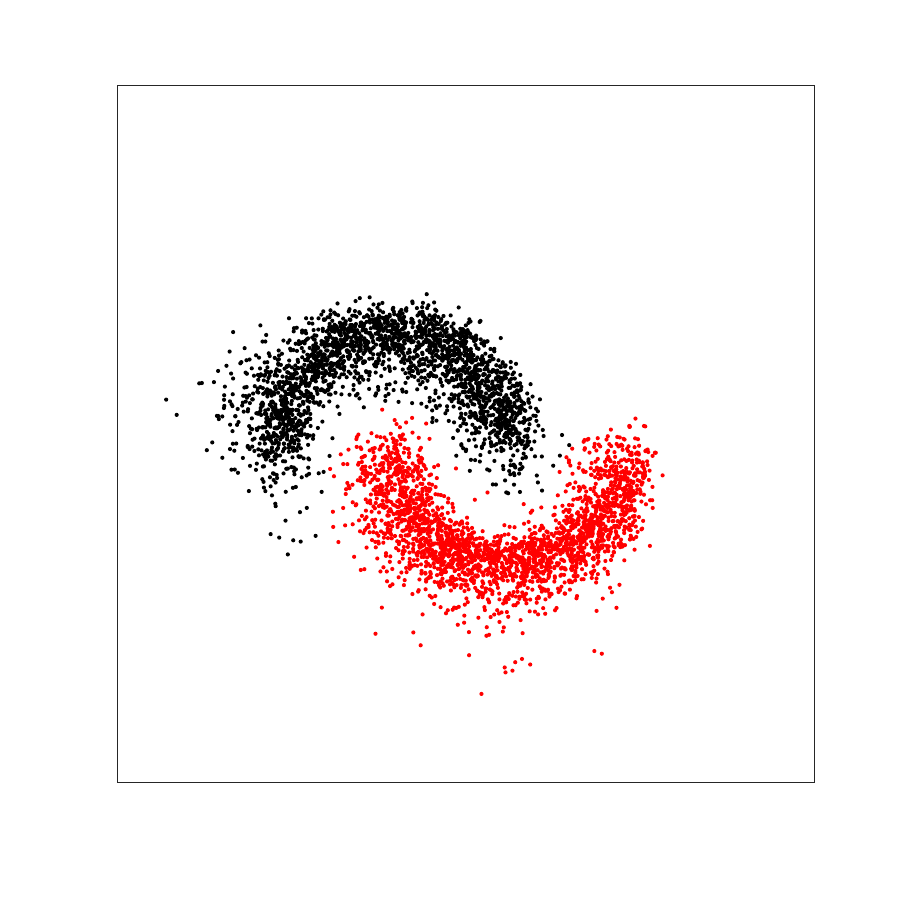}
\includegraphics[width=0.3\textwidth]{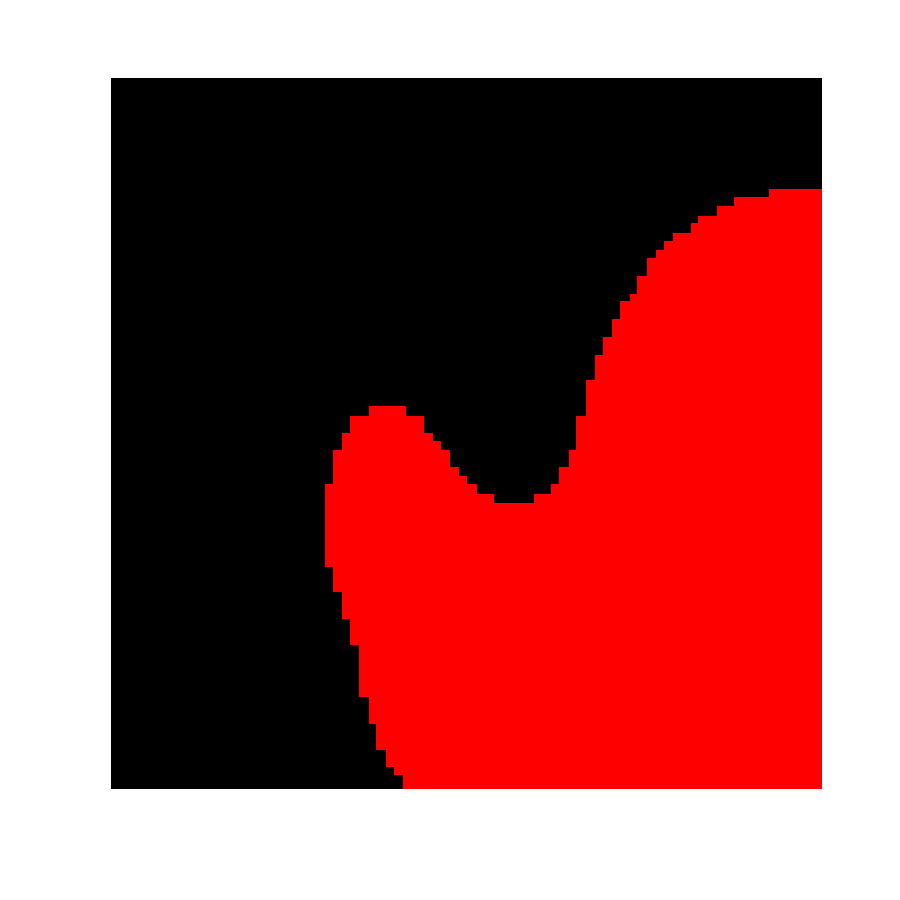}
\caption{Example of a successful clustering on the balanced \textsc{Moons} mixture model using NPMIX. (Left) Contour plot of overfitted Gaussian mixture approximation, centers marked with $\circ$'s. (Middle) Original data colour coded by the approximate Bayes optimal partition. (Right) Estimated Bayes optimal partition, visualized as the input space $\base$ colour-coded by estimated cluster membership.}
\label{fig:moons:balanced}
\end{figure}

\noindent
Visualizations of the results for our method are shown in Figures~\ref{fig:moons:unbalanced},~\ref{fig:moons:balanced}, and~\ref{fig:target}. One of the advantages of our method is the construction of an explicit partition of the entire input space (in this case, $\base=\R^{2}$), which is depicted in all three figures. Mixture models are known to occasionally lead to unintuitive cluster assignments in the tails, which we observed with the unbalanced \textsc{Moons} model. This is likely an artifact of the sensitivity of the EM algorithm, and can likely be corrected by using a more robust mixture model estimator in the first step.

We compared NPMIX against four well-known benchmark algorithms: (i) $K$-means, (ii) Spectral clustering, (iii) Single-linkage hierarchical clustering, and (iv) A Gaussian mixture model (GMM) with $K$ components. We only considered methods that classify every sample in a dataset (this precludes, e.g. density-based clustering).
Moreover, of these four algorithms, only $K$-means and GMM provide a partition of the entire input space $\base$, which allows for new samples to be classified without re-running the algorithm. 
All of the methods (including NPMIX) require the specification of the number of clusters $K$, which was set to the correct number according to the model. 
In each experiment, we sampled random data from each model and then used each clustering algorithm to classify each sample. To assess cluster accuracy, we computed the adjusted RAND index (ARI) for the clustering returned by each method. ARI is a standard permutation-invariant measure of cluster accuracy in the literature. 

The results are shown in Table~\ref{tab:results}. On the unbalanced \textsc{Moons} data, NPMIX clearly outperformed each of the four existing methods. On balanced data, $K$-means, spectral clustering, and GMM improved significantly, with spectral clustering performing quite well on average. All four algorithms were still outperformed by NPMIX. On \textsc{Target}, the results were more interesting: Both single-linkage and spectral clustering perform very well on this dataset. NPMIX shows more variance in its performance, as indicated by the high median (0.998) and lower mean (0.696). On 57/100 runs, the ARI for NPMIX was $>0.99$, and on the rest the ARI was $<0.6$. This is likely caused by sensitivity to outliers in the \textsc{Target} model, and we expect that this can be corrected by using a more robust algorithm (e.g. instead of the vanilla EM algorithm). As our motivations are mainly theoretical, we leave more detailed fine-tuning of this algorithm and thorough side-by-side comparisons to future work. For example, by using the learned mixture density to remove ``background samples'' (e.g. as in density-based clustering), this algorithm can be trivially improved.

\begin{table}[t]
\begin{center}
\caption{Average and median adjusted RAND index (ARI) for $N=100$ simulations of three different nonparametric mixture models.}
\begin{tabular}{lcccc}
\toprule
\textsc{Moons (unbalanced)} & Mean ARI & Median ARI & st. dev. \\
\midrule
NPMIX & \textbf{0.727} & \textbf{0.955} & 0.284 \\
$K$-means & 0.126 & 0.124 & 0.016 \\
Spectral & 0.197 & 0.122 & 0.232 \\
Single-linkage & 0.001 & 0.001 & 0.002 \\
GMM & 0.079 & 0.078 & $<10^{-3}$ \\
\toprule
\textsc{Moons (balanced)} & Mean ARI & Median ARI & st. dev. \\
\midrule
NPMIX & \textbf{0.934} & \textbf{0.972} & 0.188 \\
$K$-means & 0.502 & 0.503 & 0.021 \\
Spectral & 0.909 & 0.910 & 0.013 \\
Single-linkage & $<10^{-6}$ & $<10^{-6}$ & $<10^{-6}$ \\
GMM & 0.782 & 0.783 & $<10^{-3}$ \\
\toprule
\textsc{Target} & Mean ARI & Median ARI & st. dev. \\
\midrule
NPMIX & 0.696 & 0.998 & 0.354 \\
$K$-means & 0.081 & 0.072 & 0.034 \\
Spectral & \textbf{0.967} & 0.975 & 0.077 \\
Single-linkage & 0.824 & \textbf{1.000} & 0.222 \\
GMM & 0.126 & 0.124 & 0.002 \\
\bottomrule
\end{tabular}
\label{tab:results}
\end{center}
\end{table}

\section{Discussion}
\label{sec:conc}

We have established a new set of identifiability results for nonparametric mixtures that rely on the notion of \emph{clusterability}. In particular, our results allow for an arbitrary number of components and for each component to take on essentially any shape. The key assumption is separation between the components, which allows simple clustering algorithms such as hierarchical clustering to recover individual mixture components from an overfitted mixture density estimator. Furthermore, we established conditions under which identified mixtures and their partitions can be consistently estimated from data. We also discussed applications to data clustering, including a nonparametric notion of the Bayes optimal partition and an intuitive meta-algorithm for nonparametric clustering. 

The assumption that the number of components $K$ is known is of course restrictive in practice, however, this assumption can be substantially relaxed as follows: If $K$ is unknown, simply test whether or not there exists a $K$ such that the separation criterion \eqref{eq:condn:sep} holds. If such a $K$ exists and is unique, then the resulting $K$-mixture is identifiable. In practice, however, there may be more than one value of $K$ for which \eqref{eq:condn:sep} holds. Furthermore, if $\truemix$ is identifiable for some $K$, it may not be the case that $\truemix$ is identifiable for $K'<K$ owing to the separation criterion \eqref{eq:defn:sep:limsup} (cf. \eqref{eq:defn:sep}). Of course, such an exhaustive search may not be practical, in which case it would be interesting to study efficient algorithms for finding such a $K$.

As pointed out by a reviewer, there is a connection between the NPMIX algorithm introduced in Section~\ref{sec:exp} and kernel density estimation (KDE). Indeed, by choosing $L=n$, the overfitted mixture model learned in step 1 is similar to a kernel density estimate with a Gaussian kernel, although not exactly the same since KDE fixes the weights, centers, and bandwidth of each kernel unless more sophisticated adaptive bandwidth selection strategies are used. By contrast, a GMM allows these parameters to be learned from the data. Thus, in the limiting case $L=n$, NPMIX is similar to single-linkage clustering applied to a new metric defined via the Wasserstein distance between the $n$ kernels, where this new metric depends crucially on the choice of bandwidth. An important difference in practice is that by taking $L<n$, the NPMIX algorithm denoises the data in the first step, making it less sensitive to outliers. For example, \citet{priebe1994adaptive} points out that approximately $L=30$ Gaussian components suffice to approximate a log-normal density with $n=10,000$ samples; see also Corollary~\ref{cor:main:conv}.
Exploring this connection more deeply is an interesting direction for future work.

It would also be interesting to study convergence rates for the proposed estimators. In particular, there are two important quantities of interest in deriving these rates: The sample size $n$ and the number of overfitted components $L$. Interestingly, it was only recently that the minimax rate of estimation for \emph{parametric} mixtures was correctly determined \citep{heinrich2018strong}, which is $n^{1/(4(s-s_{0})+2)}$ in the $L_{1}$-Wasserstein metric, where $s_{0}$ is the true number of mixture components and $s$ is the number used in estimation. See also \citep{chen1995,nguyen2013,ho2016,ho2016singularity}. In the general case, this is also related to problems in agnostic learning \citep{li2015agnostic}. In our nonparametric setting, we expect these rates to depend on both $L$ and $n$. Furthermore, it is necessary to control the distance between the $\probm$-projection $\proj$ and $\true$, which depends on the choice of $L$ alone. This latter problem will almost certainly require imposing additional regularity conditions on $\true$, e.g. as in \citep{genovese2000,ghosal2001entropies}. 

Finally, it would be of significant interest to apply existing clustering theory to find new conditions that guarantee clusterability in the same way that Proposition~\ref{prop:pop:cluster} shows that separability is sufficient for single-linkage clustering. We have already noted that the separation constant $4\separ(\truemix)$ can be reduced. Furthermore, in simulations we have observed that complete-linkage is often sufficient when working with the proposed NPMIX algorithm. But under what precise conditions on $\true$ is complete-linkage sufficient? By applying known results from the clustering literature, it may be possible to extend our results to prove deeper identifiability theorems for nonparametric mixtures.

\bibliography{genbib} %
\bibliographystyle{abbrvnat}
\label{bib:bib}

\appendix

\section{Proofs}
\label{app:proofs}

Throughout these appendices, we denote the Hellinger metric on $\probs(\base)$ by $\dH$ and similarly the variational (a.k.a. total variation) distance by $\dTV$.

\subsection{Proof of Lemma~\ref{lem:reg:dense}}
\label{app:proofs:lem:reg:dense}

Fix $\eps>0$ and let $\gaussmix_{k,\eps}\in\mixing_{0}(\projset)$ be such that $\dTV(\mixfcn(\gaussmix_{k,\eps}),\truecmp_{k})<\eps/K$. Define $E:=\cup_{k=1}^{K}\supp(\gaussmix_{k,\eps})$ and let $g_{\eps,\ell}\sim\normalN(\mu_{\eps,\ell},\Sigma_{\eps,\ell})$ ($\ell=1,\ldots,L$) be an enumeration of $E$. If each $g_{\eps,\ell}$ is distinct, then we are done. Suppose to the contrary that $g_{\eps,\ell}=g_{\eps,\tick{\ell}}$ for some $\tick{\ell}\ne\ell$. Let $t>0$ be so small that $B(\mu_{\eps,\ell},t)\cap E=\{\mu_{\eps,\ell}\}$  for all $\ell$. Then we can find $\delta_{\ell}\ne\delta_{\tick{\ell}}\in B(\mu_{\eps,\ell},t)$ such that $\mu_{\eps,\ell}+\delta_{\ell}\ne\mu_{\eps,\tick{\ell}}+\delta_{\tick{\ell}}$. If there are multiple such atoms, repeat this process until there are no shared atoms. Let $\delta=(\delta_{1},\ldots,\delta_{L})$ and define $\gaussmix_{k,\eps,\delta}$ to be equal to $\gaussmix_{k,\eps}$ except that $g_{\eps,\ell}$ is replaced with $\normalN(\mu_{\eps,\ell}+\delta_{\ell},\Sigma_{\eps,\ell})$. Evidently, it is clear from this construction that $\dTV(\mixfcn(\gaussmix_{k,\eps,\delta}),\truecmp_{k})\to0$ as $\eps,t\to0$. Defining $\gaussmix_{\eps}:=\sum_{k=1}^{K}\truewgt_{k}\gaussmix_{\eps,k,\delta}$, we have $\mixm(\gaussmix_{\eps},\truemix)\to0$.

It remains to show that $\gaussmix_{\eps}$ is $\mixgaussians$-regular. It is clear that the $\dTV$-projection exists and is unique for $L\ge |E|$ (i.e. since $\true\in\mix(\mixgaussians_{L})$), which establishes Definition~\ref{defn:regular}\ref{defn:regular:a}. To see Definition~\ref{defn:regular}\ref{defn:regular:b}, one may check that the assignment sequence given by $\assg_{L}(\ell)=k\iff g_{\eps,\ell}\in\supp(\gaussmix_{k,\eps,\delta})$ works.
\qed

\subsection{Some metric inequalities}
\label{app:proofs:metric}

Throughout this section, we assume that $L$ and $\assg$ are fixed. Thus we will often suppress the dependence on $L$ and $\assg$, writing $\projmix_{k}=\projmix_{k}(\assg)$, $\proj_{k}=\proj_{k}(\assg)$, and so on. We will also write $\projcmp_{\ell}\in \proj_{k}\iff\assg(\ell)=k\iff\projcmp_{\ell}\in \supp(\projmix_{k})$.

Next, we introduce some new notation. Define
\begin{align}
\label{eq:defn:gensep}
\gensep(t)
:= \sup_{k}\Hdiam(\projmix_{k}) + t.
\end{align}

\noindent 
and
\begin{align}
\label{eq:defn:eps:delt}
\eps 
=\eps_{L,n}
:= \sup_{\ell}\probm(\estcmp_{\ell},\projcmp_{\ell}),
\quad
\delta
=\delta_{L}
:= \sup_{k}\probm(\proj_{k},\truecmp_{k}).
\end{align}

\noindent
Recall the following assumptions, which are restated here for reference in the proofs:
\begin{enumerate}[label=(B\arabic*)]
\item\label{lem:metric:ineq:A1} $\probm(\truecmp_{i},\truecmp_{j})\ge 4\gensep(\delta)$ for $i\ne j$; %
\item\label{lem:metric:ineq:A2} There is an assignment map such that $\probm(\proj_{k},\truecmp_{k})<\delta$ for each $k$;
\item\label{lem:metric:ineq:A3} $\estcmp_{\ell}\in B(\projcmp_{\ell},\eps)$ for all $\ell$.
\end{enumerate}

\noindent
\ref{lem:metric:ineq:A1} is \eqref{eq:condn:sep}, \ref{lem:metric:ineq:A2} is Definition~\ref{defn:regular}\ref{defn:regular:b}, and \ref{lem:metric:ineq:A3} is \eqref{eq:defn:eps:delt}.

\begin{lemma}
\label{lem:metric:ineq:true}
Let $\delta>0$ and $\eps>0$ be arbitrary. Under \ref{lem:metric:ineq:A1}-\ref{lem:metric:ineq:A3}, the following are true:
\begin{enumerate}[label=(M\arabic*)]
\item\label{lem:metric:ineq:M1} $\probm(\projcmp_{\ell},\truecmp_{k})\le\gensep(\delta)$ if $\projcmp_{\ell}\in \proj_{k}$;
\item\label{lem:metric:ineq:M2} $\probm(\estcmp_{\ell},\truecmp_{k})\le\gensep(\delta+\eps)$ if $\projcmp_{\ell}\in \proj_{k}$;
\item\label{lem:metric:ineq:M3} $\probm(\projcmp_{\ell},\truecmp_{i})\ge 3\gensep(\delta)$ if $\projcmp_{\ell}\notin \proj_{i}$;
\item\label{lem:metric:ineq:M4} $\probm(\estcmp_{\ell},\truecmp_{i})\ge 3\gensep(\delta)-\eps$ if $\projcmp_{\ell}\notin \proj_{i}$;
\end{enumerate}
\end{lemma}

\begin{proof}
We prove each claim in order below.

\medskip\noindent
\ref{lem:metric:ineq:M1}: Using \ref{lem:metric:ineq:A2}, we have 
\begin{align*}
\probm(\truecmp_{k},\projcmp_{\ell}) 
&\le \probm(\truecmp_{k},\proj_{k}) + \probm(\proj_{k},\projcmp_{\ell}) \\
&\le \delta + \Hdiam(\projmix_{k}) \\
&\le \gensep(\delta).
\end{align*}

\noindent

\medskip\noindent
\ref{lem:metric:ineq:M2}: Invoking both \ref{lem:metric:ineq:A3} and \ref{lem:metric:ineq:M1}, we have $\probm(\estcmp_{\ell},\truecmp_{k}) \le \probm(\estcmp_{\ell},\projcmp_{\ell}) + \probm(\projcmp_{\ell},\truecmp_{k})\le\eps+\gensep(\delta) = \gensep(\delta+\eps)$.

\medskip\noindent
\ref{lem:metric:ineq:M3}: Let $k$ be the index such that $\projcmp_{\ell}\in \proj_{k}$. Then, via the reverse triangle inequality,
\begin{align*}
\gensep(\delta)
\overset{\text{\ref{lem:metric:ineq:M1}}}{\ge} \probm(\truecmp_{k},\projcmp_{\ell})
\ge \probm(\truecmp_{k},\truecmp_{i}) - \probm(\truecmp_{i},\projcmp_{\ell})
\overset{\text{\ref{lem:metric:ineq:A1}}}{\ge} 4\gensep(\delta) - \probm(\truecmp_{i},\projcmp_{\ell}),
\end{align*}

\noindent
and the desired result follows after re-arrangement.

\medskip\noindent
\ref{lem:metric:ineq:M4}: Same as the proof of \ref{lem:metric:ineq:M3}, except invoke \ref{lem:metric:ineq:M2} in place of \ref{lem:metric:ineq:M1}.
\end{proof}

The previous lemma bounded the distance between $\projcmp_{\ell}$, $\estcmp_{\ell}$ and the true nonparametric measure $\truecmp_{k}$. The next lemma leverages the previous one to bound the distances between pairs of overfitted components and their estimates.

\begin{lemma}
\label{lem:metric:ineq:pair}
Let $\delta>0$ and $\eps>0$ be arbitrary. Under \ref{lem:metric:ineq:A1}-\ref{lem:metric:ineq:A3}, the following are true:
\begin{enumerate}[label=(M\arabic*)]
\setcounter{enumi}{4}
\item\label{lem:metric:ineq:M5} $\probm(\projcmp_{\ell},\projcmp_{\ell'})\le\gensep(\delta)$ if $\projcmp_{\ell},\projcmp_{\ell'}\in \proj_{k}$;
\item\label{lem:metric:ineq:M6} $\probm(\estcmp_{\ell},\estcmp_{\ell'})\le \gensep(2\eps)$ if $\projcmp_{\ell},\projcmp_{\ell'}\in \proj_{k}$;
\item\label{lem:metric:ineq:M7} If $\projcmp_{\ell}\in \proj_{k}$ and $\projcmp_{\ell'}\in \proj_{k'}$, $k\ne k'$, then $\probm(\projcmp_{\ell}, \projcmp_{\ell'}) \ge 2\gensep(\delta)$.
\item\label{lem:metric:ineq:M8} If $\projcmp_{\ell}\in \proj_{k}$ and $\projcmp_{\ell'}\in \proj_{k'}$, $k\ne k'$, then $\probm(\estcmp_{\ell}, \estcmp_{\ell'}) \ge 2\gensep(\delta) - \eps$.
\end{enumerate}
\end{lemma}

\begin{proof}
We prove each claim in order below.

\medskip\noindent
\ref{lem:metric:ineq:M5}: Since $\projcmp_{\ell},\projcmp_{\ell'}\in \proj_{k}$ it follows by definition that $\projcmp_{\ell},\projcmp_{\ell'}\in \supp(\projmix_{k})$. But then $\probm(\projcmp_{i},\projcmp_{j})\le\Hdiam(\projmix_{k})\le\gensep(\delta)$.

\medskip\noindent
\ref{lem:metric:ineq:M6}: As in the proof of \ref{lem:metric:ineq:M5}, we have $\probm(q_{\ell},q_{\ell'})\le \Hdiam(\projmix_{k})$ and hence
\begin{align*}
\probm(\estcmp_{\ell}, \estcmp_{\ell'})
\le \probm(\estcmp_{\ell}, q_{\ell}) + \probm(q_{\ell},q_{\ell'}) + \probm(q_{\ell'},\estcmp_{\ell'})
\le 2\eps + \Hdiam(\projmix_{k})
\le \gensep(2\eps).
\end{align*}

\medskip\noindent
\ref{lem:metric:ineq:M7}: We have 
\begin{align*}
\probm(\projcmp_{\ell}, \projcmp_{\ell'})
\ge \probm(\projcmp_{\ell}, \truecmp_{k'}) - \probm(\truecmp_{k'}, \projcmp_{\ell'})
\overset{\text{\ref{lem:metric:ineq:M1}}}{\ge} \probm(\projcmp_{\ell}, \truecmp_{k'}) - \gensep(\delta)
\end{align*}

\noindent
and similarly $\probm(\projcmp_{\ell}, \projcmp_{\ell'})\ge \probm(\projcmp_{\ell'}, \truecmp_{k}) - \gensep(\delta)$. Adding these inequalities,
\begin{align*}
2\probm(\projcmp_{\ell}, \projcmp_{\ell'})
&\ge \probm(\projcmp_{\ell}, \truecmp_{k'}) + \probm(\projcmp_{\ell'}, \truecmp_{k}) - 2\gensep(\delta) \\
\implies
\probm(\projcmp_{\ell}, \truecmp_{k'}) + \probm(\projcmp_{\ell'}, \truecmp_{k})
&\le 2\probm(\projcmp_{\ell}, \projcmp_{\ell'}) + 2\gensep(\delta).
\end{align*}

\noindent
Invoking \ref{lem:metric:ineq:M3} on the left, we have 
\begin{align*}
6\gensep(\delta)
\le \probm(\projcmp_{\ell}, \truecmp_{k'}) + \probm(\projcmp_{\ell'}, \truecmp_{k})
&\le 2\probm(\projcmp_{\ell}, \projcmp_{\ell'}) + 2\gensep(\delta)
\iff
\probm(\projcmp_{\ell}, \projcmp_{\ell'}) 
\ge 2\gensep(\delta).
\end{align*}

\medskip\noindent
\ref{lem:metric:ineq:M8}: In the proof of \ref{lem:metric:ineq:M7}, replace \ref{lem:metric:ineq:M1} with \ref{lem:metric:ineq:M2} and \ref{lem:metric:ineq:M3} with \ref{lem:metric:ineq:M4}.
\end{proof}

\subsection{Proof of Proposition~\ref{prop:pop:cluster}}
\label{app:proofs:prop:pop:cluster}

By definition, $\assg(i)=\assg(j)\iff\projcmp_{i},\projcmp_{j}\in \supp(\projmix_{k})$. Note also that $\separ(\projmix(\assg))=\gensep(\delta)$ (cf. \eqref{eq:defn:sep}, \eqref{eq:defn:gensep}) for $\delta:= \sup_{k}\probm(\proj_{k},\truecmp_{k})$. Now, if $\projcmp_{i},\projcmp_{j}\in \supp(\projmix_{k})$, then $\probm(\projcmp_{i},\projcmp_{j})\le\gensep(\delta)=\separ(\projmix(\assg))$ by \ref{lem:metric:ineq:M5} of Lemma~\ref{lem:metric:ineq:pair}. Conversely, suppose $\projcmp_{i}\in \supp(\projmix_{k})$ but $\projcmp_{j}\in\supp(\projmix_{k'})$ with $k\ne k'$. Then \ref{lem:metric:ineq:M7} of Lemma~\ref{lem:metric:ineq:pair} implies that 
\begin{align*}
\probm(\projcmp_{i}, \projcmp_{j}) 
\ge 2\gensep(\delta)
> \gensep(\delta)
= \separ(\projmix(\assg)),
\end{align*}

\noindent
which proves \eqref{eq:prop:pop:cluster:1}. The equivalence \eqref{eq:prop:pop:cluster:2} follows from a similar argument.

In particular, \eqref{eq:prop:pop:cluster:1} and \eqref{eq:prop:pop:cluster:2} together imply that if $\alpha(i)=k=\alpha(j)$, single-linkage hierarchical clustering will join components $q_{i}$ and $q_{j}$ before including any component $q_{\ell}$ such that $\alpha(\ell)\ne k$. Thus, cutting the resulting dendrogram at any level $t\in(\separ(\projmix(\assg)), 2\separ(\projmix(\assg)))$ will produce $K$ clusters corresponding to $\assg^{-1}(1),\ldots,\assg^{-1}(K)$.
\qed

\subsection{Proof of Theorem~\ref{thm:main:ident}}
\label{app:proofs:thm:main:ident}

We first note the following consequence of regularity (Definition~\ref{defn:regular}):
\begin{lemma}
\label{lem:proj:consistent}
If $\truemix$ is a regular mixing measure then $\mixm(\projmix(\assg_{L}),\truemix)\to 0$ as $L\to\infty$ for any $\assg_{L}\in\assignseq(\truemix)$.
\end{lemma}

\begin{proof}
This follows from Definition~\ref{defn:regularseq} and Lemma~\ref{lem:wasserstein:atom}.
\end{proof}

\begin{proof}[Proof of Theorem~\ref{thm:main:ident}]
Since $\mixset$ is clusterable by assumption, there is a function $\clustmap_{L}:\mixingmap_{L}(\mixset)\to\assignmaps$ such that 
\begin{align}
\label{eq:thm:main:ident:1}
\lim_{L\to\infty}\probm(\proj_{k}(\assg_{L}), \truecmp_{k}) = 0
\quad\text{and}\quad
\lim_{L\to\infty}|\projclustwgt_{k}(\assg_{L}) - \truewgt_{k}| = 0
\quad \forall\,k=1,\ldots,K,
\end{align}

\noindent
where $\assg_{L}=\clustmap_{L}(\projmix)$.
This defines a function $F_{L}:\mix(\mixset)\to\assignproj$ by 
\begin{align}
\label{eq:defn:mixmap}
F_{L}(\true)
= \projmix(\assg_{L}),
\quad\text{where}\quad
\projmix = \mixingmap_{L}(\projmap_{L}(\true))
\text{ and }
\assg_{L} = \clustmap_{L}(\projmix).
\end{align}

\noindent
The function $F_{L}$ defines a unique, well-defined procedure for associating to a mixture distribution $\true=\mixfcn(\truemix)$ a mixing measure. Finally, define
\begin{align*}
h(\true)
:= \lim_{L\to\infty}F_{L}(\true),
\end{align*}

\noindent
where convergence of $F_{L}(\true)$ is understood to be with respect to the Wasserstein metric $\mixm$. It remains to show that $h(\mixfcn(\truemix))=\truemix$ for all $\truemix\in\mixset$, i.e. $\mixm(F_{L}(\mixfcn(\truemix)), \truemix)\to0$.
But this follows from Lemma~\ref{lem:proj:consistent}.

Finally, to show that $\mixfcn$ is a bijection, we need to show that if $\mixfcn(\truemix)=\mixfcn(\truemix')$ for some $\truemix,\truemix'\in\mixset$, then $\truemix=\truemix'$. We just proved that $\lim_{L\to\infty}F_{L}(\mixfcn(\truemix))=\truemix$ for all $\truemix\in\mixset$, and hence
\begin{align*}
\truemix'
= \lim_{L\to\infty}F_{L}(\mixfcn(\truemix'))
= \lim_{L\to\infty}F_{L}(\mixfcn(\truemix))
= \truemix.
\end{align*}

\noindent
This proves that $m$ is injective, and surjectivity is obvious since $m$ is onto $\mix(\mixset)$ by definition. Thus $m$ is a bijection as claimed and the proof is complete.
\end{proof}

\subsection{Proof of Proposition~\ref{prop:est:cluster}}
\label{app:proofs:prop:est:cluster}

Define $\Hdiam_{0}:=\sup_{k}\Hdiam(\projmix_{k})$, so that $\gensep(\delta)=\Hdiam_{0}+\delta$.
\begin{prop}
\label{prop:est:cluster:base}
Suppose that $\eps>0$ and $\delta>0$ satisfy
\begin{align*}
3\eps - 2\delta
< \Hdiam_{0}.
\end{align*}

\noindent
Then under \ref{lem:metric:ineq:A1}-\ref{lem:metric:ineq:A3} we have $\probm(\estcmp_{i},\estcmp_{j})\le\gensep(2\eps)$ if and only if $\probm(\projcmp_{i}, \projcmp_{j}) \le \gensep(\delta)$.
\end{prop}

\begin{proof}
Suppose $\probm(\projcmp_{i}, \projcmp_{j}) \le \gensep(\delta)$, which implies $\alpha(i)=\alpha(j)$ by Proposition~\ref{prop:pop:cluster}. Then \ref{lem:metric:ineq:M6} of Lemma~\ref{lem:metric:ineq:pair} implies $\probm(\estcmp_{i}, \estcmp_{j})\le \gensep(2\eps)$.
Conversely, suppose $\probm(\estcmp_{i},\estcmp_{j})\le\gensep(2\eps)$ but $\probm(\projcmp_{i}, \projcmp_{j}) > \gensep(\delta)$. By Proposition~\ref{prop:pop:cluster}, this means that $\assg(i)\ne\assg(j)$, and invoking \ref{lem:metric:ineq:M8} of Lemma~\ref{lem:metric:ineq:pair} we deduce that $\probm(\estcmp_{i}, \estcmp_{j}) \ge 2\gensep(\delta) - \eps$. Thus, since also $\probm(\estcmp_{i}, \estcmp_{j})\le \gensep(2\eps)$, we have
\begin{align*}
&\phantom{\iff\,\,\,}
2\gensep(\delta) - \eps
\le\probm(\estcmp_{i}, \estcmp_{j})
\le \gensep(2\eps) \\
&\iff
2\Hdiam_{0}+2\delta-\eps \le \Hdiam_{0}+2\eps \\
&\iff
\Hdiam_{0}\le 3\eps-2\delta,
\end{align*}

\noindent
which contradicts the assumption that $\Hdiam_{0}> 3\eps-2\delta$.
\end{proof}

\begin{proof}[Proof of Proposition~\ref{prop:est:cluster}]
Note that $\wh{\separ}=\gensep(2\eps)$ and $\separ(\projmix(\assg))=\gensep(\delta)$, so that Propositions~\ref{prop:est:cluster:base} and~\ref{prop:pop:cluster} together imply $\probm(\estcmp_{i},\estcmp_{j})\le\gensep(2\eps)\iff\alpha(i)=\alpha(j)$. In fact, we also have $\probm(\estcmp_{i},\estcmp_{j})\ge2\gensep(\delta)-\eps\iff\alpha(i)\ne\alpha(j)$. Thus, as long as $2\gensep(\delta)-\eps > \gensep(2\eps)$, single-linkage clustering will recover the $K$ components (i.e. by cutting the dendrogram at any $t\in(\gensep(2\eps),2\gensep(\delta)-\eps)$). But $2\gensep(\delta)-\eps > \gensep(2\eps)\iff\Hdiam_{0}> 3\eps-2\delta$, which is true by assumption.
\end{proof}

\subsection{Proof of Theorem~\ref{thm:main:learn}}
\label{app:proofs:thm:main:learn}

We first need the following lemma:

\begin{lemma}
\label{lem:wasserstein:cond}
Assume $\genmix^{n},\genmix\in\mixing_{L}(\base)$. If $\mixm(\genmix^{n},\genmix)\to0$ then for any $\assg\in\assignmaps$,
\begin{align}
\label{eq:lem:wasserstein:cond:mix}
\mixm(\genmix^{n}(\assg), \genmix(\assg)) &\to 0, \\
\label{eq:lem:wasserstein:cond:proj}
\probm(\genproj^{n}_{k}(\assg), \genproj_{k}(\assg)) &\to 0, \\
\label{eq:lem:wasserstein:cond:wgt}
|\genclustwgt^{n}_{k}(\assg) - \genclustwgt_{k}(\assg)| &\to 0.
\end{align}

\noindent
Furthermore, each sequence \eqref{eq:lem:wasserstein:cond:mix}, \eqref{eq:lem:wasserstein:cond:proj}, and \eqref{eq:lem:wasserstein:cond:wgt} converges at the same rate as $\mixm(\genmix^{n},\genmix)$.
\end{lemma}

\begin{proof}
By Lemma~\ref{lem:wasserstein:atom}, we may assume without loss of generality that $\probm(\gencmp_{\ell}^{n},\gencmp_{\ell})\to0$ and $|\genwgt^{n}_{\ell} - \genwgt_{\ell}|\to 0$ at the same rate as $\mixm(\genmix^{n},\genmix)$. Then
\begin{align*}
\genclustwgt_{k}^{n}(\assg)
= \sum_{\ell\in\assg^{-1}(k)}\genwgt_{\ell}^{n}
\to \sum_{\ell\in\assg^{-1}(k)}\genwgt_{\ell}
= \genclustwgt_{k}(\assg)
\end{align*}

\noindent
and furthermore $\genwgt_{\ell}^{n}/\genclustwgt_{k}^{n}(\assg)\to\genwgt_{\ell}/\genclustwgt_{k}(\assg)$ for each $\ell$ at the same rate as $\mixm(\genmix^{n},\genmix)$. This proves \eqref{eq:lem:wasserstein:cond:wgt}, and \eqref{eq:lem:wasserstein:cond:proj} follows similarly. In particular, each of the atoms and weights in $\genmix^{n}(\assg)$ converges to an atom and weight in $\genmix(\assg)$. Invoking Lemma~\ref{lem:wasserstein:atom} once again, we deduce \eqref{eq:lem:wasserstein:cond:mix}. The proof is complete.
\end{proof}

\begin{proof}[Proof of Theorem~\ref{thm:main:learn}]
We first prove $\mixm(\estmix(\wh{\assg}_{L,n}), \projmix(\assg_{L}))\to0$. By assumption, $\estmix$ is a $\mixm$-consistent estimate of $\projmix$, i.e. $\mixm(\estmix,\projmix)\to0$ as $n\to\infty$. Proposition~\ref{prop:est:cluster} implies $\wh{\assg}_{L,n}=\assg_{L}$, and hence 
\begin{align}
\label{eq:wass:conv:est}
\mixm(\estmix(\wh{\assg}_{L,n}), \projmix(\assg_{L}))
= \mixm(\estmix(\assg_{L}), \projmix(\assg_{L}))
\to 0
\end{align}

\noindent 
by Lemma~\ref{lem:wasserstein:cond}. 
Thus,
\begin{align*}
\mixm(\estmix(\wh{\assg}_{L,n}), \truemix)
\le \mixm(\estmix(\wh{\assg}_{L,n}), \projmix(\assg_{L})) 
+ \mixm(\projmix(\assg_{L}), \truemix) 
\to0
\end{align*}

\noindent
by \eqref{eq:wass:conv:est} 
and Lemma~\ref{lem:proj:consistent}. This proves \eqref{eq:thm:main:learn:1a}, as desired.
\end{proof}

\begin{remark}
\label{rem:Lvsn}
It is interesting to study how $L$ must grow as a function of $n$ in order to achieve consistency in Theorem~\ref{thm:main:learn}. In fact, since
\begin{align*}
\mixm(\estmix(\wh{\assg}_{L,n}), \truemix)
&\le \mixm(\estmix(\wh{\assg}_{L,n}), \projmix(\assg_{L})) 
+ \mixm(\projmix(\assg_{L}), \truemix) \\
&:= \psi_{L,n} + \Psi_{L}
\end{align*}

\noindent
it is clear that $L=L_{n}$ is governed by the rate of convergence $\psi_{L,n}$ of $\mixm(\estmix(\wh{\assg}_{L,n}), \projmix(\assg_{L}))$. Assuming $\wh{\assg}_{L,n}=\assg_{L}$ (this follows, e.g. from Proposition~\ref{prop:est:cluster}) the convergence rate of $\mixm(\estmix(\wh{\assg}_{L,n}), \projmix(\assg_{L}))$ is the same as $\mixm(\estmix, \projmix)$ by Lemma~\ref{lem:wasserstein:cond}. Thus, this problem reduces to studying the rate of convergence of the estimator $\estmix$. Conveniently, by Theorem~2 in \citet{nguyen2013}, this can be bounded above as $L\to\infty$ by the variational distance $\dTV(\mixfcn(\estmix(\wh{\assg}_{L,n})), \mixfcn(\projmix(\assg_{L})))$.
\end{remark}

\subsection{Proof of Theorem~\ref{thm:main:part}}
\label{app:proofs:thm:main:part}

Write $\projdenscmp_{\ell}$ for the density of $\projcmp_{\ell}$ and $\projdens$ for the density of $\proj$, and similarly $\estdenscmp_{\ell}$ and $\estdens$ for the densities of $\estcmp_{\ell}$ and $\estproj$, respectively. To reduce notational overload in the proof, we will suppress the dependence on $\estassg$ (cf. \eqref{eq:est:dens}) below. To this end, define 

\begin{align*}
\estdens_{k}(x)
&:= \frac1{\estclustwgt_{k}(\wh{\assg}_{L,n})}\sum_{\ell\in\wh{\assg}_{L,n}^{-1}(k)}\estwgt_{\ell}\estdenscmp_{\ell}(x),&
\estclustwgt_{k}(\wh{\assg}_{L,n})
&:= \sum_{\ell\in\wh{\assg}_{L,n}^{-1}(k)}\estwgt_{\ell}, \\
\projdens_{k}(x)
&:= \frac1{\projclustwgt_{k}(\assg_{L})}\sum_{\ell\in\assg_{L}^{-1}(k)}\projwgt_{\ell}\projdenscmp_{\ell}(x),&
\projclustwgt_{k}(\assg_{L})
&:= \sum_{\ell\in\assg_{L}^{-1}(k)}\projwgt_{\ell}.
\end{align*}

\noindent
Then $\estdens_{k}(x)$ and $\projdens_{k}(x)$ are the densities of $\estproj_{k}(\wh{\assg}_{L,n})$ and $\proj$, respectively.

\begin{lemma}
\label{lem:clust:uniform}
If $\estdens_{k}\to\truedenscmp_{k}$ uniformly on $\base$ for each $k=1,\ldots,K$, then there exists a sequence $\tseq\to0$ such that $\wh{c}_{L,n}(x) = c(x)$ for all $x\in X-\eqset(\tseq)$.
\end{lemma}

\begin{proof}
Define
\begin{align}
\label{eq:lem:clust:uniform:1}
\tseq:=2\sup_{k}\sup_{x\in\base}|\estdens_{k}(x)-\truedenscmp_{k}(x)|\ge 0
\end{align}

\noindent
and note that $\tseq\to0$ since $\estdens_{k}\to \truedenscmp_{k}$ uniformly for each $k$. Then for any $x\notin\eqset(\tseq)$ and $i\ne j$, $|f_{i}(x) - f_{j}(x)|> \tseq$, which means either $f_{i}(x)>f_{j}(x)+\tseq$ or $f_{j}(x)>f_{i}(x)+\tseq$. Taking $i=c(x)$, it follows that
\begin{align}
\label{eq:lem:clust:uniform:2}
\truedenscmp_{c(x)}(x)
> \sup_{j\ne c(x)}\truedenscmp_{j}(x) + \tseq 
\quad\forall x\notin\eqset(\tseq).
\end{align}

\noindent
Thus for any $j\ne c(x)$,
\begin{align*}
\estdens_{c(x)}(x)
\overset{(a)}{>} \truedenscmp_{c(x)}(x) - \frac{\tseq}{2}
\overset{(b)}{>} \truedenscmp_{j}(x) + \frac{\tseq}{2}
\overset{(c)}{>} \estdens_{j}(x)
\end{align*}

\noindent
(a) follows from \eqref{eq:lem:clust:uniform:2}; (b) follows from \eqref{eq:lem:clust:uniform:1}, and (c) follows again from \eqref{eq:lem:clust:uniform:2}. It follows that $\wh{c}_{L,n}(x)=c(x)$ for all $x\notin\eqset(\tseq)$, as desired.
\end{proof}

\begin{lemma}
\label{lem:part:conv}
Under the assumptions of Lemma~\ref{lem:clust:uniform}, it follows that $\wh{A}_{L,n,k}\symdiff A_{k}\subset\eqset(\tseq)$ for all $k=1,\ldots,K$. In particular,
\begin{align*}
\bigcup_{k=1}^{K}\wh{A}_{L,n,k}\symdiff A_{k}
\subset\eqset(\tseq).
\end{align*}
\end{lemma}

\begin{proof}
Lemma~\ref{lem:clust:uniform} implies that $\wh{c}_{L,n}(x)=c(x):=k$ for all $x\notin\eqset(\tseq)$. 
In particular, $\wh{A}_{L,n,k}\cap \eqset(\tseq)^{c}=A_{k}\cap \eqset(\tseq)^{c}$. Lemma~\ref{lem:symdiff:contain} thus implies that $\wh{A}_{L,n,k}\symdiff A_{k}\subset\eqset(\tseq)$, as desired.
\end{proof}

\begin{proof}[Proof of Theorem~\ref{thm:main:part}]
This is an immediate consequence of Lemmas~\ref{lem:clust:uniform} and~\ref{lem:part:conv}.
\end{proof}

\section{Examples}
\label{app:ex}

In this Appendix, we prove the claims made in Examples~\ref{ex:disjmix},~\ref{ex:metamix} and~\ref{ex:identmix} from Section~\ref{sec:ident:regularity}. Define the following families of mixing measures:
\begin{align}
\label{eq:defn:disjmix}
\disjmix(\lebesgue,K)
&:= 
\Bigg\{
\sum_{k=1}^{K}\truewgt_{k}\delta_{\truecmp_{k}}
:
\truecmp_{k}\ll \lebesgue,\,
\truewgt_{k}\ge 0,\, 
\sum_{k=1}^{K}\truewgt_{k}=1
\Bigg\}, \\
\label{eq:defn:metamix}
\metamix(\projset,K)
&:= 
\Bigg\{
\sum_{k=1}^{K}\truewgt_{k}\delta_{\mixfcn(\gaussmix_{k})}
:
\gaussmix_{k}\in\projset,
\,
\truewgt_{k}\ge 0,\, 
\sum_{k=1}^{K}\truewgt_{k}=1
\Bigg\}.
\end{align}

\noindent
Here, $\lebesgue\in\probs(\base)$ is a fixed dominating measure and $\projset\subset\mixing(\base)$ is a Borel set.

\subsection{Disjoint mixtures}
The family $\disjmix(\lebesgue,K)$ defined in \eqref{eq:defn:disjmix} corresponds to Example~\ref{ex:disjmix}. For any $\truemix\in\disjmix(\lebesgue,K)$, let $\truedens=\sum_{k=1}^{K}\truewgt_{k}\truedenscmp_{k}$ denote the density of $\true=\mixfcn(\truemix)$. Let $\tick{\mixgaussians}\subset\mixgaussians$ be a compact set of Gaussian mixtures and $\projdens_{L}=\sum_{\ell=1}^{L}\projwgt_{\ell}\projdenscmp_{\ell}$ be a corresponding $\dTV$-projection onto $\mix(\tick{\mixgaussians_{L}})$, with $\projdenscmp_{\ell}\sim\normalN(\projmean_{\ell},\projvar_{\ell})$.

The following lemma implies $\tick{\mixgaussians}$-regularity of $\truemix$ as long as the $\truedenscmp_{k}$ have disjoint supports:

\begin{lemma}
\label{lem:disjmix:reg}
Suppose $\truemix\in\disjmix(\lebesgue,K)$ and $\tick{\mixgaussians}\subset\mixgaussians$ is a compact set such that $\mix(\tick{\mixgaussians})$ is dense in $\mix(\disjmix(\lebesgue,K))$. Define $E:=\cup_{k=1}^{K}E_{k}$ and an assignment sequence $\assg_{L}$ by 
\begin{align*}
\assg_{L}(\ell) = \begin{cases}
k, & \projmean_{\ell}\in E_{k}, \\
1, & \projmean_{\ell}\in E^{c}.
\end{cases}
\end{align*}

\noindent
If the $E_{k}$ are disjoint, then $\projclustwgt_{k}(\assg_{L})\to\lambda_{k}$ and $\norm{\projdens_{k}(\assg_{L}) - f_{k}}_{1}\to0$ for each $k$.
\end{lemma}

\begin{remark}
As the proof indicates, the conclusion of Lemma~\ref{lem:disjmix:reg} remains true if $\tick{\mixgaussians}$ is not compact (e.g. if $\tick{\mixgaussians}=\mixgaussians$), however, in this case the $\dTV$-projection may not be well-defined. In this case, the conclusion holds for any $\projdens_{L}\in\mix(\mixgaussians_{L})$ such that $\dTV(\projdens_{L},\truedens)\to0$ as $L\to\infty$. Furthermore, the assumption that $\cl{\mix(\tick{\mixgaussians})}=\mix(\disjmix(\lebesgue,K))$ can be trivially relaxed to the requirement that $\truedens\in\cl{\mix(\tick{\mixgaussians})} \cap \mix(\disjmix(\lebesgue,K))$. 
\end{remark}

\begin{proof}[Proof of Lemma~\ref{lem:disjmix:reg}]
To reduce notational clutter, we omit the index $L$ and write $\assg_{K}=\assg$ and $\projdens_{L}=\projdens$.
We have
\begin{align*}
\norm{\projdens - \truedens}_{1}
&= \sum_{k=1}^{K} \int_{E_{k}}|\projdens - \truedens|\dx + \int_{E^{c}}\projdens\dx \\
&= \sum_{k=1}^{K} \int_{E_{k}}|\underbrace{\sum_{j\ne k}\projclustwgt_{j}(\assg)\projdens_{j}(\assg)}_{:=h_{k,L}} + \underbrace{\varpi_{k}^{*}(\assg)\projdens_{k}(\assg) - \lambda_{k}f_{k}}_{:=-g_{k,L}}|\dx + \int_{E^{c}}\projdens_{L}\dx \\
&= \sum_{k=1}^{K} \int_{E_{k}}|h_{k,L}-g_{k,L}|\dx + \int_{E^{c}}\projdens_{L}\dx.
\end{align*}

\noindent
Since $\norm{\projdens - \truedens}_{1}\to 0$, it follows that $\int_{E^{c}}\projdens\dx\to 0$ and hence $\int_{E_{k}}\projclustwgt_{j}(\assg)\projdens_{j}(\assg)\dx\to 0$ for all $j\ne k$. 
This is because off of $E_j$, $\projclustwgt_{j}(\assg)\projdens_{j}(\assg)$ is monotonic (i.e. decreasing to the right and increasing from the left), which implies (assuming without loss of generality that $\intup_{j}<\intlo_{k}$---
i.e. $E_{j}$ lies to the left of $E_{k}$ on the real line),
\begin{align*}
\int_{\intlo_{k}}^{\intup_{k}}\projclustwgt_{j}(\assg)\projdens_{j}(\assg)
< \int_{\intup_{j}}^{\intlo_{k}}\projclustwgt_{j}(\assg)\projdens_{j}(\assg)
\le \int_{E^{c}}\projdens\to 0.
\end{align*}
Furthermore, it follows that $\int_{E_{k}}|h_{k,L}-g_{k,L}|\dx\to0$ for each $k$, which implies $\int_{E_{k}}|g_{k,L}|\dx\to 0$ since 
\begin{align*}
\int_{E_{k}}|g_{k,L}|\dx - \int_{E_{k}}|h_{k,L}|\dx
\le \int_{E_{k}}|h_{k,L}-g_{k,L}|\dx
\to 0.
\end{align*}

\noindent
Thus
\begin{align*}
\int_{E_{k}}|\varpi_{k}^{*}(\assg)\projdens_{k}(\assg) - \lambda_{k}f_{k}|\dx
\to 0
\quad\forall k,
\end{align*}

\noindent
i.e. $\varpi_{k}^{*}(\assg)\projdens_{k}(\assg) \to \lambda_{k}f_{k}$ in $L^{1}$ (on $E_{k}$), which implies
\begin{align*}
\Big|\int_{E_{k}}\big\{\varpi_{k}^{*}(\assg)\projdens_{k}(\assg) - \lambda_{k}f_{k}\big\}\dx\Big|
= \Big|\varpi_{k}^{*}(\assg)\int_{E_{k}}\projdens_{k}(\assg)\dx - \lambda_{k}\Big|
\to 0.
\end{align*}

\noindent
However, since $\int_{E^{c}}Q_{L}^{*}\dx\to 0$, we have $\int_{E_{k}^{c}}\projdens_{k}(\assg)\dx\to 0$ and hence also $\int_{E_{k}}\projdens_{k}(\assg)\dx\to 1$. It follows that $\varpi_{k}^{*}(\assg)\to\lambda_{k}$ and similarly $\norm{\projdens_{k}(\assg) - f_{k}}_{1}\to0$, as desired.
\end{proof}

\subsection{Mixtures of finite mixtures}

By restricting to finite mixtures $\projset_{0}$, the family $\metamix(\projset_{0},K)$ defined in \eqref{eq:defn:metamix} corresponds to Example~\ref{ex:metamix}. The following lemma establishes $\projset$-regularity of this family.

\begin{lemma}
\label{lem:metamix:reg}
Suppose $\projset\subset\mixing(\base)$ satisfies Condition~\ref{condn:A} and $\truemix\in\metamix(\projset_{0},K)$. If the $\gaussmix_{k}$ have disjoint supports, then $\truemix$ is $\projset$-regular under any metric $\probm$.
\end{lemma}

\begin{proof}
Let $\truemix\in\metamix(\projset_{0};\,K)$ and define $L_{0}:=\sum_{k=1}^{K}|\supp(\gaussmix_{k})|$.
For any $L\ge L_{0}$, the $\probm$-projection of $\true$ onto $\mix(\projset_{L})$ is $\proj=\true$ since $\true\in\mix(\projset_{L})$ and $\probm$ is a metric (i.e. $\probm(\proj,\true)=0\iff\proj=\true$). This verifies part \ref{defn:regular:a} of Definition~\ref{defn:regular}. To verify part \ref{defn:regular:b}, first write 
\begin{align*}
\true
= \proj
= \sum_{\ell=1}^{L_{0}}\projwgt_{\ell}\projcmp_{\ell},
\quad
\truecmp_{k}
= \sum_{\ell\in B_{k}}\projwgt_{\ell}\projcmp_{\ell}.
\end{align*}

\noindent
The set $B_{k}\subset[L_{0}]$ indexes all of the components of $\true$ that contribute to $\truecmp_{k}$, i.e. $B_{k} = \{\ell\in[L_{0}] : \projcmp_{\ell}\in\supp(\truemix_{k})\}$. Now define
\begin{align*}
\assg(\ell) = k
\iff
\ell\in B_{k}.
\end{align*}

\noindent
Then clearly $\proj_{k}(\assg)=\truecmp_{k}$ and $\projclustwgt_{k}(\assg)=\truewgt_{k}$, which establishes \ref{defn:regular:b}. 
\qedhere
\end{proof}

\subsection{Infinite mixtures}

Suppose $\projset\subset\mixing(\base)$ is compact and identifiable; i.e. restricted to $\projset$, the canonical map $\mixfcn:\projset\to\probs(\base)$ is injective. 
Then the family $\metamix(\projset,K)$ defined in \eqref{eq:defn:metamix} corresponds to Example~\ref{ex:identmix}. As the following examples illustrate, this family encompasses a wide range of nonparametric mixtures:
\begin{itemize}
\item \emph{Parametric mixture components with unbounded support.} The family $\metamix(\projset_{0},K)$ is obviously a special case of this family, including mixtures whose components are finite Gaussian mixtures and exponential family mixtures.
\item \emph{Infinite Gaussian mixtures.} Although arbitrary mixtures of Gaussians need not be identifiable \citep{teicher1960}, convolutional mixtures \citep{teicher1960,nguyen2013} and mixtures of Gaussians with fixed means are identifiable. More generally, we may allow both the means and variances to vary under certain conditions, see \citet{bruni1985identifiability}. If the $\truecmp_{k}$ come from any of these general families of Gaussian mixtures, then the resulting nonparametric mixture model is a special case of $\metamix(\projset,K)$.
\end{itemize}

\noindent
A mixing measure $\truemix\in\metamix(\projset,K)$ is of the form $\sum_{k=1}^{K}\truewgt_{k}\delta_{\truecmp_{k}}$, with $\truecmp_{k}=\mixfcn(\gaussmix_{k})$. Thus $\truemix$ uniquely defines a $\projset$-mixture by $\gaussmix:=\sum_{k}\truewgt_{k}\gaussmix_{k}$. 

Through this subsection we assume that $\probm$ is a $\phi$-divergence \citep{nguyen2013,csiszar1967information}, which we recall means that 
\begin{align}
\probm(\mu,\nu)
=\int_{\base}\phi\Big(\dod{\mu}{\nu} \Big)\,\dif\nu
\end{align}

\noindent
for some convex function $\phi$. If $\mu$ is not absolutely continuous with respect to $\nu$, we adopt the usual convention that $\probm(\mu,\nu)=+\infty$. Examples of $\phi$-divergences include the Hellinger and variational metrics.

Recall that $\mixfcn$ is a map between the metric spaces $(\mixing(\base),\mixm)$ and $(\probs(\base),\probm)$. We first record an important fact regarding the continuity of this map:
\begin{lemma}
\label{lem:mixfcn:cts}
Suppose $\projset\subset\mixing(\base)$ is compact and identifiable and $\probm$ is a $\phi$-divergence. Then $\mixfcn:\projset\to\probs(\base)$ is 1-Lipschitz continuous, and its inverse is continuous.
\end{lemma}

\begin{proof}
This follows from Lemma~1 in \citet{nguyen2013}, noting that the argument applies just as well to arbitrary measures $\mu,\nu\in\probs(\base)$. Continuity of the inverse follows from the compactness of $\projset$.
\end{proof}

The following lemma establishes $\projset$-regularity of $\metamix(\projset,K)$. Its proof uses the well-known fact that the Wasserstein metric metrizes the topology of weak convergence, which we denote by $\wconv$.

\begin{lemma}
\label{lem:identmix:reg}
Suppose $\truemix\in\metamix(\projset,K)$ and define $B_{k}:=\supp(\gaussmix_{k})$. Assume that (a) $B_{1},\ldots,B_{K}$ are disjoint, compact sets and (b) Each $B_{k}$ is a $\gaussmix$-continuity set. Then $\truemix$ is $\projset$-regular under any $\phi$-divergence.
\end{lemma}

\begin{proof} 
Let $\proj_{L}\in\mix(\projset_{L})$ be any $\probm$-projection of $\true$ and $\projmix_{L}$ its corresponding mixing measure.
We wish to show that there exists an assignment sequence $\assgseq$ such that $\probm(\proj_{k}(\assg_{L}),\truecmp_{k})\to0$ and $|\projclustwgt_{k}(\assg_{L}) - \truewgt_{k}|\to0$. 
By Lemma~\ref{lem:mixfcn:cts}, it suffices to show that $\projmix_{L,k}(\assg_{L})\wconv\gaussmix_{k}$ and $\projclustwgt_{L,k}(\assg_{L})\to\truewgt_{k}$. Consider the assignment sequence defined by $\assg_{L}(g)=k$ if and only if $g\in B_{k}$. Note that $\assg_{L}$ does not depend on $\projmix_{L}$.

We first note that $\projmix_{L}\wconv\gaussmix$: This follows again from Lemma~\ref{lem:mixfcn:cts} since $\probm(\mixfcn(\projmix_{L}),\mixfcn(\gaussmix))\to0$. Thus, in particular, since $B_{k}$ is a $\gaussmix$-continuity set we have both $\projmix_{L}(\cdot\,\given B_{k})\wconv \gaussmix(\cdot\,\given B_{k})$ and $\projmix_{L}(B_{k})\to\gaussmix(B_{k})$. But $\projmix_{L,k}(\assg_{L})=\projmix_{L}(\cdot\,\given B_{k})$ and by Lemma~\ref{lem:decomp:cond}, we have $\gaussmix_{k}=\gaussmix(\cdot\,\given B_{k})$. It follows that $\projmix_{L,k}(\assg_{L})\wconv\gaussmix_{k}$. Similarly, by Lemma~\ref{lem:decomp:cond}, we have $\projclustwgt_{L,k}(\assg_{L})=\projmix_{L}(B_{k})\to\gaussmix(B_{k})=\truewgt_{k}$. The proof is complete.
\qedhere

\end{proof}

\section{Technical lemmas}
\label{app:tech}

\begin{lemma}
\label{lem:symdiff:contain}
$A\cap Z^{c}=B\cap Z^{c}\implies A\symdiff B\subset Z$.
\end{lemma}

\begin{proof}
Suppose $x\in A\symdiff B$. There are two cases: (i) $x\in A-B$, (ii) $x\in B-A$. We prove (i); the proof for (ii) is similar. Suppose $x\in A-B$ but also $x\in Z^{c}$. Then $x\in A\cap Z^{c}= B\cap Z^{c}$, whence $x\in B$, which contradicts (i).
\end{proof}

\begin{lemma}
\label{lem:decomp:cond}
Suppose $\gaussmix:=\sum_{k}\truewgt_{k}\gaussmix_{k}\in\mixing(\base)$, $B_{k}:=\supp(\gaussmix_{k})$, and assume that $B_{1},\ldots,B_{K}$ partition $B:=\supp(\gaussmix)$. Then $\gaussmix_{k} = \gaussmix(\,\cdot\given B_{k})$ and $\truewgt_{k} = \gaussmix(B_{k})$.
\end{lemma}

\begin{proof} 
For any Borel set $A$ 
we have
\begin{align*}
\sum_{k=1}^{K}\gaussmix(A\cap B_{k})
= \gaussmix(A) 
= \sum_{k=1}^{K}\truewgt_{k}\gaussmix_{k}(A).
\end{align*}

\noindent
Choosing $A=B_{j}$, we deduce that $\truewgt_{j}=\gaussmix(B_{j})$. Then for any $A\subset B_{j}$, $\gaussmix(A) = \truewgt_{j}\gaussmix_{j}(A)$ since the $B_{k}$ form a partition. Using $\truewgt_{j}=\gaussmix(B_{j})$ and $A\subset B_{j}$, this implies that
\begin{align*}
\gaussmix_{j}(A)
= \frac{\gaussmix(A)}{\gaussmix(B_{j})}
= \frac{\gaussmix(A\cap B_{j})}{\gaussmix(B_{j})}
= \gaussmix(A\given B_{j}).
\end{align*}

\noindent
Finally, for arbitrary Borel $A$, we have
\begin{align*}
\gaussmix_{j}(A)
= \gaussmix_{j}(A\cap B_{j}) + \underbrace{\gaussmix_{j}(A\cap B_{j}^{c})}_{=0}
\overset{(i)}{=} \gaussmix_{j}(A\cap B_{j})
\overset{(ii)}{=} \gaussmix(A\given B_{j})
\end{align*}

\noindent
since (i) $\supp(\gaussmix_{j})=B_{j}$ and (ii) $A\cap B_{j}\subset B_{j}$.
\end{proof}

What follows are standard results on the Wasserstein metric. For completeness, we provide complete proofs here. 
Let $(\Theta,\basem)$ be a metric space and $G=\sum_{i=1}^{m}p_{i}\delta_{\theta_{i}}$, $G'=\sum_{j=1}^{m}p_{j}\delta_{\theta_{j}'}$ be two discrete probability measures on $\Theta$ with $m$ atoms each. Then the definition of the $L_{r}$-Wasserstein distance \eqref{eq:defn:wasserstein} is equivalent to the following:
\begin{align*}
\mixm(G,G')
= \inf\Bigg\{
\Big[\sum_{i,j}\coupling_{ij}\basem^{r}&(\theta_{i},\theta_{j}')\Big]^{1/r} :
0\le \coupling_{ij}\le 1,\, \\
&\sum_{i,j}\coupling_{ij}=1,\,
\sum_{i}\coupling_{ij}=p_{j}',\,
\sum_{j}\coupling_{ij}=p_{i}
\Bigg\}.
\end{align*}

\noindent
In the sequel, we write $G_{n}=\sum_{i=1}^{m}p_{n,i}\delta_{\theta_{n,i}}$ for a sequence of discrete distributions with exactly $m$ atoms each.

\begin{remark}
\label{rem:wass:dist}
Lemma~\ref{lem:wasserstein:atom} below encodes a crucial property of $\mixm$ that is needed in the proofs. Thus, in contrast to $\probm$ which can be any metric on $\probs(\base)$, our results are particular to the Wasserstein distance. Notwithstanding, it is likely that our results can be extended to any metric on $\mixing(\base)$ that satisfies Lemma~\ref{lem:wasserstein:atom}.
\end{remark}

\begin{lemma}
\label{lem:wasserstein:atom}
$W_{r}(G_{n},G)\to0$ if and only if there are permutations $\tau_{n}:[m]\to[m]$ such that 
\begin{enumerate}
\item[(a)] $\basem(\theta_{n,\tau_{n}(i)},\theta_{i})\to0$ for each $i$.
\item[(b)] $|p_{n,\tau_{n}(i)}-p_{i}|\to0$ for each $i$.
\end{enumerate}
\end{lemma}

\begin{proof}
The $\implies$ direction is Lemma~\ref{lem:wasserstein:convzero} below. The $\Longleftarrow$ direction follows from Lemma~\ref{lem:wasserstein:atomic} by noting that $G_{n}$ and $G$ are invariant to permutations of the index $i$.
\end{proof}

\begin{lemma}
\label{lem:wasserstein:convzero}
If $W_{r}(G_{n},G)\to0$ then there are permutations $\tau_{n}:[m]\to[m]$ such that 
\begin{enumerate}
\item[(a)] $\basem(\theta_{n,\tau_{n}(i)},\theta_{i})\to0$ for each $i$.
\item[(b)] $|p_{n,\tau_{n}(i)}-p_{i}|\to0$ for each $i$.
\end{enumerate}
\end{lemma}

\begin{proof}
We use the fact that Wasserstein convergence implies weak convergence, which is in turn equivalent to convergence of open sets, i.e. $G_{n}(U)\to G(U)$ for all open sets $U$. Choose $R>0$ so small that $B(\theta_{i},R)\cap B(\theta_{j},R)=\emptyset$ for all $i\ne j$, whence $G(B(\theta_{i},r))=p_{i}$ for all $0<r<R$. By weak convergence, we thus have for each $i$ and any $0<r<R$
\begin{align*}
G_{n}(B(\theta_{i},r))
\to G(B(\theta_{i},r))
= p_{i}.
\end{align*}

\noindent
Thus, for sufficiently large $n$, $G_{n}$ assigns positive probability to the ball $B(\theta_{i},r)$, which means that for each $i$ there is some $j\in[k]$ such that $\theta_{n,i}\in B(\theta_{j},r)$. Setting $\tau_{n}(j)=i$  and taking $r\to0$ completes the proof.
\end{proof}

\begin{lemma}
\label{lem:wasserstein:atomic}
If $\lim_{n\to\infty}\theta_{n,i}=\theta_{i}$ and $\lim_{n\to\infty}p_{n,i}=p_{i}$ for each $i$ then $W_{r}(G_{n},G)\to0$.
\end{lemma}

\begin{proof}
We use the fact that $W_{r}(G_{n},G)\to0$ is equivalent to weak convergence plus convergence of the first $r$ moments. For weak convergence, recall that $\int f\,dG_{n}= \sum_{i=1}^{k}p_{n,i}f(\theta_{n,i})$ and hence for any bounded continuous $f$,
\begin{align*}
\int f\,dG_{n}
= \sum_{i=1}^{m}p_{n,i}f(\theta_{n,i})
\longrightarrow \sum_{i=1}^{m}p_{i}f(\theta_{i})
= \int f\,dG.
\end{align*}

\noindent
Thus $G_{n}$ converges weakly to $G$. Furthermore, for any $r\ge 1$,
\begin{align*}
\int \basem(x,x_{0})^{r}\,dG_{n}(x)
= \sum_{i=1}^{m}p_{n,i}\basem(\theta_{n,i},x_{0})^{r}
\longrightarrow \sum_{i=1}^{m}p_{i}\basem(\theta_{i},x_{0})^{r}
= \int \basem(x,x_{0})^{r}\,dG(x).
\end{align*}

\noindent
Thus the first $r$ moments of $G_{n}$ converge to those of $G$. It follows that $W_{r}(G_{n},G)\to0$.
\end{proof}

\section{Wasserstein consistency of the MHDE}
\label{app:mdhe}

Assume $L$ is fixed; any dependence on $L$ will be suppressed here. Let $\{q_{\varphi}:\varphi\in\Phi\}$ be a parametric family of densities such that $\Phi$ is compact, and define
\begin{align*}
\Theta
=\Bigg\{
\theta = (p_{1},\ldots,p_{L},\varphi_{1},\ldots,\varphi_{L})
\in [0,1]^{L}\times \Phi^{L}
: \sum_{\ell=1}^{L}p_{\ell}=1
\Bigg\}.
\end{align*}

\noindent
For any $\theta\in\Theta$, let $Q_{\theta}$ denote the mixture distribution defined by $\theta$, and $\projset$ the family of mixing measures induced by $\Theta$.

Suppose that $\proj=Q_{\theta^{*}}=\sum_{\ell=1}^{L}\projwgt_{\ell}\projcmp_{\ell}$ is the Hellinger projection of $\mixfcn(\truemix)$ onto $\mix(\projset)$ and $\projmix=\mixingmap(\proj)\in\projset$, i.e.
\begin{align*}
\probm(\proj,\mixfcn(\truemix))
< \probm(Q_{\theta}, \mixfcn(\truemix)) 
\quad\forall\,\theta\ne\theta^{*}.
\end{align*}

\noindent
Given $\rv^{(1)},\ldots,\rv^{(n)}\iid\mixfcn(\truemix)$, let $\wh{\true}(\rv^{(1)},\ldots,\rv^{(n)})$ be a suitably chosen kernel density estimate of $\mixfcn(\truemix)$, and define $\wh{\theta}$ by 
\begin{align*}
\probm(Q_{\wh{\theta}},\wh{\true}(\rv^{(1)},\ldots,\rv^{(n)}))
\le \probm(Q_{\theta}, \wh{\true}(\rv^{(1)},\ldots,\rv^{(n)})) 
\quad\forall\,\theta\in\Theta.
\end{align*}

\noindent
This is the minimum Hellinger distance estimator (MHDE) defined by \citet{beran1977}. Then by the results of \citet{beran1977} \citep[see also][]{basu2011}, we conclude that $\wh{\theta}\to\theta^{*}$. Assuming that $\varphi_{n}\to\varphi$ implies $\probm(q_{\varphi_{n}}, q_{\varphi})\to 0$ (this is true, for example, when $q_{\varphi}$ is Gaussian), we deduce that $\probm(\estcmp_{\varphi_{\ell}},\projcmp_{\varphi_{\ell}})\to0$ (possibly up to re-labeling) and $\wh{p}_{\ell}\to p^{*}_{\ell}$. But this implies that $\mixm(\estmix,\projmix)\to0$, where $\estmix$ is the mixing measure induced by $\wh{\theta}$.

\section{Experiment details}
\label{app:expts}

For each experiment, we used the same simulation procedure:
\begin{enumerate}
\item Generate $n$ samples from the model (see below).
\item Use the EM algorithm with weight clipping to approximate a Gaussian mixture model with $L=100$ components. We used 20 random initializations and picked the estimate with the highest log-likelihood and terminated the algorithm at 1000 iterations (if convergence had not already occurred). Call the result $\wh{Q}=\sum_{\ell=1}^{L}\wh{\gamma}_{\ell}\wh{q}_{\ell}$.
\item Use the 2-Wasserstein distance to compute the distance matrix $D(\wh{Q})=(\probm(\wh{q}_{\ell},\wh{q}_{m}))_{\ell,m=1}^{L}$. For gaussian measures $\mu_{i}\sim\normalN(m_{i},V_{i})$, there is a closed form expression for the 2-Wasserstein distance:
\begin{align*}
\probm(\mu_{i},\mu_{j})
:= W_{2}(\mu_{i}, \mu_{j})
= \norm{m_{i}-m_{j}}_{2}^{2} + \tr(V_{i}+V_{j}-2(V_{j}^{1/2}V_{i}V_{j}^{1/2})^{1/2}).
\end{align*}
\item Use single-linkage hierarchical clustering to cluster the $\{\wh{q}_{\ell}\}$ into $K$ clusters using the distance matrix $D(\wh{Q})$, where $K$ is given by the model.
\end{enumerate}

\noindent
The details of each model are as follows (unless otherwise noted, $n=5000$ samples were drawn for the model):
\begin{itemize}
\item \textsc{GaussGamma} $(K=4)$. The data is generated from the following closed-form density:
\begin{align*}
\true &= \sum_{k=1}^{4}\truewgt_{k}\truecmp_{k},
\quad \text{where}\quad
\left\{
\begin{aligned}
\truecmp_{1} &\propto 0.22 \cdot \normalN(-7, 0.3^2) + 0.08 \cdot\normalN(-6, 0.2^2) \\
&\qquad\quad+ 0.15 \cdot \normalN(-5, 0.5^2), \\
\truecmp_{2} &\propto 0.15 \cdot \normalN(-1, 0.3^2), \\ 
\truecmp_{3} &\propto 0.15 \cdot \normalN(5, 0.3^2), \\
\truecmp_{4} &\propto 0.3 \cdot \GammaDist(18,0.5).
\end{aligned}
\right.
\end{align*}
\item \textsc{Gumbel} $(K=3)$. First, data is generated from the following mixture of Gaussians:
\begin{align*}
\mu &= \sum_{k=1}^{3}\truewgt_{k}\truecmp_{k},
\quad \text{where}\quad
\left\{
\begin{aligned}
\truecmp_{1} &\propto 0.22 \cdot \normalN(-7, 0.3^2) + 0.08 \cdot\normalN(-6, 0.2^2) \\
&\qquad\quad+ 0.15 \cdot \normalN(-5, 0.5^2), \\
\truecmp_{2} &\propto 0.15 \cdot \normalN(-1, 0.3^2), \\ 
\truecmp_{3} &\propto 0.15 \cdot \normalN(5, 0.3^2).
\end{aligned}
\right.
\end{align*}

\noindent
Note that this is the same as the previous model modulo the Gamma term. Given $Y^{(i)}\sim \mu$, $i=1,\ldots,n$, we then contaminate each sample with Gumbel noise $W\sim\GumbelDist(0,0.3)$. Thus, the final data are
\begin{align*}
Z^{(i)}
= Y^{(i)} + W^{(i)}, \quad 
Y^{(i)}\iid\mu, \quad
W^{(i)}\iid\GumbelDist(0,0.3).
\end{align*}
\item \textsc{Poly} $(K=2)$. The data is generated from a density $f$ defined as follows: Define two polynomials by
\begin{align*}
p_{1}(x) &= -2x^6 - 5x^5 - 2x^3 - 10x^2 + 7 \\
p_{2}(x) &= -x^4 + 3.5x^3 - 3x^2 + 2.
\end{align*}

\noindent
Let $p^{+}_{i}(x) = \max(0,p_{i}(x))$. Define
\begin{align*}
f(x)
&= \alpha p_{1}^{+}(x+c) + \beta p_{2}^{+}(x), \\
\alpha &= \frac12\left(\int_{-\infty}^{\infty} p_{1}^{+}(x+c)\,dx\right)^{-1}, \\
\beta &= \frac12\left(\int_{-\infty}^{\infty} p_{2}^{+}(x)\,dx\right)^{-1}.
\end{align*}

\noindent
The constants $\alpha$ and $\beta$ are defined here so that $\int f(x)\,dx=1$. The constant $c>0$ is chosen so as to separate the the support of each component by 1. In our experiments, this resulted in the values $c = -2.317$, $\alpha = 0.026$ and $\beta = 0.101$.
\item \textsc{Moons} $(K=2)$.
First, $n = n_1+n_2$ random points are randomly generated from the unit circle in $\R^{2}$, $n_1$ of these points are from the upper half of the circle (i.e. with positive $y$ value), and $n_2$ of these points are from the lower half of the circle (i.e. with negative $y$ value). Then, the samples in the upper half of the circle are moved to the left, and the lower half to the right, respectively, by a distance $b>0$. Each sample is then perturbed by Gaussian noise $W\sim\normalN(0,rI_{2\times 2})$. In our simulations we used $b=0.5$ and $r=0.015$. In the balanced case, $n_1 = n_2 = 2500$; in the unbalanced case, $n_1= 3000, n_2 = 500$.

\item \textsc{Sobolev} $(K=3)$. Random functions $g_{k}$ ($k=1,2,3$) were generated from a random expansion of the orthogonal basis of $H^{1}(\R)$ given by $\beta_{j}(x)=b_{j}\sin(x/b_{j})$ with $b_{j}=2/((2k-1)\pi)$. After appropriate normalization and truncation, the nonparametric mixture components are given by $f_{k}(x):=\exp(-g_{k}(x))/Z_{k}$, where $Z_{k}$ is a normalization constant. The weights were set to $\truewgt_{1}=\truewgt_{2}=\truewgt_{3}=1/3$.

\item \textsc{Target} $(K=6)$. The mixture model is generated as follows:

\begin{enumerate}
\item Subsample $k_1 = 15$ points from the middle cluster (denoted by $C_1$) of the original Target dataset,\footnote{\url{https://www.uni-marburg.de/fb12/arbeitsgruppen/datenbionik/data?language_sync=1}} with $k_2 = 80 $ points from the cyclic cluster (denoted by $C_2$) and $k_3 = 12$ points from outliers in four corners (denoted by $C_3$). Call these points $m_i$.
\item Define $\true = \sum_{k=1}^3 \truewgt_k \truecmp_k$, where
\begin{align*}
\truecmp_1 &\propto \sum_{m_i \in C_1} N(m_i, V_i), 
\quad 
V_i = 0.04 I_2;\\
\truecmp_2 &\propto \sum_{m_i \in C_2} N(m_i, V_i), 
\quad
V_i = 0.15 I_2;\\
\truecmp_3 &\propto \sum_{m_i \in C_3} N(m_i, V_i), 
\quad
V_i = 0.15 I_2;\\
\truewgt_1 &= 0.3, 
\quad
\truewgt_2 = 0.4, 
\quad
\truewgt_3 = 0.3.
\end{align*}

\noindent
Note that $\truecmp_{3}$ defines 4 separate components (one for each corner) with a different number of Gaussian components. Thus, the final mixture model $\true$ has $K=6$ components.

\item Generate $n = 2000$ samples from $\true$.

\end{enumerate}
\end{itemize}

\end{document}